\documentclass[12pt]{article}
\usepackage[normalem]{ulem}
\usepackage{mathpazo}
\usepackage{authblk}
\usepackage{hyperref}
\usepackage{tikz-cd}
\usepackage{tikz}
\usepackage[all]{xy}
\usepackage{lineno}
\usepackage{xcolor}
\usepackage{nomencl}
\makenomenclature
\renewcommand{\nompreamble}{The list below describes several symbols that will be used throughout the article.}

 \usepackage[utf8]{inputenc}
 \usepackage{scalerel}
 \usepackage{enumerate}
 \usepackage{authblk}
 \usepackage{wrapfig}

 \usepackage{amsmath,amscd}
 \usepackage{amsthm}
 \usepackage{amsfonts}
 \usepackage{color}
 
\usepackage{thmtools}
\usepackage{thm-restate}

\usepackage{float}

 \usepackage{amssymb}
 \usepackage[margin=1in]{geometry}

 \newtheorem{theorem}{Theorem}
 \newtheorem{proposition}{Proposition}[section]
 \newtheorem{corollary}[proposition]{Corollary}
 \newtheorem{lemma}[proposition]{Lemma}
 \newtheorem{remark}[proposition]{Remark}

 \theoremstyle{definition}
 \newtheorem{definition}[proposition]{Definition}

\newcommand{\R}{
	\mathbb{R}
}

\newcommand{\diam}{\mathrm{diam}}
\newcommand{\dgh}{d_\mathrm{GH}}

\newcommand{\dis}{\mathrm{dis}}
\newcommand{\distance}{\mathrm{distance}}

\newcommand{\betti}{b_1}
\newcommand{\hyp}{\mathrm{hyp}}

\newcommand{\length}{\mathrm{L}}

\newcommand{\Reeb}[1]{\mathrm{R}_{#1}}
    \newcommand{\Rf}{\Reeb{f}}
    \newcommand{\Rg}{\Reeb{g}}
    \newcommand{\Rh}{\Reeb{h}}
    \newcommand{\Rp}{\Reeb{p}}
    \newcommand{\Rq}{\Reeb{q}}
    \newcommand{\Rphi}{\Reeb{\phi}}
    \newcommand{\Rpsi}{\Reeb{\psi}}
\newcommand{\barf}{\bar{f}}
\newcommand{\barg}{\bar{g}}
\newcommand{\pf}{\pi_f}

\newcommand{\pFunc}[3]{\pi_{#1}^{{#2},{#3}}}
    \newcommand{\pfrs}{\pFunc{f}{r}{s}}
\newcommand{\calB}{\mathcal{B}}
\newcommand{\Reebs}[2]{\mathrm{R}_{#1}^{#2}}
    \newcommand{\Rfr}{\Reebs{f}{r}}
    \newcommand{\Rfs}{\Reebs{f}{s}}
    \newcommand{\Rgr}{\Reebs{g}{r}}

\newcommand{\lr}{\mathrm{length}_\R}

\newcommand{\DistF}[1]{\mathrm{d}_{#1}}
    \newcommand{\df}{\DistF{f}}
    \newcommand{\dg}{\DistF{g}}
\newcommand{\DistFr}[2]{\mathrm{d}_{#1}^{#2}}
    \newcommand{\dfr}{\DistFr{f}{r}}

\usepackage{chngcntr}
\usepackage{apptools}
\AtAppendix{\counterwithin{theorem}{section}}
\definecolor{darkblue}{rgb}{0.0, 0.0, 0.8}
\definecolor{darkred}{rgb}{0.8, 0.0, 0.0}
\definecolor{darkgreen}{rgb}{0.0, 0.8, 0.0}

\providecommand{\keywords}[1]
{
  \small	
  \textbf{\textbf{Keywords---}} #1
}

\begin{document}

\title{Metric graph approximations of geodesic spaces\thanks{We acknowledge funding from these sources: NSF AF 1526513, NSF DMS 1723003, NSF CCF 1740761, NSF CCF 2310412, and BSF 2020124.}}

\author[1]{Facundo M\'emoli}
\author[2]{Osman Berat Okutan}
\author[3]{Qingsong Wang}

\affil[1]{Department of Mathematics,
    The Ohio State University.\\

    \texttt{facundo.memoli@gmail.com}}
\affil[2]{Max Planck Institute for Mathematics in the Sciences\\

    \texttt{osman.okutan@mis.mpg.de}}
\affil[3]{Department of Mathematics,
University of Utah\\

    \texttt{qswang@math.utah.edu}}

\maketitle

\begin{abstract}
    We study the question of approximating a compact geodesic metric space by metric graphs satisfying a uniform upper bound on  their first Betti number. We prove that, up to a suitable multiplicative constant, Reeb graphs of distance functions to a point provide optimal approximation in the Gromov-Hausdsorff sense. 
\end{abstract}

\keywords{Metric graphs, Reeb graphs, Gromov-Hausdorff distance}

\newpage
\tableofcontents
\newpage

\section{Introduction}\label{sec:intro}

How can one faithfully represent a shape by a graph? 
\begin{figure}
	\centering
	\includegraphics[width=0.5\linewidth]{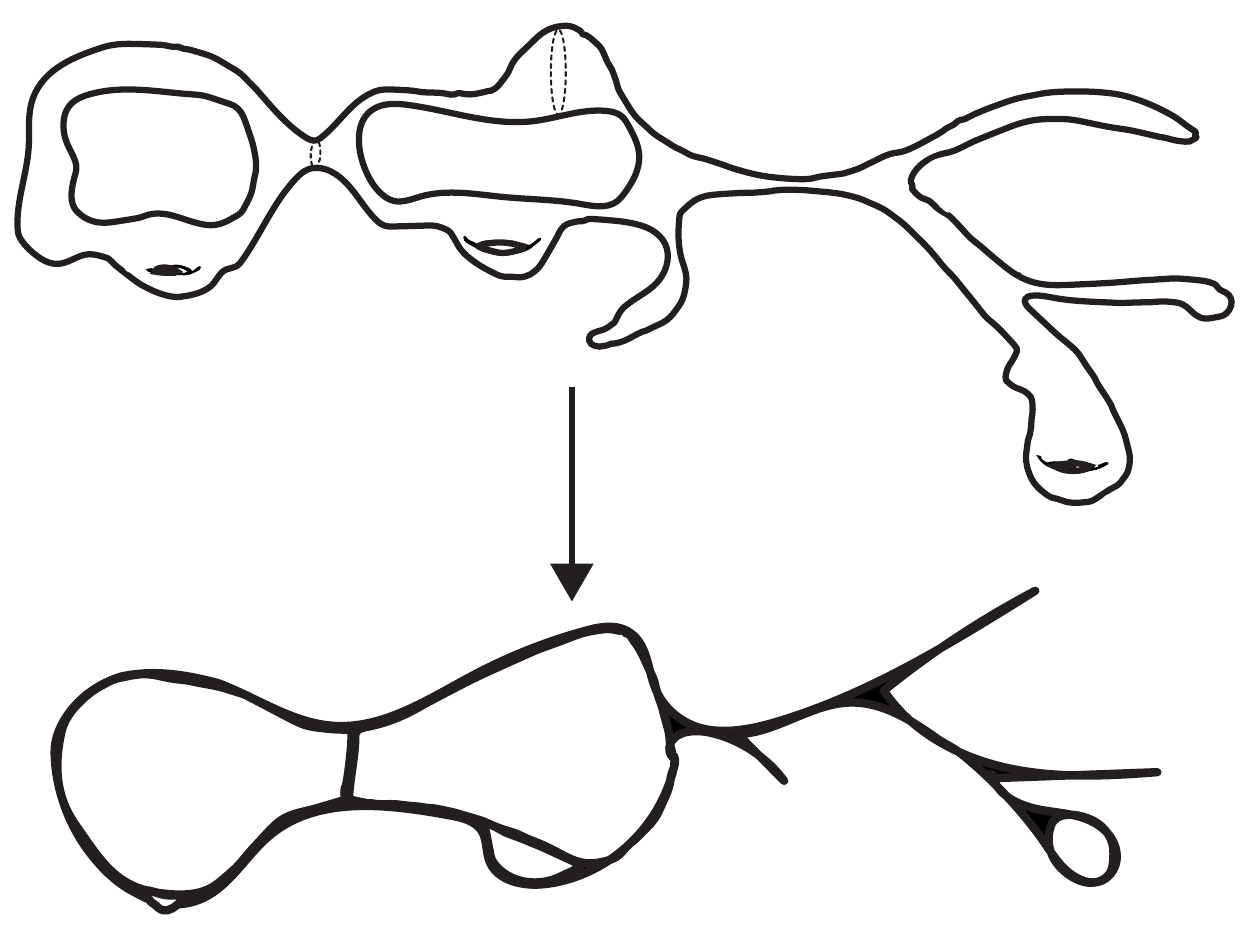}
	\caption{Approximating a  geodesic space by a metric graph.}
\end{figure}
\begin{figure}[h]
	\centering
\includegraphics[width=0.4\linewidth]{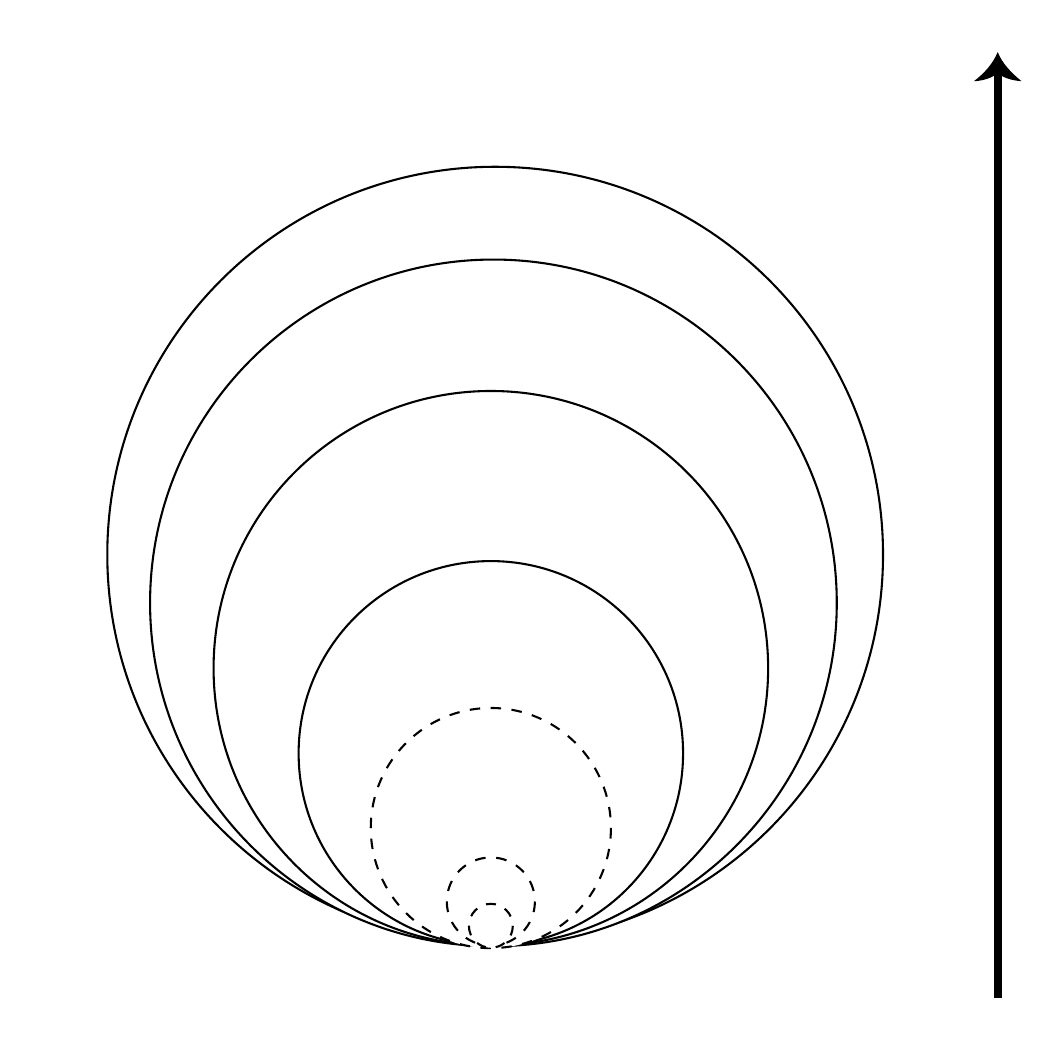}
	\caption{The Hawaiian earring is a compact geodesic space. The Reeb graph of the height function depicted yields an infinite graph.}
	\label{fig:hawaii}
\end{figure}
Assuming that the original shape $X$ is a geodesic space, if we sample  points from it densely enough  and connect sufficiently close pairs of sample points, then the resulting graph $G$ will be a good approximation of the original shape $X$ in Gromov-Hausdorff sense. However, the  graph $G$ tends to be very complicated; for example, it can have a very high first Betti number $\betti(G)$. What if we put an upper bound on the first Betti number? In this case, \emph{Reeb graphs} of suitable functions $f:X\to \R$ come to mind, as their first Betti numbers are always bounded above by that of the original shape. There is a potential problem with this approach: Reeb graphs are known to be finite graphs whenever the given function is constructible in a certain sense, for example Morse or piecewise linear. However, in general, a geodesic space does not necessarily carry a simplicial or smooth structure and, furthermore, there are cases for which the resulting Reeb graph is not finite; see Figure \ref{fig:hawaii}.

 Still, a geodesic space (just as any metric space), provides one with a rich class of natural real-valued functions reflecting the geometric structure of the space, namely  distance functions to individual points, and other functions built using convex combinations and/or maxima-minima of such functions. All those functions are 1-Lipschitz. Now, the question is, if we start with a 1-Lipschitz function on a geodesic space, is the quotient space induced through the Reeb construction really a graph? Or: in what sense can it be regarded as a graph? The Hawaiian ring example shows that this space is not always a finite graph.

\begin{figure}[h]
	\centering \includegraphics[width=0.7\linewidth]{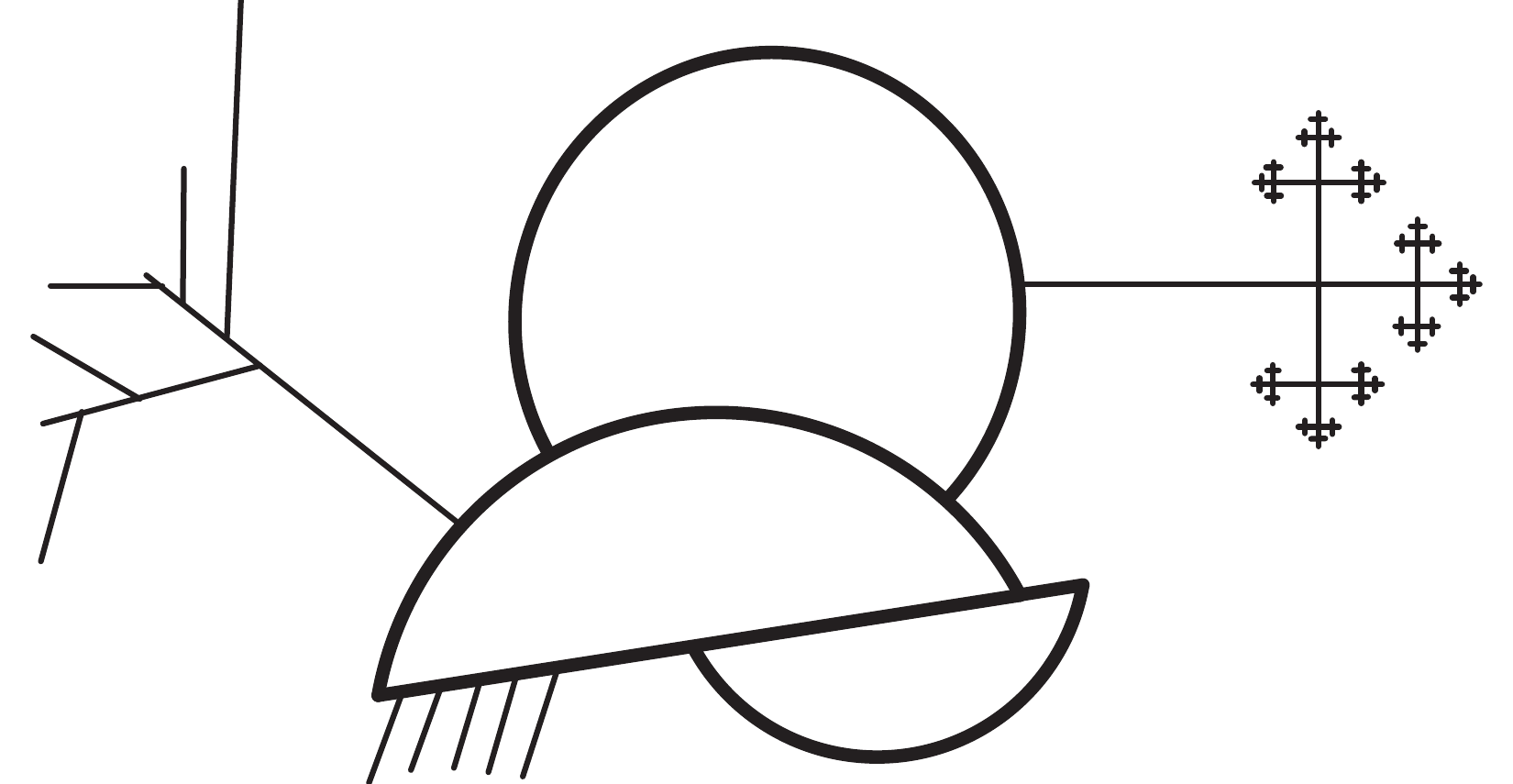}
	\caption{An illustration of a ``tame" metric graph: it is allowed to exhibit ``thorny" metric tree-like parts while having a finite metric graph as a ``core".}
	\label{fig:tame-graph}
\end{figure}

In this paper, we show that given a compact geodesic space $X$ with $\betti(X)<\infty$ and a 1-Lipschitz function $f: X \to \R$, the Reeb graph $\Rf$ endowed with the intrinsic metric induced by length structure on $\R$ pulled back through $f$, is a ``tame" metric graph in the following  sense: It has a geodesic subspace $G$ (the ``core") which is isometric to a finite graph such that the quotient $X/G$ is a (possibly infinite) tree; see Figure \ref{fig:tame-graph}. We call $\Rf$ with the metric described above a \emph{Reeb metric graph}. As in the Morse and piecewise linear case, $\betti(\Rf) \leq \betti(X)$. We also show that, over geodesic spaces endowed with 1-Lipschitz functions, the Reeb metric graph construction is Gromov-Hausdorff stable, where the constant of stability depends on the first Betti numbers of the geodesic spaces and the number of local minima of the functions (cf. Theorem \ref{thm:stability}).

Now, after establishing that 1-Lipschitz functions produce metric graphs (in the sense above), we can ask how faithfully do they represent the original space.  Let $\delta_n(X)$ denote the infimal Gromov-Hausdorff distance to $X$ that can be achieved by a finite metric graph with the first Betti number less than or equal to $n$. Let $\rho_{p}(X):=\dgh(X,\Rf)$ where $p \in X$, and $f: X \to \R$, $ x \mapsto d_X(p,x)$. In Theorem \ref{thm:approximation} we show that $\rho_{p}(X) \leq C(n,X) \, \delta_n(X)$ where $C(n,X)$ depends linearly on $n$ and $\betti(X)$. This means that Reeb metric graphs of distance to a point functions are optimal metric graph representations up to a multiplicative constant depending linearly on the upper bound on the first Betti numbers of the approximating metric graphs. For compact geodesic spaces with Hausdorff dimension 2, we provide an upper bound to $\rho_p(X)$ in terms of the area of the space $X$.

We investigate tree approximations separately. We define merge metric trees induced by 1-Lipschitz $f: X \to \R$, denoted by $T_f$. In Theorem \ref{thm:metric_tree_approximation} we show that the construction produces a metric tree in a Gromov-Hausdorff stable manner. We show that if $f: X \to \R$ is distance to a point, then $\dgh(X,T_f) \leq 13\, \delta_0(X)$. This means that the merge metric tree induced by distance to a point is the best metric tree approximation up to a multiplicative constant.

\section{Topological Reeb Graphs}\label{sec:reeb}

In this section $X, Y$ denote compact Hausdorff topological spaces, and $f: X \to \R$, $g: Y \to \R$ denote continuous functions.

\begin{definition}[Topological Reeb graph]
    The \emph{topological Reeb graph}, or Reeb graph in short,
    associated with $f: X \to \R$, denoted by $\Rf$, is the quotient space of $X$ that identifies points in the same \emph{connected component} of the fibers $f^{-1}(c)$ for all $c \in \R$. The quotient map from $X$ to $\Rf$, denoted $\pf: X \to R_f$, is called the \emph{Reeb quotient map}. We alternatively denote $\pf(x)$ by $[x]$. Note that $f$ induces a continuous function on $\Rf$, which we denote by $\barf: \Rf \to \R$ and define as $\barf([x])=f(x)$.
\end{definition}

The following lemma provides an alternative characterization of the equivalence class induced by the Reeb quotient map.

\begin{lemma}\label{lem:alternative_char_quotient}
    Let $X$ be a compact, Hausdorff topological space and $f: X \rightarrow \mathbb{R}$ be a continuous function. For a given $x \in X$, we introduce the following notation related to the preimage of $f$:
    \begin{enumerate}
        \item $C(x)$ represents the connected component of $x$ in $f^{-1}(f(x))$.
        \item For any $\epsilon > 0$, $C_\epsilon(x)$ is the connected component of $x$ in the preimage of the closed interval $f^{-1}([f(x)-\epsilon, f(x)+\epsilon])$.
        \item $\mathring{C}_\epsilon(x)$ stands for the connected component of $x$ in the preimage of the open interval $f^{-1}((f(x)-\epsilon, f(x)+\epsilon))$. Note that, if $X$ is locally path connected, then $\mathring{C}_\epsilon(x)$ is a open connected subset in $X$ and hence also path connected.
        \item $PC_\epsilon(x)$ denotes the path component of $x$ in the preimage of the closed interval $f^{-1}([f(x)-\epsilon, f(x)+\epsilon])$.
    \end{enumerate}

    Then, the following properties hold:

    \begin{enumerate}
        \item $C(x) = \bigcap_{\epsilon > 0} C_\epsilon(x)$.
        \item $C(x) = \bigcap_{\epsilon > 0} \mathring{C}_\epsilon(x)$.
        \item if, additionally, $X$ is locally path connected, then $C(x) = \bigcap_{\epsilon > 0} PC_\epsilon(x)$.
    \end{enumerate}
\end{lemma}

\begin{proof}
    We first show $C(x) = \bigcap_{\epsilon > 0} C_\epsilon(x)$.
    By definition, $C(x) \subseteq C_\epsilon(x)$ for all $\epsilon > 0$, so $C(x) \subseteq \bigcap_{\epsilon > 0} C_\epsilon(x)$.
    Conversely, it is clear that $\bigcap_{\epsilon > 0} C_\epsilon(x) \subseteq f^{-1}(f(x))$.
    Note that $C_\epsilon(x)$ is closed in $X$ for all $\epsilon > 0$ and hence $C_\epsilon(x)$ is compact for all $\epsilon > 0$.
    Therefore, $C_\epsilon(x)$ is a decreasing sequence of compact connected sets and hence $\bigcap_{\epsilon > 0} C_\epsilon(x)$ is connected by \cite[Corollary~6.1.19]{engelking1989general}.
    As a connected subset of $f^{-1}(f(x))$ that contains the connected component $C(x)$, $\bigcap_{\epsilon > 0} C_\epsilon(x)$ must be equal to $C(x)$.

    For the second item, it is clear from the inclusions
    \[
        \mathring{C}_\epsilon(x) \subseteq C_\epsilon(x) \subseteq \mathring{C}_{2\epsilon}(x)
    \]
    that $\bigcap_{\epsilon > 0} \mathring{C}_\epsilon(x) = \bigcap_{\epsilon > 0} C_\epsilon(x)$.
    Then, the result follows from the first item.

    For the third item, as each $\mathring{C}_\epsilon$ is also the path component of $x$ in $f^{-1}((f(x)-\epsilon, f(x)+\epsilon))$, we have the following inclusion maps
    \[
        \mathring{C}_\epsilon(x) \subseteq PC_{\epsilon}(x) \subseteq \mathring{C}_{2\epsilon}(x).
    \]
    Then we have $\bigcap_{\epsilon > 0} PC_\epsilon(x) = \bigcap_{\epsilon > 0} \mathring{C}_\epsilon(x)$, and the result follows from the second item.
\end{proof}

As a corollary, we show that the Reeb graph is compact Hausdorff whenever $X$ is compact Haudorff.

\begin{corollary}\label{cor:reeb_hausdorff}
    Let $X$ be a compact, Hausdorff topological space and $f: X \rightarrow \mathbb{R}$ be a continuous function. Then $\Rf$ is a compact, Hausdorff topological space.
\end{corollary}

\begin{proof}
    As $X$ is compact, its quotient space $\Rf$ is also compact.
    To show that $\Rf$ is Hausdorff, it suffices to show that for any $[x], [x'] \in \Rf$ such that $[x] \neq [x']$, there exist disjoint open sets $U, V$ such that $[x] \in U$ and $[x'] \in V$.
    Let $f(x)=c$ and $f(x')=c'$.
    If $c \neq c'$, then we can find disjoint open sets $U, V$ of the real line such that $c \in U$ and $c' \in V$.
    Then $\pf^{-1}(U \cap f(X))$ and $\pf^{-1}(V \cap f(X))$ of $X$ such that $[x] \in \pf^{-1}(U \cap f(X))$ and $[x'] \in \pf^{-1}(V \cap f(X))$.
    We now turn attention to the case where $c = c'$.
    Let $\mathring{C}_\epsilon(x)$ and $\mathring{C}_\epsilon(x')$ be the connected components in $f^{-1}((c-\epsilon, c+\epsilon))$ containing $x$ and $x'$ respectively. Then for each $\epsilon>0$ either $\mathring{C}_\epsilon(x) \cap \mathring{C}_\epsilon(x') = \emptyset$ or $\mathring{C}_\epsilon(x) = \mathring{C}_\epsilon(x')$.
    We claim that there exists some $\epsilon > 0$ such that $\mathring{C}_\epsilon(x) \cap \mathring{C}_\epsilon(x') = \emptyset$ and hence concludes the proof.
    If not, then for any $\epsilon > 0$, $\mathring{C}_\epsilon(x) = \mathring{C}_\epsilon(x')$.
    In this case, by Lemma~\ref{lem:alternative_char_quotient}, $[x]=[x']$, which is a contradiction.
\end{proof}

\begin{corollary}\label{cor:reeb_path_component}
    Let $X$ be a locally path connected, compact, Hausdorff topological space and $f: X \rightarrow \mathbb{R}$ be a continuous function. Let $U$ be a path component of $f^{-1}(I)$ for an open interval $I \subseteq \R$.
    Then we have that $U = \pf^{-1}(\pf(U))$, $\pf(U)$ is open and $\pf(U)$ is a path component of $\barf^{-1}(I)$.
\end{corollary}

\begin{proof}
    To show that $U = \pf^{-1}(\pf(U))$, it suffices to show that $U$ is closed under the equivalence relation defining $\Rf$. Let $x, x'\in U$ be such that $[x]=[x']$. Then $f(x)=f(x')=c$ for some $c\in \R$. Then by Lemma~\ref{lem:alternative_char_quotient}, $x,x'$ are in the same path component of $f^{-1}([c-\epsilon,c+\epsilon])$ for all $\epsilon>0$.
    Note that, there exists some $\epsilon>0$ such that $[c-\epsilon,c+\epsilon]\subseteq I$ and hence $f^{-1}([c-\epsilon,c+\epsilon])\subseteq f^{-1}(I)$.
    Since $y$ belongs to the same path component of $f^{-1}([c-\epsilon,c+\epsilon])$ as $x$, $x'$ also belongs to $U$. This shows that $U$ is closed under the equivalence relation defining $\Rf$.

    Since $\pf$ is a quotient map, a set $V$ in $\Rf$ is open if and only if $\pf^{-1}(V)$ is open in $X$. Since $\pf^{-1}(\pf(U))$ is open in $X$, $\pf(U)$ is open in $\Rf$.

    To show the last claim, let $f^{-1}(I)=\cup_i U_i$ be the decomposition of $f^{-1}(I)$ into the disjoint union of its path components. Then $\cup_i \pf(U_i)$ is a union of connected sets and $\cup_i \pf(U_i) = \barf^{-1}(I)$. Since $U_i$ is path-connected, $\pf(U_i)$ is connected. Hence, it suffices to show that $\cup_i \pf(U_i)$ is a disjoint union. This comes from the fact that $U_i$ are closed under the equivalence relation defining the quotient map $\pf$.
\end{proof}

The following proposition describes some topological properties of Reeb graphs.

\begin{proposition}\label{prop:hausdorff}
    Let $X$ be a locally path connected, compact, Hausdorff topological space and let $f: X \rightarrow \mathbb{R}$ be a continuous function. $R_f$ is a locally path connected, compact, Hausdorff space with a topological basis given by $$\calB :=\{\pf(U): U \text{ is a path component of } f^{-1}(I) \text{ for an open interval } I \subseteq \R \}.$$
\end{proposition}

\begin{proof}
    We first show that $\calB$ is a topological basis for $\Rf$.
    By Corollary~\ref{cor:reeb_path_component}, we know that $\calB$ consists of open sets. Then it suffices to show that for any $[x] \in \Rf$ and any neighborhood $V$ of $[x]$ in $\Rf$, there exists $U \in \calB$ such that $[x] \in U \subseteq V$.
    Let $f(x)=c$ and $C(x)$ be the connected component of $x$ in $f^{-1}(c)$.
    As in Lemma~\ref{lem:alternative_char_quotient}, we use $C_\epsilon(x)$ to denote the connected component of $x$ in $f^{-1}([c-\epsilon,c+\epsilon])$ for $\epsilon>0$.
    Lemma \ref{lem:alternative_char_quotient} shows that
$\cap_{\epsilon>0} C_\epsilon(x) = C(x) \subseteq \pf^{-1}(V)$.
    We use the notation $C_\epsilon^c(x)$ to denote the complement of $C_\epsilon(x)$ in $X$.
        Since $\cap_{\epsilon>0} C_\epsilon(x) = C(x) \subseteq \pf^{-1}(V)$, we have that, for any fixed $x\in X$, the sets
        $\{ (C_\epsilon^c(x))_\epsilon, \pf^{-1}(V)\}$ form an open cover of $X$.
    
    Since $X$ is compact, $\{(C_\epsilon^c(x))_\epsilon, \pf^{-1}(V)\}$ admits a finite subcover.
    Note that $C_\epsilon^c(x)$ is a decreasing family of open sets as $\epsilon \to 0$, 
        as, for any $0< \epsilon < \epsilon'$, the inclusion $C_{\epsilon}^c(x) \subseteq C_{\epsilon'}^c(x)$ holds, we have that
     there exists $\epsilon_0>0$ such that $C_{\epsilon_0}^c(x) \cup \pf^{-1}(V)=X$. This implies that $C_{\epsilon_0}(x) \subseteq \pf^{-1}(V)$. Now if we let $U$ to be the path component of $x$ in $f^{-1}(c-\epsilon_0,c+\epsilon_0)$, then $U \subseteq C_{\epsilon_0}(x)$, and $\pf(U) \in \calB$ and $\pf(U) \subseteq V$. This shows that $\calB$ is a basis for $\Rf$.

    According to Corollary~\ref{cor:reeb_path_component}, each $U \in \calB$ is path connected and hence $\Rf$ is locally path connected.
    Since $X$ is compact, $\Rf$ is also compact.
        The fact that $\Rf$ is Hausdorff follows from Corollary~\ref{cor:reeb_hausdorff}.
    
\end{proof}

The following theorem describes how maps between $X$ and $Y$ induce maps between the Reeb graphs associated with $f: X\to \R$ and $g: Y\to \R$. The proof utilizes some basic results about quotient topologies and quotient maps that can be found, for example, Section 2.4 of \cite{engelking1989general}

\begin{theorem}[Naturality]\label{thm:naturality}
    Let $f: X \to \R$, $g: Y \to \R$ and $\psi: X \to Y$ be a continuous maps such that $f= g \circ \psi$. Then $\bar{\psi}: \Rf \to \Rg$, $[x] \mapsto [\psi(x)]$ is a well-defined continuous map. Furthermore, if $\psi$ is onto and has connected fibers, then $\bar{\psi}$ is a homeomorphism.
\end{theorem}

\begin{proof}
    Since $f = g\circ \psi$, we have that for any $x, x'$ in the fiber $f^{-1}(c), c\in \R$, $\psi(x), \psi(x')$ belong to the fiber $g^{-1}(c)$.
    As $\psi$ is continuous, it maps connected sets to connected sets, and hence we have $[\psi(x)] = [\psi(x')]$. Then by the universal property of quotient maps, the map $\bar{\psi}$ is well-defined and is the unique continuous map that satisfies $\psi \circ \pi_g = \bar{\psi}\circ \pi_f$.

    For the second claim, first note that as $\psi$ is a surjective map between compact Hausdorff spaces, it is closed and is a quotient map. Then the composition $\pi_g \circ \psi$ is also a quotient map. Since $\pi_g \circ \psi = \bar{\psi} \circ \pi_f$ and both $\pi_g \circ \psi$ and $\pi_f$ are quotient maps, $\bar{\psi}$ is a quotient map. Then to show that $\bar{\psi}$ is a homeomorphism, it suffices to show that $\bar{\psi}$ is injective.
    That is, for any $x, x'\in X$ such that $[\psi(x)] = [\psi(x')]$, we need to show that $[x] = [x']$. Let $C$ be the connected component in $g^{-1}(f(x))$ containing $\psi(x)$ and $\psi(x')$. Then $\psi^{-1}(C)$ belongs to the fiber $f^{-1}(x)$ and contains $x$ and $x'$. As $\psi$ is a quotient map that is closed and with connected fibers, \cite[Theorem~6.1.29]{engelking1989general} guarantees that $\psi^{-1}(C)$ is connected. Then $[x] = [x']$.
\end{proof}

\begin{corollary}\label{cor:idempotencyAndProduct}
    Let $f: X \to \R$, $g: Y \to \R$ be continuous maps  inducing the Reeb graphs $\Rf$ and $\Rg$ respectively. Then, the following properties hold:    \begin{enumerate}[i)]
        \item (Idempotency) $\Rf$ is homeomorphic to $\Reeb{\barf}$.
        \item (Product) Let $h_1: X \times Y \to \R$, $(x,y) \mapsto f(x)+g(y)$ and $h_2 : \Rf \times Y \to \R$, $([x],y) \mapsto f(x)+g(y)$. Then the map $\Reeb{h_1} \to \Reeb{h_2}$, $[x,y] \mapsto [[x],y]$ is a homeomorphism.
    \end{enumerate}
\end{corollary}
\begin{proof}
    We consider the quotient map $\pi_f: X \to \Rf$. Then, it holds that $\barf \circ \pi_f = f$ where $\barf: \Rf \to \R$ is the map induced by $f$ on $\Rf$.
    . Addditonly, as $X$ is compact Hausdorff, so is $\Rf$ by Corollary \ref{cor:reeb_hausdorff}. Note that the map $\pi_f: X \to \Rf$ is onto and has connected fibers. Therefore, by Theorem \ref{thm:naturality} $\barf: \Rf \to \Rf$ is a homeomorphism. This proves $i)$.

    For the second item, we consider the map $\pi_{f}\times \mathrm{id}:X \times Y \to \Rf \times Y$ that sends $(x,y)$ to $([x],y)$.
    Meanwhile, it holds that $h_2 \circ (\pi_{f}\times \mathrm{id}) = h_1$.
    Moreover, both $X\times Y$ and $\Rf \times Y$ are compact Hausdorff spaces, and the map $\pi_{f}\times \mathrm{id}$ is onto and has connected fibers.
    Therefore, by Theorem \ref{thm:naturality}, the map $\Reeb{h_1} \to \Reeb{h_2}$ by $[x,y] \mapsto [[x],y]$ is a homeomorphism.

\end{proof}

We now introduce the notion of smoothing of a Reeb graph~\cite{dmp16}, which is a crucial construction we use for proving the stability of Reeb graphs later on.

\begin{definition}[Smoothing of a Reeb graph~\cite{dmp16}]
    Given $r\geq 0$, the \emph{r-smoothing} $\Rfr$ is the Reeb graph of the map $f^r: X \times [-r,r] \to \R$, $(x,t) \mapsto f(x)+t$. The map on $\Rf^r$ induced by $f^r$ is denoted by $\barf^r$. The image of $(x,t)$ under the Reeb quotient is denoted by $[x,t]$. We identify $\Reebs{f}{0}$ with $\Rf$.
\end{definition}

\begin{remark}
    The``smoothing" construction we present is also mentioned in \cite{dmp16}, though it is not the primary definition used in that reference. In \cite{dmp16} (and also in the subsequent  literature), the predominant definition of  $r$-smoothing of a Reeb graph $\Rf$ is the one described by $\Rf \times [-r,r] \to \R$, $([x],t) \mapsto f(x)+t$. While our approach may at first seem distinct, both are fundamentally equivalent. This equivalence is implied by Theorem 4.2.1 of \cite{dmp16} and further demonstrated by setting $s=0$ in item ii) of Theorem~\ref{thm:flow}.
    \end{remark}

\begin{theorem}[Naturality of smoothings]\label{thm:naturality_smoothing}
    Let $f: X \to \R$, $g: Y \to \R$ and $\psi: X \to Y$ be continuous maps such that $f= g \circ \psi$. Then $\bar{\psi}^r: \Rfr \to \Rgr$, $[x,t] \mapsto [\psi(x),t]$ is a well defined continuous map. Furthermore, if $\psi$ is onto and has connected fibers, then $\bar{\psi}^r$ is a homeomorphism.
\end{theorem}
\begin{proof}
    Let $\psi^r: X \times [-r,r] \to Y \times [-r,r]$, $(x,t) \mapsto (\psi(x),t)$. Then $\psi^r$ is a map between compact Hausdorff spaces and it is straightforward to check $f^r=g^r \circ \psi^r$.
    Then by Theorem \ref{thm:naturality}, $\bar{\psi}^r$ is a well defined continuous map.
    If $\psi$ is also onto and has connected fibers, then $\psi^r$ is also onto and has connected fibers.
    Then by Theorem \ref{thm:naturality}, $\bar{\psi}^r$ is a homeomorphism.
\end{proof}

The family of smoothings of a Reeb graph has some nice properties that we describe in the following theorem. To this end, we define a \emph{persistent family} of topological spaces
 to be a family of topological spaces $(X_r)_{r \geq 0}$ together with continuous structure maps $\pi^{r, s}: X_r \to X_s$ for $0 \leq r \leq s$ such that $\pi^{r, r} = \mathrm{id}$ and $\pi^{r, s} \circ \pi^{s, t} = \pi^{r, t}$ for $0 \leq r \leq s \leq t$.

\begin{theorem}[Flow of smoothings]\label{thm:flow}
    Let $f: X \to \R$ be a continuous map. Then, the following properties hold:
    \begin{enumerate}[i)]
        \item Smoothings form a persistent family under the continuous structure maps $\pi_f^{r,s}: \Rfr \to \Rfs$, $[x,t] \mapsto [x,t]$ for $0 \leq r \leq s$.
        \item The map $\Reebs{\barf^r}{s} \to \Reebs{f}{r+s}$, $[[x,t],u] \mapsto [x,t+u]$ is a well defined homeomorphism.
    \end{enumerate}
    These two properties show that the family of smoothings $(\Rfr)_{r \geq 0}$ is the result of applying a \emph{flow} to $\Rf$ in the sense of~\cite{de2017theory}.
\end{theorem}
\begin{proof}
    The family of spaces $(X \times [-r,r])_{r \geq 0}$ is a persistent family where the structure maps are inclusions. In particular, the inclusion maps commute with the maps $f^r: X \times [-r,r] \to \R$. Hence by Theorem \ref{thm:naturality_smoothing}, the family of smoothings $(\Rfr)_{r \geq 0}$ is a persistent family. This proves part $i)$.

    We now prove part $ii)$.
    To this end, we will first introduce an intermediate Reeb graph $\Reeb{h}$, where $h: X \times [-r,r] \times [-s,s] \to \R$, $(x,t,u) \mapsto f(x)+t+u$. By item $ii)$ of Corollary~\ref{cor:idempotencyAndProduct}, the map $\Reeb{h} \to \Reebs{\barf^r}{s}$, $[x,t,u] \mapsto [[x,t],u]$ is a homeomorphism. Next, we will show that $\Reeb{h}$ is also homeomorphic to $\Reebs{f}{r+s}$. To see this, we consider the map $\psi: X \times [-r,r] \times [-s,s] \to X \times [-r-s,r+s]$, $(x,t,u) \mapsto (x,t+u)$.
    Then $\psi$ is onto and the fibers of $\psi$ are homeomorphic to closed line segments and hence connected. Additionally, we have $h = f^{r+s} \circ \psi$.
    Then by Theorem \ref{thm:naturality_smoothing}, the map $\Reeb{h} \to \Reebs{f}{r+s}$, $[x,t,u] \mapsto [x,t+u]$ is a homeomorphism. This concludes the proof of part $ii)$.
\end{proof}

\section{Metrizing Reeb Graphs}

From now on, $(X,d_X)$ will denote a \textbf{compact geodesic space}, and $f: X \to \R$ will denote a $1$-Lipschitz function.
$\Rf$ is the Reeb graph of $f$ with induced function $\barf: \Rf \to \R$.

We aim to introduce a metric on $\Rf$. Our first step involves assigning a length to every continuous path in $\Rf$, a process referred to as establishing a 'length structure' (see Definition~\ref{def:length_structure} in Appendix~\ref{app:metric_geometry} for details). In the subsequent step, we assign the distance between two points as the infimum of the lengths of all paths connecting them.
Our goal is to prove that this derived metric metrizes the quotient topology on $\Rf$, see Proposition~\ref{prop:reebMetric}.

\begin{definition}[Metric on Reeb graph]\label{def:df}
    Let $\mathcal{A}$ denote the set of continuous paths in $\Rf$ para\-metrized by closed intervals $[a, b]\subset \R$. Let $L_f: \mathcal{A} \to [0,\infty]$, $L_f(\beta):=\lr(\barf \circ \beta)$. The \emph{Reeb metric} $\df$ on $\Rf$ is defined as
    $$\df([x],[y]):=\inf\{L_f(\beta): \beta \in \mathcal{A}, \beta(a)=[x], \beta(b)=[y] \}.$$
\end{definition}

In what follows, we verify that $L_f$ is a length structure on $\Rf$ and that $\df$ is its associated length metric.
In particular, this justifies the giving the name `Reeb metric' to $\df$.

{
\begin{proposition}\label{prop:reebMetric}
    $L_f$ is a length structure on $\Rf$. Furthermore, 
    $\df$ is a metric on the topological space $\Rf$ and the topology induced by $\df$ is no coarser than the topology of $\Rf$.
\end{proposition}

\begin{proof}
    According to Proposition~\ref{prop:length_metric}, it suffices to check that $L_f$ is indeed a length structure on $\Rf$.
    That is, $L_f$ satisfies the following properties:
    \begin{enumerate}
        \item \textbf{Additivity of Path Length}: For any path \( \gamma:[a, b] \rightarrow X \) and any \( c \in[a, b] \),
        \[ L_f\left(\gamma_{\left.\right|_{[a, b]}}\right)=L_f\left(\gamma_{\left.\right|_{[a, c]}}\right)+L_f\left(\gamma_{\left.\left.\right|_{[c, b]}\right]}\right) \]
  \item \textbf{Continuity of Length}:
We have
  $L_f(\gamma|_{[a,\cdot]}):[a,b] \to [0,\infty]$ is continuous whenever $L_f(\gamma)<\infty$,
  \item \textbf{Invariance under Reparameterizations}: The length is invariant under reparameterizations,
        \[ L_f(\gamma \circ \varphi)=L_f(\gamma) \]
        for any homeomorphism \( \varphi : [a, b] \to [a, b]\).
  \item \textbf{Compatibility with Topology}: $L_f$ is compatible with
    the topology of $X$ in the sense that for any neighborhood $U$ of
    a point $x \in X$, the length of any path connecting $x$ with
    points of the complement of $U$ is strictly positive:
  \[
      \inf \left\{L(\gamma): \gamma(a)=x, \gamma(b) \in X \backslash U\right\}>0
  \]
    \end{enumerate}
    The first three properties are clear from the definition of
    $L_f$. It remains to show the last one.  To see this. Let
    $c:=f(x)$. By Proposition \ref{prop:hausdorff}, there exists $\epsilon>0$ such that the path component $\barf^{-1}((c-\epsilon,c+\epsilon))$ of $[x]$ is contained in $U$. Hence, for any $\beta \in \mathcal{A}, \beta(a)=[x], \beta(b) \notin U $, one has $|\barf(\beta(b))-c|>\epsilon$, which implies that $L_f(\beta)>\epsilon$.
\end{proof}}

Now we are ready to define \emph{Reeb metric graph}.

\begin{definition}[Reeb metric graph and smoothings]
    For any given metric space $(X, d_X)$ and a function $f$ on it, we define its associated
     \emph{Reeb metric graph} as $(\Rf,\df)$ where $\df$ is the metric induced by $L_f$. Note that under the metric $\df$, $\barf: \Rf \to \R$ is $1$-Lipschitz. We denote the Reeb metric of the $r$-smoothing $\Rfr$ by $\dfr$.
\end{definition}

The following proposition shows that $\df$ metrizes $\Rf$:

\begin{proposition}\label{prop:reebMetrization}
    The topology induced by the Reeb metric $\df$ coincides with the quotient topology on $\Rf$.
    In particular, $(\Rf,\df)$ is a compact metric space as
      $\Rf$ is compact (cf. Corollary~\ref{cor:reeb_hausdorff}).
\end{proposition}
\begin{proof}
    The topology induced by $\df$ is no coarser than the quotient topology by Proposition~\ref{prop:reebMetric}. It remains to show that the identity map from $R_f$ with quotient topology to $(R_f,d_f)$ is continuous. This is equivalent to showing that $(X,d_X) \to (\Rf,\df)$ is continuous which, with slight abuse of notation, we still denote it as $\pf$.
    Indeed, $\pf$ is $1$-Lipschitz and hence continuous.
    As $(X, d_X)$ is a geodesic space, for any $x, x'\in X$, we have that
    \[
        d_X(x, x') = \inf_\gamma \mathrm{length}_X(\gamma),
    \]
    where $\gamma$ ranges over all continuous paths from $x$ to $x'$.
    Note that $\pf \circ \gamma$ is a continuous path from $[x]$ to $[x']$ in $\Rf$.
    Then we have $\lr(\pf \circ \gamma) \leq \mathrm{length}_X (\gamma)$ as $f$ is $1$-Lipschitz.
    Therefore, we have $d_f([x], [x']) \leq d_X(x, x')$ and the map $\pf$ is $1$-Lipschitz.
\end{proof}

Properties like naturality and idempotency of topological Reeb graph in Section~\ref{sec:reeb}
can be upgraded to the metric setting by using the following lemma.

\begin{lemma}\label{lem:metricReebGraphCategorical}
    Let $f: X \to \R$, $g : Y \to \R$ be $1$-Lipschitz maps on compact geodesic spaces with Reeb graphs $\Rf$, $\Rg$ respectively
    . Let $\bar{\psi}: \Rf \to \Rg$ be a continuous map such that $\barf= \barg \circ \bar{\psi}$. Then $\bar{\psi}$ is $1$-Lipschitz. Furthermore, if $\bar{\psi}$ is a homeomorphism, then it is an isometry.
\end{lemma}

\begin{proof}
    For any $[x], [x'] \in \Rf$, we have
    \[
        \df([x], [x']) = \inf_\gamma \lr(\barf \circ \gamma),
    \]
    where $\gamma$ ranges over all continuous paths from $[x]$ to $[x']$ in $\Rf$.
    Note that $\bar{\psi} \circ \gamma$ is a continuous path from $\bar{\psi}([x])$ to $\bar{\psi}([x'])$ in $\Rg$.
    Then we have $\lr(\barf \circ \gamma) = \lr(\barg \circ \bar{\psi} \circ \gamma) \geq \dg(\bar{\psi}([x]), \bar{\psi}([x'))$ which implies $\df([x], [x']) \geq \dg(\bar{\psi}([x]), \bar{\psi}([x'))$. This shows that $\bar{\psi}$ is $1$-Lipschitz.

    If $\bar{\psi}$ is a homeomorphism, then its inverse is also $1$-Lipschitz by the argument above. Therefore, it is an isometry.
\end{proof}

\begin{proposition}\label{prop:metricReebGraphCategorical}
    Let $f: X \to \R$, $g: Y \to \R$ be $1$-Lipschitz maps on compact geodesic spaces with Reeb graphs $\barf: \Rf \to \R$, $\barg: \Rg \to \R$, 
 respectively.
    Let $\psi: X \to Y$ be continuous maps such that $f= g \circ \psi$. Then, the following properties hold:
    \begin{enumerate}[i)]
        \item (Naturality) $\bar{\psi}^r: \Rfr \to \Rgr$, $[x,t] \mapsto [\psi(x),t]$ is a $1$-Lipschitz. Furthermore, if $\psi$ is onto and has connected fibers, then $\bar{\psi}^r$ is an isometric bijection.
        \item (Idempotency) $(\Rf,\df)=(\Reeb{\barf},\DistF{\barf})$.
        \item (Product) Let $h_1: X \times Y \to \R$, $(x,y) \mapsto f(x)+g(y)$ and $h_2 : \Rf \times Y \to \R$, $([x],y) \mapsto f(x)+g(y)$. Then the map $\Reeb{h_1} \to \Reeb{h_2}$, $[x,y] \mapsto [[x],y]$ is an isometry.
        \item Smoothings form a persistent family under the $1$-Lipschitz structure maps $\pi_f^{r,s}: \Rfr \to \Rfs$, $[x,t] \mapsto [x,t]$ for $0 \leq r \leq s$.
        \item The map $\Reebs{\barf^r}{s} \to \Reebs{f}{r+s}$, $[[x,t],u] \mapsto [x,t+u]$ is an isometry.
    \end{enumerate}
\end{proposition}
\begin{proof}
    The results follow from applying Lemma~\ref{lem:metricReebGraphCategorical} to the results in Theorem \ref{thm:naturality_smoothing}, Corollary \ref{cor:idempotencyAndProduct}, and Theorem \ref{thm:flow}.
\end{proof}

The following result shows that smoothings of Reeb graphs are right continuous in the sense of Gromov-Hausdorff distance.

\begin{proposition}\label{prop:rightContinuity}
    For any $r \geq 0$, and let $\pi_f^{r,s}: \Rf \to \Rfs$ be the structure map.
    Then, we have $\dis(\pi_f^{r,s})$ goes to $0$ as $s\to r^+$.
    Furthermore, 
    $\Rfs$ converges to $\Rfr$ in the sense of Gromov-Hausdorff distance as $s \to r^+$.
\end{proposition}

\begin{proof}
    We first assume that $r=0$.
    Note that, the image $\pi_f^{0,s}(\Rf)$ is an $s$-net in $\Rfs$.
    Let $\dis(\pi_f^{0, s})$ denote the distortion of $\pi_f^{0,s}$.
    Then we have the bound $\dgh(\Rf, \Rfs) \leq 2\dis(\pi_f^{0, s}) + 2s$, see for example \cite[Corollary~7.3.28]{bbi01}.
    To bound $\dis(\pi_f^{0, s})$, we can alternative consider the pullback metric $\hat{d}^{s}_f:=(\pi_f^{0, s})^* d_f$ on $\Rf$, that is, for every $[x], [x'] \in \Rf$, we define $\hat{d}^{s}_f([x], [x']) := d_f(\pi_f^{0, s}([x]), \pi_f^{0, s}([x']))$.
    Therefore, we can represent $\dis(\pi_f^{0, s})$ as
    \[
        \dis(\pi_f^{0, s}) = \sup(\||d_f - \hat{d}^{s}_f\|_\infty)
    \]
    where $\|\cdot\|_\infty$ denotes the sup norm.
    Note $\hat{d}^{s}_f$ is a pseudometric on $\Rf$.
    Since the structure maps $\pi_f^{t, s}$ are $1$-Lipschitz, we have $\hat{d}^{s'}_f \leq \hat{d}^{s}_f$ for every $s'> s\geq 0$.
    We further define the limiting pseudometric $\hat{d}^{0}_f$ on $\Rf$ as
    \[
        \hat{d}_f([x], [x']) = \lim_{s \to 0} \hat{d}^{s}_f([x], [x']) = \sup_{s>0} \hat{d}^{s}_f([x], [x']).
    \]
    for any $[x], [x'] \in \Rf$.
    As the supremum of pseudometrics, $\hat{d}_f$ is a pseudometric on $\Rf$.
    We will furthermore show that $\hat{d}_f$ is indeed a metric on $\Rf$. To see this, we need to show that $\hat{d}_f([x], [x']) = 0$ implies $[x] = [x']$.
    For any integer $n>0$, if $\hat{d}^{1/n}_f([x], [x']) = 0$, then $d_f^{1/n}([x, 0], [x', 0]) = 0$ and $(x, 0), (x', 0)$ are in the same connected component $C_n$ of a level set of $f^{1/n}: X \times [-1/n, 1/n] \to \R$.
    We use $C_n'$ to denote the image of $C_n$ under the inclusion $X \times [-1/n, 1/n] \to X \times [-1, 1]$. Then as $n$ goes to infinity, $C_n'$ form a decreasing sequence of compact connected sets in $X \times [-1, 1]$.
    By \cite[Corollary~6.1.19]{engelking1989general}, $\cap_n C_n'$ is connected and contained in $X \times \{0\}$.
    Note that, $f^1$ on $X \times \{0\}$ restricts to $f$ on $X\times \{0\}$, and hence $\cap_n C_n' = C \times \{0\}$ for some connected subspace $C$ of the level set $f$ on $X$ containing $x$ and $x'$.
    Therefore, $[x] = [x']$.
    We have shown that $\hat{d}_f$ is a metric on $\Rf$.

    Meanwhile, we have $\hat{d}_f \leq d_f$ as $\hat{d}_f^s\leq d_f$ for all $s>0$, and hence the identity map $(\Rf, \hat{d}_f) \to (\Rf, d_f)$ is a continuous bijection between compact metric spaces, so it is a homeomorphism.
    Therefore, $\hat{d}_f$ metrizes $\Rf$.
    Next, we will show that $\hat{d}_f \geq d_f$. To this end, we will first show that $(\Rf, \hat{d}_f)$ is a length space.
   
        Since $(\Rf, d_f)$ is a compact metric space
        (cf. Proposition~\ref{prop:reebMetrization}), we can apply
        \cite[Theorem~2.4.16]{bbi01} (Theorem~\ref{thm:midpoint} in
        the Appendix) so that it suffices to show that for every $[x], [x'] \in \Rf$ , there exists $[x''] \in \Rf$ such that $\hat{d}_f([x], [x'']) = \hat{d}_f([x], [x'])/2$ and $\hat{d}_f([x'], [x'']) = \hat{d}_f([x], [x'])/2$.
    
    Note that each $(\Rf^{1/n}, d_f^{1/n})$ is a length space, hence there exists $[x_n'', t_n]$ such that $d_f^{1/n}([x, 0], [x_n'', t_n]) = d_f^{1/n}([x, 0], [x', 0])/2$ and $d_f^{1/n}([x', 0], [x_n'', t_n]) = d_f^{1/n}([x, 0], [x', 0])/2$.
    By the compactness of $X$, we can assume that $x_n''$ converges to $x''$ in $X$ when $n$ goes to infinity.
    Also, note that $t_n\leq 1/n$ and hence $t_n$ converge to $0$ as $n$ goes to infinity.
    Then we have
    \begin{align*}
        \hat{d}_f([x], [x'']) &\leq \lim_{n \to \infty} d_f^{1/n}([x, 0], [x'', 0]) \\
        &\leq \limsup_{n \to \infty} \left( d_f^{1/n}([x, 0], [x_n'', t_n]) +d_f^{1/n}([x_n'', t_n], [x_n'', 0]) + d_f^{1/n}([x_n'', 0], [x'', 0]) \right)\\
        &\leq \limsup_{n \to \infty} \frac{d_f^{1/n}([x, 0], [x', 0])}{2} + 1/n + d_X(x'', x_n'')
    \end{align*}
    Therefore, one has $\hat{d}_f([x], [x'']) \leq \hat{d}_f([x], [x'])/2$ by letting $n$ goes to infinity.
    Similarly, we have $\hat{d}_f([x'], [x'']) \leq \hat{d}_f([x], [x'])/2$.
    Therefore, by triangle inequality, both inequalities are equalities and $[x'']$ is a midpoint between $[x]$ and $[x']$ with respect to $\hat{d}_f$.
    To show that $\hat{d}_f \geq d_f$, it suffices to show that, for any continuous path $\gamma$ connecting $[x]$ and $[x']$ in $\Rf$, we have $$\hat{d}_f([x], [x'])= \mathrm{length}_{(\Rf, \hat{d}_f)} (\gamma)\geq \lr(\barf \circ \gamma).$$

        The inequality follows from the fact that $d_f^s([x'', 0], [x''', 0]) \geq |\barf([x''])-\barf([x'''])|$. Consequently, $\hat{d}_f^s([x''], [x''']) \geq |\barf([x''])-\barf([x'''])|$ for any $s>0$ and $[x''], [x'''] \in \Rf$. Therefore, when we consider the supremum over $s$, the inequality $\hat{d}_f([x], [x'])= \mathrm{length}_{(\Rf, \hat{d}_f)} (\gamma)\geq \lr(\barf \circ \gamma)$ holds.
    
    Therefore, we have that $\hat{d}_f = d_f$ and $\hat{d}_f^s$ converges to $d_f$ pointwisely as $s$ goes to $0$.
    By Dini's theorem \cite[Theorem~2.4.10]{dudley2002real}, the convergence is uniform.
    Therefore, we have
    \[
        \lim_{s \to 0} \dis(\pi_f^{0,s}) = \lim_{n \to \infty} \dis(\pi_f^{0,1/n}) = \lim_{n \to \infty}\sup(|\df-\hat{d}^{1/n}_f|)=0.
    \]
    This concludes the proof for the case $r=0$.
    For general $r>0$, there is a commutative diagram from Proposition \ref{prop:metricReebGraphCategorical} part $v)$,  with bijective and isometric vertical maps as follows:
    $$\begin{tikzcd}
            \Reeb{\barf^r} \arrow[r,"\pFunc{\barf^r}{0}{s-r}"] \arrow[d] & \Reebs{\barf^r}{s-r} \arrow[d] \\
            \Rfr \arrow[r,"\pfrs"] & \Rfs
        \end{tikzcd}$$
    Hence, by replacing $\Rf$ with $\Reeb{\barf^r}$, it reduces to the case $r=0$.
\end{proof}

\begin{remark}\label{rem:non_geodesic}
    Note that for a $1$-Lipschitz function $\phi: E \to \R$ on a compact metric space $E$, the Reeb metric $d_\phi$, in general, may not metrize $\Rphi$. For example, let $E$ be the topologist's sine curve
    \begin{equation*}
        T:=\left\{\left(x, \sin \frac{1}{x}\right): x \in(0,1]\right\} \cup\{(0,0)\},
    \end{equation*}
    and $\phi: E \to \R$ be the height function that sends $(x, \sin \frac{1}{x})$ to $\sin \frac{1}{x}$ and $(0,0)$ to $0$.
    Then, $\Rphi$ is homeomorphic to $E$. However, the Reeb metric $d_\phi$ does not metrize $\Rphi$ as the point $(0,0)$ is an isolated point under the metric of $d_\phi$.

    On the other hand, one can show that if $\phi$ factors through a geodesic space with connected fibers, then $d_\phi$ metrizes $\Rphi$.
    To see this, assume that there is a $1$-Lipschitz map $\pi: E \to X$ onto a compact geodesic space, $\pi$ has connected fibers, $f: X \to \R$ $1$-Lipschitz and $\phi=f \circ \pi$.
    Then, by Theorem \ref{thm:naturality}, the map $\bar{\pi}: \Rphi \to \Rf$, $[x] \mapsto [\pi(x)]$ is a homeomorphism.
    By Lemma~\ref{lem:metricReebGraphCategorical}, $\bar{\pi}$ is an isometry.
    The Proposition~\ref{prop:reebMetrization} shows that $\df$ metrizes $\Rf$ and Lemma~\ref{lem:metricReebGraphCategorical} identifies $(\Rphi, d_\phi)$ with $(\Rf, \df)$ and hence $d_\phi$ metrizes $\Rphi$.
\end{remark}

The following is a key result providing a relationship between paths in $X$ and paths in $\Rf$.

\begin{lemma}\label{lem:pathLifting}
    Assume that $E$ is a compact locally path connected space, and $\pi: E \to X$ is a continuous surjection with connected fibers.
    Then for any $\epsilon>0$, arbitrary points $p, p'\in E$, and any path $\beta$ that goes from $\pi(p)$ to $\pi(p')$ in $X$, there exists a path $\alpha$ that goes from $p$ to $p'$ in $E$
    such that $d_X(\pi \circ \alpha, \beta) \leq \epsilon$.
    In particular $\|f \circ \pi \circ \alpha - f \circ \beta\|_\infty \leq \epsilon$. Note that this result applies to the Reeb quotient map $\pi_f : X \to \Rf$, as it has connected fibers.
\end{lemma}
\begin{proof}
    By the Lebesgue number lemma, there exists a partition $0 = t_0 < \dots < t_n=1$ and open sets $(V_i)_{i=1}^n$ in $X$ such that $\beta|_{[t_{i-1},t_i]} \subseteq V_i$, $\diam(V_i) \leq \epsilon$, and each $V_i$ is path connected. Let $U_i:=\pi^{-1}(V_i)$, then $U_i$ is connected by \cite[Theorem~6.1.29]{engelking1989general}, hence $U_i$ is also path connected by the local path connectedness of $E$. Let $(p_i)_{i=0}^n$ be such that $p_0=p$, $p_n=p'$, and $\pi(p_i)=\beta(t_i)$ for all $i$. There exists continuous $\alpha_i:[t_{i-1},t_i] \to U_i$ such that $\alpha_i(t_{i-1})=x_{i-1}$ and $\alpha_i(t_i)=p_i$. Let $\alpha: [0,1] \to X$ be the map defined piecewisely by $\alpha|_{[t_{i-1},t_i]}=\alpha_i$. Then $\alpha$ is a path in $E$ that goes from $p$ to $p'$ and  $d_X(\pi \circ \alpha, \beta) \leq \epsilon$.
\end{proof}

\section{Structure of Reeb Metric Graphs}\label{sec:structure}

In this section, we delve deeper into the structure of the Reeb metric graph. First, we define what it means to be a \emph{core} of a metric space.

\medskip
\noindent\textbf{Notation:} Given a subspace $A$ of a metric space $(X,d_X)$, $(X/A,d_A)$ denotes quotient metric space (see \cite[Definition~3.1.12]{bbi01}) of $X$ under the equivalence relation $x \sim y$ if $x=y$ or $x,y \in A$. We denote the $1$-Lipschitz metric quotient map by $\pi_A:X \to X/A$. See Appendix \ref{app:quotient} for more details about the quotient metric spaces.

\begin{definition}[Core of a metric space]
    Let $(X,d_X)$ be a compact metric space. We call a closed subspace $A$ of $X$ a \emph{core} of $X$ if $X/A$ is a metric tree. See Appendix \ref{app:hyperbolicity}, \ref{app:core} for more details about metric trees and cores.
\end{definition}

    The main result in this section is Theorem~\ref{thm:structure} showing the Reeb metric graph is a compact geodesic space with Lebesgue covering dimension at most one. The Lebesgue covering dimension of a topological space is defined as follows:
\begin{definition}[Lebesgue covering dimension]\label{defn:Lebesgue_dim}
	Let $X$ be a topological space. We say that $X$ has Lebesgue
        dimension $\leq n$ if, for any open covering $\{U_i\}$ of $X$,
        there is a refinement $\{V_i\}$ of $\{U_i\}$ with the property that every $x\in X$ lies in at most $n+1$ of the $V_i$s. We say that $X$ is $n$-dimensional if $X$ has dimension $\leq n$, but does not have dimension $\leq n-1$.
\end{definition}

\begin{theorem}[Structure of Reeb Metric Graph]\label{thm:structure}
    Let $(X,d_X)$ be a geodesic space, $f: X \to \R$ be a $1$-Lipschitz function. Then $\Rf$ is a compact geodesic space with Lebesgue covering dimension at most one. Furthermore, if $\betti(X) < \infty$, then $\Rf$ has a core $G$ such that $G$ is a finite metric graph and $\betti(G) \leq \betti(X)$.
\end{theorem}

We prove this theorem at the end of the section.

\medskip

    We now need to introduce a notion of distance that compares two
    pairs of metric spaces with real valued functions defined on
    them. Such notion has been considered
    in~\cite{carlsson2010multiparameter, chazal2009gromov} and has found applications in topological data analysis.

    Let us first recall the definition of correspondence and also that of the standard Gromov-Hausdorff distance.

\begin{definition}[Correspondence]\label{def:correspondence}
    A correspondence $R$ between two given sets $X,Y$, is a relation between them such that for all $x$ in $X$, there exists a $y_0$ in $Y$ such that $x \, R \, y_0$ and for each $y$ in $Y$, there exists an $x_0$ in $X$ such that $x_0 \, R \, y$.
\end{definition}

\begin{definition}[Distortion of a correspondence]
	Let $X,Y$ be bounded metric spaces and $R$ be a correspondence between $X,Y$. The metric distortion $\dis(R)$ of the correspondence $R$ is defined as
	$$\dis(R):=\sup_{(x,y),(x',y') \in R} |d_X(x,x')-d_Y(y,y')| .$$
\end{definition}

If both $R, S$ are correspondences between $X,Y$, we say that $S$ is a subcorrespondence of $R$ if $S \subset R$.

There are several equivalent ways of defining the Gromov-Hausdorff distance (see~\cite[Section 7.3]{bbi01}). In this paper, we use the following:

\begin{definition}[Gromov-Hausdorff distance]\label{def:gh_distance}
    Let $X,Y$ be bounded metric spaces.
    The Gromov-Hausdorff distance $\dgh(X,Y)$ between them is defined as $$\dgh(X,Y):=\frac{1}{2} \, \inf \{ \dis(R): R \text{ is a correspondence between } X,Y \}.$$
\end{definition}

    Then, the Gromov-Haudorff distance between two metric spaces with functions defined on them is defined as follows:

\begin{definition}[{\cite[Definition~2]{carlsson2010multiparameter}}]
    Let $X$ and $Y$ be metric spaces and $f: X \to \R$ and $g: Y \to \R$ be functions.
     Given a \emph{correspondence} \( R \) between \( X \) and \( Y \), we define:
    \[ \dis_{f,g}(R):=\max\left(\dis(R), \sup_{(x,y) \in R} |f(x)-g(y)|\right). \]
    Using this, we define
    $$ \dgh((X,f),(Y,g)):= \frac12 \inf_R \dis_{f,g}(R),$$
    where the infimum is taken over all correspondences between $X$ and $Y$. If $f=d_X(p,\cdot)$ and $g=d_Y(q,\cdot)$, then we use the notation $\dgh((X,p),(Y,q))$.
\end{definition}
\begin{remark}
    As in the case of Gromov-Hausdorff distance \cite[Chapter
    7]{bbi01}, by considering the composition of correspondences, one can show that $\dgh((X,f),(Z,h)) \leq \dgh((X,f),(Y,g))+\dgh((Y,g),(Z,h))$.
\end{remark}

\begin{lemma}\label{lem:correspondence_nbhd}
    Assume $\dgh((X,f),(Y,g)) <r$. Then, there exists a compact metric space $E$ with $1$-Lipschitz maps $\pi_X: E \to X$, $\pi_Y: E \to Y$ with connected fibers, and $1$-Lipschitz functions $\phi: E \to \R$, $\psi: E \to \R$ such that $\phi=f \circ \pi_X$, $\psi= g \circ \pi_Y$, $\|\phi-\psi\|_{\infty} < 6r$.
\end{lemma}
\begin{proof}
    Let $R$ be a correspondence between $X$ and $Y$ such that $\dis_{f,g}(R) < 2r$. Let $E$ be the closed $2r$ neighborhood of $R$ in $X \times Y$ with respect to the $\ell^\infty$-product metric, that is, $$E:=\{(p,q) \in X \times Y: d_X(p,x)\leq 2r, d_Y(q,y) \leq 2r \text{ for some } (x,y) \in R\}.$$
    Let $\pi_X: E \to X$, $(p, q) \mapsto p$, $\pi_Y: E \to Y$, $(p, q) \mapsto q$. Let $\phi=f \circ \pi_X$ and $\psi=g \circ \pi_Y$.
    Let $(p,q) \in E$, and let $(x,y) \in R$ such that $d_X(p,x), d_Y(q,y) \leq 2r$. Then, we have
    $$|\phi(p,q)-\psi(p,q)|=|f(p)-g(q)| \leq |f(p)-f(x)|+|f(x)-g(y)|+|g(y)-g(q)| < 6r. $$

    It remains to show that $\pi_X, \pi_Y$ have connnected fibers. Let $p \in X$, and $C$ be the closed subset of $Y$ such that $\pi_X^{-1}(p)=\{p\} \times C$. We are going to show that $C$ is path connected. Let $q_0 \in Y$ such that $(p,q_0) \in R$. Note that $q_0 \in C$. Let $q \in C$. Let us construct a path in $C$ from $q$ to $q_0$. There exists $(x,y) \in R$ such that $d_X(p,x), d_Y(q,y) \leq 2r$. Since $(x,y), (p,q_0) \in R$, we have $d_Y(y,q_0) \leq d_X(p,x)+2r \leq 4r$. Let $y'$ be a midpoint between $y$ and $q_0$. The following diagram shows how the points are related:
    $$\begin{tikzcd}
            q_0 \arrow[r,"\leq 2r",no head] & y' \arrow[r,"\leq 2r",no head]& y \arrow[r,"\leq 2r",no head] & q \\
            p \arrow[r,"\leq 2r",no head] \arrow[u,"R",leftrightarrow] & x \arrow[ur,"R",leftrightarrow] &  &
        \end{tikzcd}$$
    For any point $z$ that is a midpoint between
    $y$ and $q$ (i.e $d_Y(y,q)=d_Y(y,z)+d_Y(z,q)$) or $y$ and $y'$,
    $(p,z)$ is contained in the closed $2r$ neighborhood of $(x,y) \in R$, hence $z \in C$. For any point $z$ in between $q_0$ and $y'$, $(p,z)$ is contained in the closed $2r$ neighborhood of $(p,q_0) \in R$, hence $z \in C$. Hence the path from $q$ to $q_0$ obtained by concatenating geodesics from $q$ to $y$, $y$ to $y'$ and $y'$ to $q_0$ is contained in $C$. This shows $C$ is path connected, hence fibers of $\pi_X$ are connected. Similarly fibers of $\pi_Y$ are connected.
\end{proof}

The following proposition shows how smoothings of the Reeb graph of 
one space approximates the Reeb graph of another space.
\begin{proposition}\label{prop:approximation_smoothing}
    Let $r>0$ such that $\dgh((X,f),(Y,g)) < r$. Then,
    $$\dgh(\Rf,\Rg^{6r}) \leq 12r + \dis(\pi_f^{0,12r})/2.$$
\end{proposition}
\begin{proof}
    Let $E$, $\pi_X: E \to X$, $\pi_Y$, $\phi: E \to \R$, $\psi: E \to \R$ be as in Lemma \ref{lem:correspondence_nbhd}. By Proposition \ref{prop:metricReebGraphCategorical}, $\Rphi$ is isometric to $\Rf$, $\Rpsi^{6r}$ is isometric to $\Rg^{6r}$, and $\dis(\pi_f^{12r})=\dis(\pi_\phi^{0,12r})$ (see Remark \ref{rem:non_geodesic}). Hence, 
    it is enough to show that
    $$\dgh(\Rphi,\Rpsi^{6r}) \leq 12r + \dis(\pi_\phi^{0,12r})/2.$$
    
    Let $\eta: E \to \R$, $\eta := \phi - \psi$, $F: E \to E \times [-6r,6r]$, $z \mapsto (z,\eta(z))$, $G: E \times [-6r,6r] \to E \times [-12r,12r]$, $(z,t) \mapsto (z,t-\eta(z)) $. The following diagram commutes:
    $$\begin{tikzcd}
            E \arrow[r,"F"] \arrow[rd,"\phi"'] & E \times [-6r,6r] \arrow[r,"G"] \arrow[d,"\psi^{6r}"] &  E \times [-12r,12r] \arrow[ld, "\phi^{12r}"] \\
            & \R &
        \end{tikzcd}$$
    By Proposition \ref{prop:metricReebGraphCategorical}, $F$ and $G$ induce $1$-Lipschitz maps between the corresponding Reeb metric graphs, hence we have, for any $[w],[z] \in \Rphi$,
    $$d_\phi([w],[z]) \geq d_\psi^{6r}([w,\eta(w)],[z,\eta(z)]) \geq d_\phi^{12r}([w,0],[z,0]).  $$
    Let $\bar{R}$ be the correspondence between $\Rphi, \Rpsi^{6r}$ given by $\bar{R}:=\{([z],[z,t]): z \in E, t \in [-6r,6r] \}$.
    Then, we have
    \begin{equation*}
        \begin{split}
            |d_\phi([w],[z])-d_\psi^{6r}([w,s],[z,t])| & \leq |d_\phi([w],[z])-d_\psi^{6r}([w,\eta(w)],[z,\eta(t)])| + \\ &|d_\psi^{6r}([w,\eta(w)],[z,\eta(z)])-d_\psi^{6r}([w,s],[z,t])| \\
            & \leq \dis(\pi_\phi^{0,12r}) + d_\psi^{6r}([w,\eta(w)],[w,s]) + d_\psi^{6r}([z,\eta(z)],[z,t]) \\
            & \leq \dis(\pi_\phi^{0,12r}) + |\eta(w)-s| + |\eta(z)-t| \\
            & \leq \dis(\pi_\phi^{0,12r}) + 24r.
        \end{split}
    \end{equation*}
    This shows that $\dis(\bar{R}) \leq \dis(\pi_\phi^{0,12r}) + 24r$ and completes the proof.
\end{proof}

Before proving Theorem \ref{thm:structure}, we need one more result about the first Betti number 
of smoothings of a Reeb graph.

\begin{proposition}\label{prop:smoothingBetti}
     Let $X$ and $Y$ be compact metric spaces and $f: X \to \R$ and $g: Y \to \R$ be $1$-Lipschitz functions.
        Assume that the Reeb quotient map $\pi_{g^s}: Y \to \Rg^s$ induces a surjection between first homology groups for all $s>0$.
    Then, for any $r>0$ such that $\dgh((X,f),(Y,g))<r$, we have $$\betti(\Rg^{12r}) \leq \betti(X).$$
\end{proposition}

\begin{remark}
    It is shown in~\cite[Theorem~3.2]{Dey2013} that if $Y$ is a finite
    simplicial complex and $g$ is a piecewise linear function, then
    the Reeb quotient map induces a surjection between the first
    homology groups. Hence, the above proposition applies in this setting.
\end{remark}

\begin{proof}[Proof of Proposition \ref{prop:smoothingBetti}]
    Since $\pi_{g^s}: Y \to \Rg^s$ induces a surjection between first homology groups for all $s>0$ and the map $\pi_{g^{s'}}$ satisfies $\pi_{g}^{s, s'} \circ \pi_{g^s} = \pi_{g^{s'}}$ for $s \leq s'$, the map $\pi_g^{s, s'}$ must also induce a surjection on first homology groups
     $H_1(\Rg^s) \to H_1(\Rg^{s'})$ as well.
\medskip

    Let $E$, $\pi_X: E \to X$, $\pi_Y$, $\phi: E \to \R$, $\psi: E \to \R$ be as in Lemma \ref{lem:correspondence_nbhd}. Let $s=\sup|\phi-\psi|<6r$. As in the proof of Proposition \ref{prop:approximation_smoothing}, we have a commutative diagram of continuous maps
    $$\begin{tikzcd}
            E \arrow[r,"F"] \arrow[rd,"\psi"'] & E \times [-s,s] \arrow[r,"G"] \arrow[d,"\phi^{s}"] &  E \times [-2s,2s] \arrow[ld, "\psi^{2s}"] \\
            & \R &
        \end{tikzcd}$$
    where $G \circ F$ maps $z$ to $[z,0]$. Since $\psi = g \circ \pi_Y$ and the map $\pi_Y$ is $1$-Lipschitz with connected fibers,
    the naturality of the Reeb graph construction (Theorem \ref{thm:naturality}) applies and, as a consequence, the map $H_1(\Rpsi) \to H_1(\Rpsi^{2s})$ induced by $\overline{G\circ F}: \Rpsi \to \Rpsi^{2s}$ is surjective.
    Also note that the map $H_1(\Rpsi) \to H_1(\Rpsi^{2s})$ splits through $H_1(\Rphi^s)$.
    By considering an analogous construction involving $\pi_X: E\to X$, one sees that the map $H_1(\Rpsi) \to H_1(\Rpsi^{12r})$ is surjective and splits through $H_1(\Rphi^{6r})$.

    We have the following diagram:
    $$\begin{tikzcd}
        H_1(\Rphi^s) \arrow[r] \arrow[d, two heads] & H_1(\Rphi^{6r}) \arrow[d, two heads] \\
        H_1(\Rpsi^{2s}) \arrow[r, two heads] & H_1(\Rpsi^{12r})
    \end{tikzcd}$$
    where all the morphisms with the possible exception of the top horizontal one are guaranteed to be surjections. By Theorem \ref{thm:naturality_smoothing}, we have
    $$\betti(\Rg^{12r}) = \betti(\Rpsi^{12r}) \leq \dim(\mathrm{im}(H_1(\Rphi^s) \to H_1(\Rphi^{6r}))) = \dim(\mathrm{im}(H_1(\Rf^s) \to H_1(\Rf^{6r}))) .$$
    Hence, it is enough to show that 
    $$\mathrm{im}(H_1(\Rf^s) \to H_1(\Rf^{6r})) \subseteq \mathrm{im}(H_1(X \times [-6r,6r]) \to H_1(\Rf^{6r})). $$
    By Proposition \ref{prop:metricReebGraphCategorical}, smoothings form a flow, hence the above inclusion is equivalent to
    $$\mathrm{im}(H_1(R_{f^s}) \to H_1(R_{f^s}^{6r-s})) \subseteq \mathrm{im}(H_1((X \times [-s,s]) \times [-6r+s,6r-s]) \to H_1(R_{f^s}^{6r-s})). $$
    We will prove this by showing a more general fact that for any geodesic space $(Z,d_Z)$ with a $1$-Lipschitz map $h: Z \to \R$ and $\epsilon>0$,
    $$\mathrm{im}(H_1(\Rh) \to H_1(\Rh^\epsilon)) \subseteq \mathrm{im}(H_1(Z \times [-\epsilon,\epsilon]) \to H_1(\Rh^{\epsilon})). $$
    Assume $d_h([w],[z]) \leq \epsilon$. Let $\gamma$ be the geodesic between $[w],[z]$. Then $(\gamma,h(z) - h \circ \gamma)$ is a path $R_h \times [-\epsilon,\epsilon]$, where $\bar{h}^\epsilon$ takes constant value $f(z)$, impliying that $[w,f(z)-f(w)]=[z,0]$ in $R_h^\epsilon$. Now, given a loop $\beta$ in $\Rh$, by Lemma \ref{lem:pathLifting}, there exists a loop $\alpha$ in $Z$ such that $d_h(\pi_h \circ \alpha,\beta) \leq \epsilon$. By the argument above, $[\alpha,\bar{h}(\beta)-h \circ \alpha]=[\beta,0]$ in $\Rh^\epsilon$. Hence $\pi_{h}^\epsilon \circ (\alpha, \bar{h}(\beta)- h \circ \alpha) = [\beta,0] $ in $R_h^\epsilon$. This completes the proof.
\end{proof}

We are ready to prove Theorem \ref{thm:structure}.

\begin{proof}[Proof of Theorem \ref{thm:structure}]
    By idempotency of the Reeb graph construction, Corollary \ref{cor:idempotencyAndProduct}, the fibers of $\barf: \Rf \to \R$ are totally disconnected (i.e. connected components are singletons). By \cite[II.4.A]{hurewicz2015dimension}, fibers of $\barf$ are $0$-dimensional. Then, by \cite[Theorem~VI~7]{hurewicz2015dimension}, $\dim(R_f) \leq 1$.

    Given $\epsilon>0$, let $G$ be a finite metric graph such that $\dgh(X,G) < \epsilon$, whose existence is given by \cite[Proposition~7.5.5]{bbi01}. Let $R$ be a correspondence between $X$ and $G$ such that $\dis(R)<2 \epsilon$. Pick a $1$-dimensional CW complex structure on $G$ with vertex set $V={v_1,\dots,v_n}$ such that each edge has length at most $\epsilon$. Let $x_1,\dots,x_n$ be points in $X$ such that $(x_i,v_i) \in R$. Let $g: G \to \R$ be the function obtained by letting $g(v_i)=f(x_i)$, and extend $g$ to edges piecewise linearly. Let $(x,y) \in R$. There exist vertices $v_i,v_j$ such that $y$ is contained in an edge between $v_i,v_j$. Then, we have
    \begin{equation*}
        \begin{split}
            |f(x)-g(y)| &\leq |f(x)-f(x_i)|+|f(x_i)-g(v_i)|+|g(v_i)-g(y)| \\
                        &\leq d_X(x,x_i)+|g(v_i)-g(v_j)| \\
                        &< d_G(y,v_i)+2\epsilon + |f(x_i)-f(x_j)| \\
                        &\leq 3\epsilon+d_X(x_i,x_j) \\
                        &\leq 5\epsilon+d_G(v_i,v_j) \leq 6\epsilon.
        \end{split}
    \end{equation*}
    Hence, $\dis_{f,g}(R) < 6\epsilon$, implying $\dgh((X,f),(G,g)) < 3\epsilon$. Note that $g:G \to \R$ is piecewise linear, therefore $g^s: G \times [-s,s] \to \R$, $(x,t) \mapsto g(x)+t$ is also piecewise linear, implying that $\Rg^s$ is a finite graph by \cite[Proposition~7]{di2017reeb} for all $s \geq 0$. By Proposition \ref{prop:smoothingBetti}, $H :=\Rg^{36\epsilon}$ is a finite metric graph with $\betti(H) \leq \betti(X)$. Furthermore, by Proposition \ref{prop:approximation_smoothing}, we have
    $$\dgh(\Rf,H) \leq 72\epsilon + \dis(\pi_f^{0,\,72\epsilon})/2. $$
    By Proposition \ref{prop:rightContinuity}, the quantity above converges to $0$ as $\epsilon \to 0$. Hence we can construct a sequence $(G_n)$ of finite metric graphs with $\betti(G_n) \leq \betti(X)$ such that $(G_n)$ Gromov-Hausdorff converges to $\Rf$.

    By Lemma \ref{lem:simple_core}, there exists a core $A_n$ of $G_n$ such that $A_n$ is finite metric graph with number of edges bounded above by $3\betti(X)$. Since $A_n$ is homotopy equivalent to $G_n$ by Proposition \ref{prop:deformation_retraction}, $\betti(A_n) \leq \betti(X)$. By passing to a subsequence, by Lemma \ref{lem:graph_convergence}, we can assume that $A_n$ Gromov-Hausdorff converges to
    a finite metric graph
    $A$ such that $\betti(A) \leq \betti(X)$.
    Let $R_n$ be a correspondence between $G_n$  and $\Rf$ such that $\dis(R_n) \leq 2\dgh(\Rf,G_n)+2/n$. Let $B_n$ be the closure of the set of points in $\Rf$ corresponding to $A_n$ under $R_n$. Note that $\dgh(B_n,A_n) \leq \dgh(\Rf,G_n)+1/n$. By \cite[Theorem~7.3.8]{bbi01}, by passing to a subsequence if necessary, we can assume that $B_n$ Hausdorff converges to a closed subspace $B$ of $\Rf$. Note that $\dgh(A,B) \leq \dgh(A,A_n)+\dgh(A_n,B_n)+d_{\mathrm{H}}(B,B_n)$, which converges to $0$ as $n \to \infty$, hence $A$ is isometric to $B$. It remains to show that $B$ is a core of $\Rf$.

    By Corollary \ref{cor:quotient_hausdorff}, $\Rf/{B_n}$ Gromov-Hausdorff converges to $\Rf/B$.
        As in Appendix~\ref{app:quotient}, we use the notation $D_{\Rf}(x,B_n)$
        to denote the distance from $x\in \Rf$ to $B_n\subseteq \Rf$ in $\Rf$, and similarly use $D_{G_n}(y,A_n)$,
        for $y\in G_n$ and $A_n\subseteq G_n$.
    Let $f_n: \Rf \to \R$, $x \mapsto D_{\Rf}(x,B_n)$, and $g_n: G_n \to \R$, $y \mapsto D_{G_n}(y,A_n)$. Let $(x,y) \in R_n$. Let $a \in A_n$. There exists $b \in B_n$ such that $(a,b) \in R_n$. Hence, we have
     $$g_n(y)-\df(x,b) \leq d_{G_n}(y,a)-\df(x,b) \leq 2\dgh(\Rf,G_n)+2/n. $$
    Taking supremum of this inequality over $b \in B$, we get
    $$g_n(y)-f_n(x) \leq 2\dgh(\Rf,G_n)+2/n. $$
    Similarly,
    $$f_n(x)-g_n(y) \leq 2\dgh(\Rf,G_n)+2/n. $$
    This shows that $\dgh((\Rf,f_n),(G_n,g_n)) \leq  \dgh(X,G_n)+1/n.$ Now, by Lemma \ref{lem:quotient_gh} and Corollary \ref{cor:quotient_hausdorff}, we have 
    \begin{equation*}
        \begin{split}
            \dgh(\Rf/B, G_n/A_n) &\leq \dgh(\Rf/B,\Rf/B_n) + \dgh(\Rf/B_n,G_n/A_n) \\
                                 &\leq d_\mathrm{H}(B,B_n)+2\dgh(\Rf,G_n)+2/n,
        \end{split}
    \end{equation*}
    which converges to $0$ as $n$ goes to infinity. By Lemma \ref{lem:limit_tree}, $\Rf/B$ is a metric tree. This shows that $B$ is a core of $\Rf$ and concludes the proof.
\end{proof}

\section{Gromov Hausdorff Stability and Approximation of Reeb metric graphs}\label{sec:stability}

In this section, $(X,d_X)$ and $(Y,d_Y)$ are compact geodesic spaces,
$p \in X$, $q \in Y$, and $f: X \to \R$, $g: Y \to \R$ are
$1$-Lipschitz functions.

    Our main focus will be on the metric Reeb graph of the distance function $x\to d_X(p, x)$ to a fixed point $p\in X$.
   To ease notation we will use the shorthand $\Rp$ for the metric
   Reeb graph $R_{d_X(p,\cdot)}$ and similarly $d_p$ for the Reeb metric on $\Rp$.

The main theorems we prove in this section are the following:

\begin{theorem}[Stability]\label{thm:stability}
    Let $X$ and $Y$ be compact geodesic spaces, $f: X \to \R$ and $g: Y \to \R$ be $1$-Lispchitz functions. Let $M_f,M_g$ be the number of local minima of $f,g$ respectively. Then,
    $$ \dgh(\Rf,\Rg) \leq (8\betti(X)+8\betti(Y)+8M_f+8M_g-4)\,\dgh((X,f),(Y,g)). $$
\end{theorem}

    The following result shows that the distance between a finite metric graph $G$ and $X$ is lower bounded by the distance between $\Rp$ and $X$, up to a multiplicative constant depending on the Betti numbers of $X$ and $G$.

\begin{theorem}[Approximation]\label{thm:approximation}
    Let $G$ be a finite metric graph. Then,
    $$\dgh(X,\Rp) \leq (8\betti(X)+8\betti(G)+13)\,\dgh(X,G). $$
\end{theorem}

We postpone the proofs of these two theorems to the end of the section. First we need to prove some lemmas.

\begin{lemma}\label{lem:functional-gh}
    Let $X$ and $Y$ be compact metric spaces, $f: X \to \R$ and $g: Y \to \R$ be $1$-Lispchitz functions. Then $\dgh((X,f),(Y,g))=0$ if and only if there is an isometry $\psi: X \to Y$ such that $f=g \circ \psi$. Note that in this case, $\Rf$ is isometric to $\Rg$.
\end{lemma}
\begin{proof}

        The ``if'' part is straightforward as the isometry $\psi$ induces a correspondence between $X$ and $Y$ with distortion $0$ and the equation $f= g \circ \psi$ implies that $\dgh((X,f),(Y,g))=0$.
        Hence $\dgh((X,f),(Y,g))=0$.
    
    Now we assume that $\dgh((X,f),(Y,g))=0$.
    Let $R_n \subseteq X \times Y$ be a correspondence between $X$ and $Y$ such that $\dis_{f,g}R_n < 1/n$. Without loss of generality we can assume that $R_n$ is closed, as it has the same distortion with its closure. By \cite[Theorem~7.38]{bbi01}, by passing to a subsequence, we can assume that $R_n$ Hausdorff converges to a closed subset $R \subseteq X \times Y$ equipped with the sup-metric. Note that it is enough to show that $R$ is a correspondence with $\dis_{f,g}(R)=0$, hence in this case $R$ is the graph of an isometry $\psi: X \to Y$ and $f=g \circ \psi$. Let $\epsilon_n:=d_\mathrm{H}(R,R_n)$. Let $x \in X$. We need to show there exists $y \in Y$ such that $(x,y) \in R$. Let $y_n \in Y$ such that $(x,y_n) \in R_n$. There exists $(x_n,y_n')$ in $R$ such that $d_X(x,x_n), d_Y(y_n,y_n') \leq \epsilon_n$. Without loss of generality, we can assume that $y_n'$ converges to a point $y$. Hence, $(x_n,y_n')$ converges to $(x,y) \in R$. Similarly, for any $y \in Y$, there exists $x \in X$ such that $(x,y) \in R$. This shows that $R$ is a correspondence. It remains to show that $\dis_{f,g}(R)=0$.

    Let $(x,y), (x',y') \in R$. There exists $(x_n,y_n), (x_n', y_n')$ in $R_n$, which are $\epsilon_n$ close to $(x,y), (x',y')$ respectively. We have,
    \begin{align*}
        \begin{split}
            |f(x)-g(y)| & \leq |f(x)-f(x_n)|+|f(x_n)-g(y_n)|+|g(y_n)-g(y)| \leq 2\epsilon_n+1/n
        \end{split} 
    \end{align*}
    Hence, letting $n \to \infty$, we get $|f(x)-g(y)|=0$. We also have,
    \begin{align*}
        \begin{split}
            |d_X(x,x')-d_Y(y,y')| &\leq |d_X(x,x')-d_X(x_n,x_n')|+|d_X(x_n,x_n')-d_Y(y_n,y_n')| \\ &+|d_Y(y_n,y_n')-d_Y(y,y')| \\
            &\leq d_X(x,x_n)+d_X(x',x_n')+1/n+d_Y(y,y_n)+d_Y(y',y_n')\\
            & \leq 4\epsilon_n+1/n.
        \end{split}
    \end{align*}
    Letting $n \to \infty$, we get $|d_X(x,x')-d_Y(y,y')|=0$. This shows that $\dis_{f,g}(R)=0$, and the proof is complete.
\end{proof}

\begin{lemma}\label{lem:graph_paths}
    Assume that $\Rf$ has a core which is a finite metric graph and $\barf: \Rf \to \R$ has finitely many local minima. Let $M$ be the number of local minima of $\barf$. Then for each path $\beta: [0,1] \to \Rf$ from $[x]$ to $[x']$, there exists $0=t_0 \leq t_1 \leq  \dots  \leq t_m=1$ where $m \leq 2(M+\betti(\Rf))$ such that
    $$\df([x],[x']) \leq \Sigma_{i=1}^m|\barf(\beta(t_i))-\barf(\beta(t_{i-1}))|. $$
    
        Furthermore, if $\beta$ is a geodesic, then the equality holds.
    
\end{lemma}
\begin{proof}
    By Lemma \ref{lem:graph_core_extension}, we can assume that the finite metric graph core $G$ of $\Rf$ contains $[x],[x']$. By Lemma \ref{lem:core_path}, without loss of generality we can assume that $\beta$ is contained in $G$. By Lemma \ref{lem:local_minima}, number of local minima of $\barf|_G$ is less than or equal to $M$. By Lemma \ref{lem:edge_injective}, there is a $1$-dimensional CW-complex structure on $G$ such that $\barf$ is injective over each edge. Add $[x],[x']$ to the vertex set. In the remaining part of the proof we will use this CW-complex structure on $G$.
    
    Without loss of generality, we can assume that $\beta$ does not visit $[x]$ twice since we can just consider the part of $\beta$ after its last visit of $[x]$. Define a finite sequence of points in $G$ inductively as follows: $[x_0]=[x]$, $[x_{i+1}]$ is the first vertex visited by $\beta$ after its last visit of $[x_i]$, stop when $[x_n]=[x']$. Let $s_n$ be the last time $\beta$ visits $[x_i]$. We have $0=s_0 < s_1 < \dots s_n =1 $. Note that the $[x_i]$s are all distinct, as $\beta$ visits $[x_{i+1}]$ after its last visits of $[x_0],\dots,[x_i]$. Furthermore, there is an edge between $[x_i],[x_{i+1}]$, since otherwise $\beta$ has to visit other vertices after its last visit of $[x_i]$ before reaching $[x_{i+1}]$. Since $\barf$ is $1$-Lipschitz, $\df([x_i],[x_{i-1}]) \geq |f(x_i)-f(x_{i-1})|$. By using the edge between $x_{i-1},x_i$, we see that $\df([x_i,x_{i-1}]) \leq |f(x_i)-f(x_{i-1})|$. Hence $\df([x_i],[x_{i-1}])=|f(x_i)-f(x_{i-1})|$. Now, we have
    $$\df([x],[x']) \leq \sum_{i=1}^n \df([x_i],[x_{i-1}]) = \sum_{i=1}^n |f(x_i)-f(x_{i-1})|.$$

    Consider $G$ as a directed graph, where an edge is directed from the vertex taking the smaller $\barf$ value to the higher one, and we arbitrarily choose an orientation for the edges between vertices with the same $\barf$ value.

             Let $i_1,\dots,i_N$ be the indices of the vertices in the sequence $[x_1],\dots,[x_{n-1}]$
            whose indegree is greater than $1$ sorted in increasing order. We then extend them to  a sequence of indices $S := \{i_0=0, i_1,\dots,i_N,i_{N+1}=n\}$.
            Then, for any two consecutive indices $i_k$ and $i_{k+1}$ in $S$, the function values $f([x_{i_k}]), f([x_{i_k + 1}]), \ldots, f([x_{i_k}])$ on the vertices corresponding to the consecutive integers $i_k, i_k + 1, \ldots, i_{k+1}$ between $i_k$ and $i_{k+1}$
            is either monotonic or initially decreasing and then increasing with respect to these indices. This is because otherwise there would be a vertex with indegree greater than $1$ and whose index is between $i_k$ and $i_{k+1}$.
            Let $j_k$ be the index in $S$ where the function value on the corresponding vertex achieves a minimum among all indices in $S$.

    We have
    $$\Sigma_{i=i_{k}+1}^{i_{k+1}} |f(x_i)-f(x_{i-1})|=f(x_{i_k})-f(x_{j_k}) + f(x_{i_{k+1}})-f(x_{j_k}). $$
    Let $t_{2k}=s_{i_k}, t_{2k+1}=s_{j_k}$. Note that $0=t_0 \leq t_1 \leq \dots \leq t_{2N+2}=1$, and
    $$\Sigma_{i=1}^n |f(x_i)-f(x_{i-1})|= \Sigma_{i=1}^{2N+2} |\barf(\beta(t_i)-\barf(\beta(t_{i-1}))|.$$

    It remains to show that $N \leq \betti(\Rf)+M-1$. Note that by Proposition \ref{prop:deformation_retraction}, $\betti(G)=\betti(\Rf)$. Let us denote the indegree of a vertex $v$ by $\iota(v)$. Since the number of edges is equal to the sum of indegrees, by Euler's formula
    applied to $G$, we have
    $$\betti(G)-1= \Sigma_{v} \iota(v) - \Sigma_{v} 1, $$
    where the summation is taken over the vertices of $G$. The vertices with $\iota(v)=0$ are exactly the local minima of $\barf|_G$. Also, $N$ is less than equal to the number of vertices with $\iota(v)>1$. Hence, we have
    $$\betti(\Rf)-1=\betti(G)-1= \Sigma_{v: \iota(v)>1} (\iota(v)-1) - M \geq N-M, $$
    which implies that $N \leq \betti(\Rf)+M-1$.
    This completes the proof of the first part of the lemma.
    For the second part, note that the right-hand side of the
      inequality is always less than or equal to the length of
      $\beta$. Hence, if $\beta$ is a geodesic realizing
      $\df([x],[x'])$, the inequality is actually an  equality.
\end{proof}

\begin{lemma}\label{lem:isometry}
    Let $G$ be a finite metric graph and $q \in G$. Then the Reeb quotient map $G \to \Rq$ is an isometry.
\end{lemma}
\begin{proof}
    By Lemma \ref{lem:distance_to_point}, there is a $1$-dimensional CW complex structure on $G$ such that $g:=d_G(q,\cdot): G \to \R$ maps each edge isometrically into $\R$. This implies the level sets of $g$ are finite, as each edge contains at most one point in a level. Hence the Reeb quotient map is the identity. If $x,x'$ are contained in an edge in $G$, then by above $d_G(x,x')=d_G(p,x)-d_G(p,x') \leq \dg(x,x') \leq d_G(x,x')$, hence $d_G(x,x')=\dg(x,x')$. Since the metric structure of a finite metric graph is completely determined by the metric structure of its edges \cite[Exercise~3.2.12]{bbi01}, $d_G = \dg$.
\end{proof}

We are ready to prove Theorem \ref{thm:stability} and Theorem \ref{thm:approximation}.

\begin{proof}[Proof of Theorem \ref{thm:stability}]
    By Lemma \ref{lem:functional-gh}, we can assume that $\betti(X),\betti(Y),M_f,M_g < \infty$. By Theorem \ref{thm:structure}, $\Rf,\Rg$ have cores which are finite metric graphs with $\betti(\Rf) \leq \betti(X)$, $\betti(\Rg) \leq \betti(Y)$. If $x$ is not a local minimum of $f$, then $[x]$ is not a local minimum of $\barf$.
    Hence the number of local minima of $\barf$ is less than or equal to $M_f$. Similarly, the number of local minima of $\barg$ is less than or equal to $M_g$.

    Assume $r>2\dgh((X,f),(Y,g))$. Let $R$ be a correspondence between $X$ and $Y$ with distortion $\dis_{f,g}(R) < r$. Let $\epsilon>0$. Let $(x,y),(x',y') \in R$. Let $\beta$ be a geodesic from $[x]$ to $[x']$. By Lemma \ref{lem:pathLifting}, there exists a path $\alpha$ in $X$ from $x$ to $x'$ such that $|\barf\circ \beta - f \circ \alpha| < \epsilon$. Let $0=s_0 < \dots < s_n=1$ be such that $\diam(\alpha|_{[s_{i-1},s_i]}) \leq \epsilon$ for all $i=1,\dots,n$. Let $x_i=\alpha(s_i)$ for all $i$. Let $(y_i)_{i=0}^n$ be points in $Y$ such that $y_0=y$, $y_n=y'$ and $(x_i,y_i) \in R$ for all $i$. Let $\gamma$ be a path in $Y$ such that $\gamma|_{[s_{i-1},s_i]}$ is a geodesic from $y_{i-1}$ to $y_i$ for $i=1,\dots,n$. Such a path can always be constructed by concatenating said geodesics.
    
    Let $t \in [0,1]$, and $i$ be an index such that $t \in [s_{i-1},s_i]$. We have
    \begin{equation*}
        \begin{split}
            |\barf(\beta(t))-g(\gamma(t))| &\leq |\barf(\beta(t))-f(\alpha(t)|+|f(\alpha(t))-f(x_i)|+|f(x_i)-g(y_i)|+|g(y_i)-g(\gamma(t))| \\
            & \leq 2\epsilon + r + d_Y(y_{i-1},y_i) \leq 2r + 3\epsilon.
        \end{split}
    \end{equation*}
    Note that the inequality in Lemma \ref{lem:graph_paths} still holds if we consider any refinement of the partition of $[0,1]$.
    
        Let $0=:s_0 \leq  \dots \leq s_k:=1$ and $0=:u_0 \leq \dots \leq
        u_l:=1$ be the partitions of $[0,1]$ for $\beta$, and $\pi_g
        \circ \gamma$ respectively, given by Lemma
        \ref{lem:graph_paths}. Taking the common refinement $(t_i)$ of these partitions without repeating $0,1$, we get $0=:t_0 \leq \dots \leq t_m:=1$, where $m:=k+l-1 \leq C$ where $C:= 2\betti(X)+2\betti(Y)+2M_f+2M_g-1$. Since the partition $(t_i)$ is a refinement of both $(s_i)$ and $(u_i)$, the inequality given in Lemma \ref{lem:graph_paths} holds for both $\beta$, and $\pi_g \circ \gamma$ with the partition $(t_i)$. Then we have
        \begin{align*}
            d_g([y], [y']) \leq \Sigma_{i=1}^m |g(\gamma(t_i))-g(\gamma(t_{i-1}))|
        \end{align*}
        and
        \begin{align*}
            d_f(\alpha(t),\beta(t)) &= \Sigma_{i=1}^m |\barf(\beta(t_i)-\barf(\beta(t_{i-1}))|
        \end{align*}
        where the equality holds in the second formula because $\beta$ is a geodesic.

    Therefore, we have
    \begin{equation*}
        \begin{split}
            \dg([y],[y'])-\df([x],[x']) &\leq \Sigma_{i=1}^m \left(|g(\gamma(t_i))-g(\gamma(t_{i-1}))|-|\barf(\beta(t_i)-\barf(\beta(t_{i-1}))| \right)\\
            &\leq \Sigma_{i=1}^m \left(|\barf(\beta(t_i))-g(\gamma(t_{i-1}))|+|\barf(\beta(t_{i-1})-g(\gamma(t_{i-1}))| \right) \\
            &\leq 4Cr+6C\epsilon.
        \end{split}
    \end{equation*}
    Since $\epsilon>0$ was arbitrary, we have $\dg([y],[y'])-\df([x],[x']) \leq 4Cr.$ Similarly, starting from a geodesic between $[y],[y']$, we can show that $\df([x],[x'])-\dg([y],[y']) \leq 4Cr$. Hence, if we consider the correspondence $\bar{R}$ between $\Rf$ and $\Rg$ given by $([x],[y]) \in \bar{R}$ if $(x,y) \in R$. then $\dgh(\Rf,\Rg) \leq \dis(\bar{R})/2 \leq 2Cr $. This completes the proof, since $r>2\dgh((X,f),(Y,g))$ was arbitrary.
\end{proof}

\begin{remark}
    Note that, in the statement and proof of Theorem \ref{thm:stability}, we could use $\betti(\Rf),\betti(\Rg), M_{\barf},M_{\barg}$ instead of $\betti(X),\betti(Y),M_f,M_g$, which are always respectively smaller, by Theorem \ref{thm:structure}.
\end{remark}

\begin{proof}[Proof of Theorem \ref{thm:approximation}]
    Let $r>\dgh(X,G)$. Let $R$ be a correspondence between $X$ and $G$ such that $\dis(R) < r$. Let $q \in G$ such that $(p,q) \in R$. Note that $\dgh((X,p),(G,q)) < r$. Let $f:X \to \R$, $x \mapsto d_X(p,x)$, and $g: Y \to \R$, $y \mapsto d_Y(p,y)$. Both $f$ and $g$ have uniqe local minimums, $p,q$ respectively.  By Theorem \ref{thm:stability} and Lemma \ref{lem:isometry}, we have
    \begin{equation*}
        \begin{split}
            \dgh(\Rp,G)=\dgh(\Rp,\Rq) &\leq (8\betti(X)+8\betti(G)+12)\,\dgh((X,p),(G,q)) \\
            &\leq (8\betti(X)+8\betti(G)+12)r.
        \end{split}
    \end{equation*}
    Since $r>\dgh(X,G)$ was arbitrary, we get
    $$\dgh(\Rp,G) \leq (8\betti(X)+8\betti(G)+12)\,\dgh(X,G). $$
    Therefore,
    $$\dgh(X,\Rp) \leq \dgh(X,G)+\dgh(\Rp,G) \leq (8\betti(X)+8\betti(G)+13)\,\dgh(X,G). $$
\end{proof}

\section{Merge Metric trees}

In this section, we define \emph{merge metric trees} and establish some of their fundamental properties. In the next section, we are going to use merge metric trees to obtain metric tree approximations of geodesic spaces.

\begin{definition}[Merge function]
    Let $X$ be a compact connected topological space and $f: X \to \R$ be a continuous function. Define $m_f: X \times X \to \R$, the merge function induced by $f$, as
    $$m_f(x,x') :=\inf\{\sup(f|_A): A \subseteq X, \, x, x' \in A, \, A \text { is connnected } \}. $$
\end{definition}

Before defining merge metric trees, we need to establish some properties of the merge function.

\begin{proposition}\label{prop:merge_function}
    Let $X$ be a compact connected topological space, $f: X \to \R$
    continuous, and $x,x',x'' \in X$. Then, 
    \begin{enumerate}[i)]
        \item $f(x),f(x') \leq m_f(x,x')=m_f(x',x)$.
        \item $m_f(x,x'') \leq \max(m_f(x,x'),m_f(x',x'')). $
        \item The map ${\mathrm{t}_f}: X \times X \to \R$, $(x,x') \mapsto 2m_f(x,x')-f(x)-f(x')$ defines a pseudometric on $X$.
    \end{enumerate}
\end{proposition}
\begin{proof}
    For i) The claim follows directly from the definition of $m_f$.
    For ii) This can be deduced by taking the union of connected sets used in the definition of $m_f(x,x')$ and $m_f(x',x'')$.
    For iii) By definition, it is clear that ${\mathrm{t}_f}(x,x)=0$ and ${\mathrm{t}_f}(x,x')={\mathrm{t}_f}(x',x) \geq 0$. We need to show that ${\mathrm{t}_f}$ satisfies the triangle inequality. We have
    \begin{equation*}
        \begin{split}
            {\mathrm{t}_f}(x,x')+{\mathrm{t}_f}(x',x'') &= 2m_f(x,x')+2m_f(x',x'')-f(x)-2f(x')-f(x'') \\
                   &= 2\max(m_f(x,x'),m_f(x',x''))-f(x)-f(x'') + 2(\min(m_f(x,x'),m_f(x',x''))-f(x'))  \\
                   &\geq 2m_f(x,x'')-f(x)-f(x'') = {\mathrm{t}_f}(x,x'').
        \end{split}
    \end{equation*}
\end{proof}

\begin{definition}[Merge metric tree]
    Let $X$ be a compact connected topological space and $f: X \to \R$ be a continuous function. The \emph{merge metric tree} induced by $f:X \to \R$, denoted by $(T_f,\mathrm{t}_f)$, is the metric space induced by the pseudometric
    $$\mathrm{t}_f:X \times X \to \R, \,(x,x') \mapsto 2m_f(x,x')-f(x)-f(x').$$
    In other words, $T_f = X/\sim$, where $x \sim x'$ if and only if ${\mathrm{t}_f}(x,x')=0$ and $T_f$ is endowed with the quotient metric $\mathrm{t}_f$.
\end{definition}
    Note that one difference between $\mathrm{t}_f$ and $\mathrm{d}_f$ is that $\mathrm{t}_f$ can be defined on general compact connected spaces, while $\mathrm{d}_f$ requires path connectedness.
    With a slight abuse of notation, as in the case of Reeb graphs, we
    denote the point in $T_f$ induced by a point $x$ in $X$ by $[x]$
    and by $\barf: T_f \to \R$ we mean the function $[x] \mapsto
    f(x)$. This should not cause any confusion since we never invoke
    Reeb graphs and merge trees in the same statement. We define
    $\tau_f: X \to T_f$, $x \mapsto [x]$, which a priori is not
    necessarily continuous. Since $m_f(x,x') \geq f(x),f(x')$, we have
    $\mathrm{t}_f([x],[x']) \geq f(x)-f(x')$, hence $\barf: T_f \to
    \R$, $[x] \mapsto f(x)$ is well defined and $1$-Lipschitz. If $X$
    is a geodesic space and $f: X \to \R$ is $1$-Lipschitz, by
    considering the geodesic between $x$ and $x'$ in $X$, one sees
    that $\mathrm{t}_f([x],[x']) \leq d_X(x,x')$, hence $\tau_f: X \to
    T_f$ is a $1$-Lipschitz quotient map in this case. Also, in this case, $T_f$ is connected.

\medskip

Let us show that $T_f$ is a metric tree (cf. Appendix \ref{app:hyperbolicity}).

\begin{proposition}\label{prop:zero_hyp}
    Let $X$ be a compact connected topological space and $f: X \to \R$ be a continuous function. Then, $(T_f,\mathrm{t}_f)$ has hyperbolicity zero. If $X$ is a geodesic space, then $T_f$ is a metric tree.
\end{proposition}
\begin{proof}
    Let $p$ be a point where $f$ achieves its maximum. Then $m_f(p,x)=f(p)$ for any $x$ in $X$. So, $\mathrm{t}_f([p],[x])=f(p)-f(x)$. Considering Gromov product $g_{[p]}$, by Proposition \ref{prop:merge_function} we have
    $$g_{[p]}([x],[x'])=(\mathrm{t}_f([p],[x])+\mathrm{t}_f([p],[x'])-\mathrm{t}_f([x],[x'])/2=f(p)-m_f(x,x'). $$
    Hence $g_{[p]}(x,x'') \geq \min(f(p)-m_f(x,x'),f(p)-m_f(x',x''))= \min(g_{[p]}(x,x'),g_{[p]}(x',x''))$ for any $x,x',x''$ in $X$. Therefore, $\hyp_{[p]}(T_f)=0$, which implies that $\hyp(T_f)=0$.

    If $X$ is geodesic, then $T_f$ is connected as $\tau_f: X \to T_f$ is continuous and surjective. By Lemma \ref{lem:tree_connected}, $T_f$ is a metric tree.  
\end{proof}

The following result shows which points are identified under the quotient map $\tau_f: X \to T_f$.

\begin{proposition}\label{prop:merge_quotient}
    Let $x,x' \in X$. Then, $\mathrm{t}_f([x],[x'])=0$ if and only if $f(x)=f(x')=:c$ and $x,x'$ are in the same connected component of $\{f \leq c \}$.
\end{proposition}
\begin{proof}
    The if part follows from the definition of $\mathrm{t}_f$. We can assume that $f(x)=f(x')=c$, since otherwise $\mathrm{t}_f([x],[x']) \neq 0$. Note that $\mathrm{t}_f([x],[x'])=0$ if and only if $m_f(x,x')=c$. Let $C_n$ be the connected component of $x$ in $\{f \leq c+1/n \}$. Note that $y \in C_n$. $(C_n)_n$ forms a decreasing family of compact connected sets in $X$. By \cite[Corollary~6.1.19]{engelking1989general}, $C:=\cap_n C_n$ is a connected set containing $x,x'$. Since $C \subseteq \{f \leq c \}$, $x$ and $x'$ are in the same connected component of $\{ f \leq c \}$.
\end{proof}

As in the case of Reeb graphs, we establish  idempotency and naturality of the merge metric tree construction.

\begin{theorem}[Merge metric tree naturality]\label{thm:merge_tree_naturality}
    Let $X, Y$ be compact, connected topological spaces and $f: X \to \R$, $g: Y \to \R$ be continuous functions. Let $\phi: X \to Y$ be continuous maps such that $f = g \circ \phi $. Then $\bar{\phi}: T_f \to T_g$, $[x] \mapsto [\phi(x)]$ is $1$-Lipschitz.
    {If additionally, $X,Y$ are Hausdorff and $\phi$ is onto and has connected fibers, then $\bar{\phi}$ is an isometry.}
\end{theorem}
\begin{proof}
    If $x,x' \in X$ and $A$ is a connected subset of $X$ containing $x$ and $x'$, then $\phi(A)$ is a connected subset of $Y$ containing $\phi(x),\phi(y)$, and $\sup f|_A=\sup g|_{\phi(A)}$. This shows that $m_g(\phi(x),\phi(x')) \leq m_f(x,x')$, implying $\mathrm{t}_g([\phi(x)],[\phi(x')]) \leq \mathrm{t}_f([x],[x'])$.
    
    For the second claim, it is enough to show that $m_g(\phi(x),\phi(x')) \leq m_f(x,x')$. Let $B$ be a connected subset of $Y$ containing $\phi(x),\phi(x')$.
  
        The assumptions that $X, Y$ are Haudorff and $\phi$ is onto and has connected fibers allow us to apply \cite[Theorem~6.1.28]{engelking1978dimension} to conclude that $C=\phi^{-1}(B)$ is a connected subset of $X$ containing $x, x'$.
    
    Note that $\sup f|_C=\sup g|_A$. Hence $m_f(x,x') \leq m_g(\phi(x),\phi(x'))$.
\end{proof}

Up to this point, we only established that $T_f$ is a metric tree when $X$ is a geodesic space and $f: X \to \R$ is $1$-Lipschitz. The result above allows us to extend this to a more general setting.

\begin{corollary}
    Let $E$ be a compact, connected, Hausdorff space and $\phi: E \to \R$ be a continuous function. Assume there is a geodesic space $X$ with a $1$-Lipschitz function $f: X \to \R$, and a continuous map $\pi: E \to X$ with connected fibers. Then $(T_\phi,\mathrm{t}_\phi)$ is a metric graph and $\tau_\phi: E \to T_\phi$ is continuous.
\end{corollary}
\begin{proof}
    Consider the following commutative diagram:
    $$\begin{tikzcd}
        E \arrow[r, "\pi"] \arrow[d,"\tau_\phi"] & X \arrow[d, "\tau_f"]\\
        T_\phi \arrow[r, "\bar{\pi}"] & T_f
    \end{tikzcd}$$
    By Theorem \ref{thm:merge_tree_naturality}, $\bar{\pi}$ is an isometry, hence by Proposition \ref{prop:zero_hyp}, $T_\phi$ is a metric tree. $\tau_\phi$ is continuous since it is equal to $\bar{\pi}^{-1} \circ \tau_f \circ \pi$.
\end{proof}

\begin{theorem}[merge metric tree idempotency]\label{thm:merge_tree_idempotency}
    Assume $X$ is a geodesic space and $f: X \to \R$ is $1$-Lipschitz. Let $\barf: T_f \to \R$. $[x] \mapsto f(x)$. Then $T_f \to T_{\barf}$, $[x] \mapsto [[x]]$ is an isometry.
\end{theorem}

We need the following lemma:
\begin{lemma}\label{lem:sublevel_component}
    Assume $X$ is a geodesic space and $f: X \to \mathbb{R}$ is $1$-Lipschitz. Let $x,{x'} \in X$. Then $x,{x'}$ are contained in the same connected component of $\{ f \leq r \}$ in $X$ if and only if $[x], [{x'}]$ are contained in the same connected component of $ \{ \bar{f} \leq r\}$ in $T_f$.
\end{lemma}

\begin{proof}
    The ``only if'' part follows from the fact that the map $\tau_f: X \to T_f$ is continuous, and hence the image of the connected component containing $x$ and ${x'}$ in $\{ f \leq r \}$ under $\tau_f$ is connected in $T_f$.

    Let $p \in X$ be a point where $f$ achieves its maximum and set $M:=f(p)$.
    We define $E$ to be the set
    $$E:=\{(x,t): x \in X,\, t \in \mathbb{R},\, f(x) \leq t \leq M \}.$$
    We equip $E$ with the $\ell_1$ product metric so that $E$ is then a compact metric space.
    We now define a map $\pi$ from $E$ to $T_f$ through the following construction.
    For each point $x$ in $X$, we fix $\gamma_x: [0, 1]\to X$, a path from $x$ to $p$.
    Given such $\gamma_x$, we define a function $\lambda_{\gamma_x}:=\mathbb{R}\rightarrow[0,1]$, such that for a given $t \in \mathbb{R}$, $\lambda_{\gamma_x}(t) = \inf\{s \in [0, 1]\mid f(\gamma_x(s))=t\}$.
    Essentially, $\lambda_{\gamma_x}(t)$ gives us the smallest 'time' parameter $s$ along the path $\gamma_x$ where the function $f$ takes the value $t$.
    Now, for any given $t \in \mathbb{R}$ with $f(x)\leq t \leq M$, we define $x_{t}$ as the point on the path $\gamma_x$ where $f$ equals $t$ for the first time. Using our newly introduced function, we express this formally as
    \[
        x_{t}:=\gamma_x(\lambda_{\gamma_x}(t)).
    \]
    \noindent\textbf{Claim:} the construction of $x_{t}$ is independent of the choice of $\gamma_x$.\\
    Indeed, let $\gamma'_x$ be another path from $x$ to $p$, and we obtain $x'_{t}$ in the same manner as $x_{t}$.
    Then, $x_{t}$ and $x'_{t}$ both satisfy that $f(x_{t})=f(x'_{t})=t$, and they are both in the same connected component of $\{f \leq t \}$.
    Therefore, Proposition~\ref{prop:merge_quotient} implies that $[x_{t}]=[x'_{t}]$.

\smallskip
    With these definitions in place, we define the map $\pi: E \to T_f$, where $(x,t) \mapsto [x_{t}]$ by fixing a path $\gamma_x$ for each $x \in X$.
    Our next step is to show that this map $\pi$ is $1$-Lipschitz.

    \noindent\textbf{Claim:} $\mathrm{t}_f([x_{s}],[x_{t}])=|s-t|$: Assume $s \leq t$.
    We consider a segment, $\alpha$, of the path $\gamma_x$ that goes from $x_{s}$ to $x_{t}$, which is formally represented by $\gamma_x|_{[\lambda_{\gamma_x}(s), \lambda_{\gamma_x}(t)]}$. It is clear that $t =f(x_{t}) \leq m_f(x,x_{t}) \leq \sup f\circ \alpha=t$. As a result, we find $m_f(x_{s},x_{t})=t$, and hence, by definition of the pseudometric $\mathrm{t}_f$, we have $\mathrm{t}_f([x_{s}],[x_{t}])=2t-t-s=t-s$.

\smallskip
    \noindent\textbf{Claim} $m_f(x_{s},{x'}_{s})=\max(s,m_f(x,{x'}))$: Since $x=x_{f(x)}$ and $s\geq f(x)$, by the above discussion, there is
    $m_f(x,x_{s})= m_f(x_{f(x)}, x_s)=s$.
    We have,
    \begin{equation*}
        \begin{split}
            m_f(x_{s},{x'}_{s}) & \leq \max(m_f(x_{s},x),m_f(x,{x'}),m_f({x'},{x'}_{s})) \\
                         & \leq \max(s,m_f(x,x_{s}),m_f(x_{s},{x'}_{s}),m_f({x'}_{s},{x'})) \\
                         & = \max(s,m_f(x_{s},{x'}_{s}))=m_f(x_{s},{x'}_{s}).
        \end{split}
    \end{equation*}
    This proves the claim.

\smallskip
    If $m_f(x,{x'}) \leq s$, then $m_f(x_{s},{x'}_{s})=s$, which implies that $\mathrm{t}_f([x_{s}],[{x'}_{s}])=0$. If $m_f(x,{x'}) > s$, then $m_f(x_{s},{x'}_{s})=m_f(x,{x'})$, hence $\mathrm{t}_f([x_{s}],[{x'}_{s}])=2m_f(x,{x'})-2s \geq \mathrm{t}_f([x],[{x'}])$. Therefore, in any case, $\mathrm{t}_f([x_{s}],[{x'}_{s}]) \leq \mathrm{t}_f([x],[{x'}])$. Combining this with the first claim, we get
    $$\mathrm{t}_f([x_{s}],[{x'}_{t}]) \leq \mathrm{t}_f([x_{s}],[x_{t}])+\mathrm{t}_f([x_{t}],[{x'}_{t}]) \leq |s-t|+ \mathrm{t}_f([x],[{x'}]) \leq |s-t|+d_X(x,{x'}).$$ Therefore, $\pi: E \to T_f$ is $1$-Lipschitz. It is onto since $x_{f(x)}=x$.

    For $(x,s) \in E$, let $A_x^s$ denote the connected component of $x$ in $\{f \leq s \}$. It's worth noting that $x$ and $x_s$ are contained in $A_x^s$.
    By Proposition \ref{prop:merge_quotient}, $[{x'}_{t}]=[x]$ if and only if $f(x)=t$ and ${x'}$ is in the connected component of $x_{t}$ in $\{f \leq t \}$. Therefore, $[{x'}_{t}]=[x]$ if and only ${x'} \in A_x^{f(x)}$. Hence  $$\pi^{-1}([x]) = A_x^{f(x)} \times \{f(x)\}.$$
    This shows that $\pi$ has connected fibers. By \cite[Theorem~6.1.29]{engelking1989general}, preimages of connected sets under $\pi$ are connected.

    Let $\phi: E \to \mathbb{R}$, $(x,t) \mapsto t$. Note that $\phi=\bar{f} \circ \pi$. Assume $[x],[{x'}]$ are contained in the same component $C$ of the sublevel set $\{\bar{f} \leq r \}$. Then $D=\pi^{-1}(C)$ is a connected set contained in the sublevel set $\{\phi \leq r \}$. Hence $x,{x'}$ are contained in the same connected component of $\{f \leq r\}$.
\end{proof}

\begin{proof}[Proof of Theorem \ref{thm:merge_tree_idempotency}]
    By taking the image of a connected set containing $A$ under $\tau_f$, one sees that $m_{\barf}([x],[{x'}]) \leq m_f(x,{x'})$. Assume $m_{\barf}([x],[{x'}]) <r$. Then $[x],[{x'}]$ are in the same connected component of $\{\barf \leq r \}$. By Lemma \ref{lem:sublevel_component}, $x,{x'}$ are in the same connected component of $\{f \leq r \}$. Hence $m_f(x,{x'}) \leq r$. Since $r>m_{\barf}([x],[{x'}])$ was arbitrary, we have $m_f(x,{x'}) \leq m_{\barf}([x],[{x'}])$. Therefore, $m_f(x,{x'})=m_{\barf}([x],[{x'}])$, which implies that $\mathrm{t}_f([x],[{x'}])=t_{\barf}([[x]],[[{x'}]])$.
\end{proof}

\section{Gromov-Hausdorff Stability of Merge Metric Trees and Tree Approximations}\label{sec:tree_approximation}

In this section, $(X,d_X)$ and $(Y,d_Y)$ are compact geodesic spaces, $p \in X$, $q \in Y$, and $f: X \to \R$, $g: Y \to \R$ are $1$-Lipschitz functions.

\begin{theorem}[Merge metric tree stability]\label{thm:merge_tree_stability}
    Let $(X,d_X)$ and $(Y,d_Y)$ be compact geodesic spaces, and $f: X \to \R$, $g: Y \to \R$ be $1$-Lipschitz functions. Then,
    $$\dgh(T_f,T_g) \leq 12\, \dgh((X,f),(Y,g)).$$
\end{theorem}
\begin{proof}
    Let $r>\dgh((X,f),(Y,g))$. Let $E$, $\pi_X: E \to X$, $\pi_Y$, $\phi: E \to \R$, $\psi: E \to \R$ be as in Lemma \ref{lem:correspondence_nbhd}. By Theorem \ref{thm:merge_tree_naturality}, $T_f$ is isometric to $T_\phi$ and $T_g$ is isometric to $T_\psi$.
    For any $e\in E$, we use $[e]_\phi$ to denote the point in $T_\phi$ induced by $e \in E$  and $[e]_\psi$ to denote the point in $T_\psi$ induced by $e \in E$.
    Since $|\phi-\psi|<6r$, then for $e,e' \in E$, $|m_\phi(e,e')-m_\psi(e,e')|<6r$.
   Therefore, we have $$|\mathrm{t}_\phi([e]_\phi,[e']_\phi)-\mathrm{t}_\psi([e]_\psi,[e']_\psi)| \leq 2|m_\phi(e,e')-m_\psi(e,e')|+|\phi(e)-\psi(e)|+|\phi(e')-\psi(e')| < 24r.$$
    Hence, if we consider the correspondence between $T_\phi,T_\psi$ induced by $E$, it has distortion less than $24r$. Therefore,
    $$\dgh(T_f,T_g)=\dgh(T_\phi,T_\psi) < 12r.$$
    This completes the proof, as $r>\dgh((X,f),(Y,g))$ was arbitrary.
\end{proof}

    \begin{remark}\label{rem:compare_perez}
        While \cite{perez2020c} recently investigated the stability of
        merge trees, their analysis is confined to functions $f,g:X\to
        \R$ sharing
        the same domain $X$. Notably, Theorem~4.21 in \cite{perez2020c}
        establishes that $\dgh(T_f, T_g)\leq 2
        \|f-g\|_{L^\infty(X)}$. Since $\dgh((X,f),(X,g))\leq
        \|f-g\|_{L^\infty(X)}$, their bound is tighter than ours. However, our approach offers greater generality by accommodating functions with distinct domains.
    \end{remark}

\noindent\textbf{Notation:} If $X$ is a compact geodesic space and $p \in X$, $(T_p,t_p)$ denotes the merge metric tree induced by the function $X \to \R$, $x \mapsto -d_X(p,x)$. We denote the merge metric tree quotient map by $\tau_p: X \to T_p$.

The following result shows that $T_p$ is the best tree approximation of $X$ up to a factor of $13$.
\begin{theorem}[Metric tree approximation]\label{thm:metric_tree_approximation}
    $(X,d_X)$ be a compact geodesic space and $p \in X$. Let $T$ be an arbitrary compact metric tree. Then,
    $$\dgh(X,T_p) \leq 13\, \dgh(X,T).$$
\end{theorem}

We first prove the following lemma:
\begin{lemma}\label{lem:Tq}
    Let $(T,d_T)$ be a compact metric tree and $q \in T$. Then $\tau_q: T \to T_q$ is an isometry.
\end{lemma}
\begin{proof}
    Let $x,x' \in T$. We already know that $t_q([x],[x']) \leq d_T(x,x')$. Let us show that $d_T(x,x') \leq t_q([x],[x'])$. Let $\gamma$ be the unique simple path from $x$ to $x'$. Let $x''$ be the point closest to $q$ on $\gamma$. Let $\gamma_1$ be the part of $\gamma$ from $x''$ to $x$, and $\gamma_2$ be the part of $\gamma$ from $x''$ to $x'$. Let $\alpha$ be the unique simple path from $q$ to $x''$. Since $\alpha$ is a geodesic
    , and $x''$ was the closest point on $\gamma$ to $q$, $\alpha \cdot \gamma_1$ and $\alpha \cdot \gamma_2$ are simple paths, which are the geodesics from $q$ to $x$ and $x'$ respectively. The following diagram shows the points and paths we constructed
    $$\begin{tikzcd}
        & q \arrow[d,"\alpha"] & \\
        & x''  \arrow[dl, "\gamma_1"'] \arrow[dr, "\gamma_2"] & \\
        x & & x'
    \end{tikzcd}$$
    
    We have
    $$d_T(x,x')=d_T(x,x'')+d_T(x',x'')=d_T(q,x)+d_T(q,x')-2d_T(q,x'').$$
    As any connected set in $T$ containing $x$ and $x'$ is itself also a metric tree (since it is connected and has hyperbolicity $0$), the unique simple path it contains between $x$ and $x'$ has to be $\gamma$. This implies that $m_p(x,x') \geq -d_T(p,q)$. Hence,
    $$t_p([x],[x']) \geq -2d_T(p,q)+d_T(p,x)+d_T(p,x')=d_T(x,x'). $$
\end{proof}
\begin{proof}[Proof of Theorem \ref{thm:metric_tree_approximation}]
    Let $r>\dgh(X,T)$. Let $R$ be a correspondence between $X$ and $T$ with $\dis(R)<2r$. Let $q$ be a point in $T$ such that $(p,q) \in R$. If we let $f: X \to \R$, $x \mapsto -d_X(p,x)$, $g: T \to \R$, $y \mapsto -d_T(q,y)$, then $\dis_{f,g}(R) < 2r$. Hence $\dgh((X,f),(T,g))<r$. By Theorem \ref{thm:merge_tree_stability} and Lemma \ref{lem:Tq}, we have
    $$\dgh(T_p,T)=\dgh(T_p,T_q) \leq 12\,\dgh((X,f),(T,g)) \leq 12r.$$
    Since $r>\dgh(X,T)$ was arbitrary, we have $\dgh(T_p,T) \leq 12 \dgh(X,T)$. Therefore,
    $$\dgh(X,T_p) \leq \dgh(X,T)+\dgh(T,T_p) \leq 13\,\dgh(X,T).$$
\end{proof}

\section{Metric Distortion of the Reeb Quotient Map}\label{sec:reeb_distortion}

In this section we generalize the main result of Zinov'ev \cite{z06} from the setting of compact oriented two-dimensional Riemannian manifolds\footnote{A similar question, still in the setting of Riemannian manifolds, was studied by Gromov in \cite[Appendix~1]{frm}.} to the setting of compact geodesic spaces $X$ with Hausdorff dimension two and $\betti(X)<\infty$. We follow a similar line of proof, with the difference that some of the differential topology arguments in \cite{z06} are replaced by topological arguments  directly at the level of Reeb graphs. The main results we prove in this section are the following.

\begin{proposition}\label{prop:reeb-diameter-distortion}
    Let $X$ be a compact geodesic space, $p \in X$ and $f:X \to \R$ defined by $x \mapsto d_X(p,x)$. Let $D:=\sup_{[x] \in \Rp} \diam(\pi_f^{-1}([x]))$. Then, $$\dis(\pi_f) \leq (2\betti(X)+1) D $$
\end{proposition}

\begin{theorem}\label{thm:area_bound}
    Let $X$ be a compact geodesic space with $\betti(X)<\infty$ and Hausdorff dimension two.
    Let $p$ be a point in $X$, and $f:X \to \R$, $x \mapsto d_X(p,x)$. Then,
    $$\dis(\pi_f) \leq 8 (\betti(X)+1)^{3/2} \sqrt{\mathcal{H}^2(X)}. $$
    Here $\mathcal{H}^2(X)$ denotes the two-dimensional Hausdorff measure of $X$.
\end{theorem}

Proofs are deferred to the end of the section. We need some preliminary definitions and results.

\begin{definition}[p-Merging point]
    Let $X$ be a compact geodesic space and $p \in X$. A point $x$ in $X$ is called a \emph{p-merging point} if there are two simple geodesics $\alpha,\beta:[0,1] \to X$ from $p$ to $x$ such that the segment
    $\alpha((t,1)) \cap \beta((t,1)) = \emptyset$ for 
        some $t \in (0,1)$.
\end{definition}

\begin{lemma}\label{lem:merging_points_in_core}
    Let $X$ be a compact geodesic space, $A$ be a closed core of $X$, and $p \in A$. Then, $p$-merging points of $X$ are equal to $p$-merging points of $A$.
\end{lemma}

    \begin{proof}
    Since $A$ is a retract of $X$ by Lemma \ref{lem:retraction}, $A$ is a geodesic subspace of $X$. Therefore, a $p$-merging point in $A$ is also a $p$-merging point in $X$.
    
    Let a point $x$ be in $X\backslash A$. Consider a geodesic \(\gamma: [0,1] \to X\) such that \(\gamma(0) = p\) and \(\gamma(1) = x\). Let \(x' := \gamma(t_{x'})\), for \(t_{x'} \in [0,1]\), be the last point on \(\gamma\) that lies in $A$. 
        Then, the path $\eta: [0, 1- t_{x'}] \to X$ defined by $\eta(t) := \gamma(1-t)$ is a simple path from $x$ to $A$, which intersects $A$ only at its right endpoint $\eta(1-t_{x'})$. The path $\eta$ therefore is a reparametrization of the unique path mentioned in Lemma \ref{lem:unique_path}. In particular, this implies that for any two geodesics from $p$ to $x$ there exists reparameterizations  of them such that they have domain  $[0,1]$ and there exists a point $t_0 \in (0,1)$ such that the two geodesics coincide on $[t_0,1]$.
    Hence, $x$ cannot be a $p$-merging point in this case.

    Finally, let $x$ be a $p$-merging point in $X$. From the above discussion, we deduce that $x \in A$. By Proposition \ref{prop:totally_geodesic}, any geodesic from $p$ to $x$ lies entirely within $A$. Therefore, $x$ is also a $p$-merging point of $A$.

\end{proof}

Before proving the next lemma, we need the following definition:

\begin{definition}[Edge path]\label{def:edge_path}
        Let \( G \) be a finite metric graph. An \emph{edge path} in \( G \) is a sequence of (bijective) parametrizations of directed edges of \( G \) such that for any two consecutive edges \( e_1 \) and \( e_2 \) in the sequence, the terminal vertex of \( e_1 \) coincides with the initial vertex of \( e_2 \). This ordered combination, respecting their given orientations, is referred to as "concatenation".
\end{definition}

\begin{lemma}\label{lem:graph_merging_points}
    Let $G$ be a finite metric graph and $p \in G$. Let $(V,E)$ be the
    $1$-dimensional CW-complex structure on $G$ given in Lemma
    \ref{lem:distance_to_point}. Make $G$ a directed graph by giving
    {an} orientation to each edge along the direction {in which}
    $d_G(p,\cdot)$ increases. Then, the $p$-merging points of $G$ are
    exactly the vertices with indegree at least $2$. Furthermore, the
    number of $p$-merging points of $G$ is less than or equal to $\betti(G)$.
\end{lemma}
\begin{proof}
    Consider a directed edge $[v,w]$. Note that if $x$ is in the interior $(v,w)$, then any geodesic from $p$ to $x$ goes through $v$ and contains $[v,x]$, as otherwise it would include $w$ contradicting the maximality of $f$ at $w$ along $[v,w]$. This shows that interior points of edges are not $p$-merging points.

    Let $w$ be a vertex which is a merging point. Let $\alpha, \beta$ be geodesics from $p$ to $w$ as in the definition of merging points. Note that $\alpha,\beta$ are edge paths, hence the last edges they contain $[v_1,w]$ and $[v_2,w]$ are distinct. This edges are directed, so indegree of $w$ is at least $2$.

    Now let $w$ be a vertex with indegree at least two. Let $v_1,v_2$ be two distinct vertices such that $[v_1,w], [v_2,w]$ are directed edges in $E$. Let $\alpha_i$ be a geodesic from $p$ to $v_i$ for $i=1,2$.
    {
        Then, the concatenation $\alpha_i \ast [v_i, w]$ of $\alpha_i$ and $[v_i,w]$ is a geodesic from $p$ to $w$.
    }
    , as its length is $d_G(p,v_i)+d_G(v_i,w)=d_G(p,v_i)+d_G(p,w)-d_G(p,v_i)=d_G(p,w)$. This shows vertices with indegree two are merging points.

    Let $M$ be the number of $p$-merging points. The number of edges in $G$ is equal to the sum of indegrees of vertices. Let us denote the indegree by $indeg(v)$. Note that only vertex with indegree $0$ is $p$. By Euler's formula we have,
    \begin{equation*}
        \begin{split}
            \betti(G)-1&=\sum_v indeg(v) - \sum_v 1 \\
                       &= \left(\sum_{v: indeg(v) \geq 2 } (indeg(v)-1) \right)-1 \geq M-1
        \end{split}
    \end{equation*}
    Hence $M \leq \betti(G)$.
\end{proof}

\begin{lemma}\label{lem:path_partition}
    Let $G$ be a finite metric graph and $p,x,x' \in G$. Let $\gamma:[0,1] \to G$ be a simple path from $x$ to $x'$. Then there are points $x_1,\dots,x_n$ which $\gamma$ visits in order such that $n \leq 2\betti(G)+1$, and the following holds:
    \begin{enumerate}
        \item and $d_G(p,\cdot)$ is monotone along the part of $\gamma$ between $x_i,x_{i+1}$ for all $i=1,\dots,n-1$.
        \item any geodesic from $x$ to $p$ contains $x_1$,
        \item any geodesic from $x'$ to $p$ contains $x_n$,
    \end{enumerate}
\end{lemma}
\begin{proof}
     Let $V$ be a vertex set of a $1$-dimensional CW-complex structure on $G$ as in Lemma \ref{lem:distance_to_point} with repsect to $p$, containing $x,x'$.  Let $v_1,\dots,v_m$ be the $p$-merging vertices $\gamma$ visits in order. Let $v_0=x,v_{m+1}=x'$. Note that $m \leq \betti(G)$ by Lemma \ref{lem:graph_merging_points}, and part $\gamma_i$ of $\gamma$ in between $v_i,v_{i+1}$ does not contain any $p$-merging vertex for $i=0,\dots,m$. Note that $\gamma_i$ is an edge path, and $d_G(p,\cdot)$ should be either monotone, or initially decreasing then increasing along $\gamma_i$, since otherwise there would be a $p$-merging vertex inside it. Let $w_i$ denote the vertex in $\gamma_i$ where $d_G(p,\cdot)$ achieves its minimum. The following diagram shows how the vertices sits in $\gamma$:
     $$\begin{tikzcd}[column sep=small]
         x=v_0 \arrow[dr, no head] & & v_1 \arrow[dl, no head] \arrow[dr, no head]  & & \cdots \arrow[dr, no head] \arrow[dl, no head]& & v_m  \arrow[dl, no head] \arrow[dr, no head]& & v_{m+1}=x' \arrow[dl, no head] \\
          & w_0  & & w_1 & & w_{m-1} & & w_m
     \end{tikzcd}$$
     Let $(x_1,\dots,x_n)=(w_0,v_1,w_1,\dots,v_m,w_m)$. Note that $n=2m+1 \leq 2\betti(G)+1$ and $d_G(p,\cdot)$ is monotone along the part of $\gamma$ between $x_i,x_{i+1}$ for all $i=1,\dots,n-1$.
    {This proves the first part of the lemma.}
    If $x$ is a merging vertex, then $v_1=v_0=x$, so $x_1=w_0=x$, any geodesic from $x$ to $p$ contains $x_1=x$ trivially. If $x$ is not a merging vertex. Then $[x_1,x]=[w_0,x]$ is the only incoming edge into $x$ where $d_G(p,\cdot)$ is increasing. This shows that any geodesic from $p$ to $x$ contains $x_1$.
    {This proves the second part of the lemma. The third part is proved similarly.}
\end{proof}

\begin{lemma}\label{lem:component-distance}
    Let $X$ be a compact geodesic space and $p,x,y \in X$. Let $f: X \to \R$, $x \mapsto d_X(p,x)$. Let $C_x,C_y$ denote connected component of $x,y$ in $f^{-1}(f(x))$, $f^{-1}(f(y))$ respectively. Assume there is a path $\alpha:[0,1] \to \Rp$ from $[x]$ to $[y]$ such that $\barf \circ \alpha$ is increasing, and $\alpha((0,1))$ does not contain any $[p]$-merging point. Then,
    $$\distance(C_x,C_y)=d_p([x],[y])=d_X(p,y)-d_X(p,x).$$
\end{lemma}
\begin{proof}
    We have $d_X(p,y)-d_X(p,x) \leq d_p([x],[y]) \leq \length_f(\alpha)=d_X(p,y)-d_X(p,x)$. Note that if $x' \in C_x$ and $y' \in C_y$, then $d_X(x',y') \geq d_X(p,y')-d_X(p,x')=d_X(p,y)-d_X(p,x).$ Hence $\distance(C_x,C_y) \geq d_X(p,y)-d_X(p,x) $ Let us show the reverse inequality.

    Let $[z]=\alpha(s)$ for some $s \in (0,1)$. Let us show that any geodesic $\beta$ from $[p]$ to $[z]$ contains $[x]$. Let $[z']=\beta(t)$ be the first point $\beta$ visits in $\gamma$. Assume $[z'] \neq [x]$. Let $\gamma$ be a geodesic from $[p]$ to $[x]$. Then $\gamma \cdot \alpha|_{[0,s]}$ and $\beta|_{[0,t]}$ are simple geodesics from $[p]$ to $[z']$, as $\barf$ is monotone along them. As $[z']$ is not a merging point, these two paths agree at their ends, but this is a contradiction $\beta|_{[0,t]}$ intersects $\alpha$ at its endpoint. Hence $[z']=[x]$.

    Let $t_n$ be an increasing sequence in $(0,1)$ converging to $1$. Let $\alpha(t_n)=[y_n]$. Without loss of generality $y_n$ is converging to $y'$. Note that $[y]=[y']$, so $y' \in C_y$. Let $\gamma_n$ be a geodesic in $X$ from $y_n$ to $p$. Then, by above discussion, $\pi_f \circ \gamma_n$ contains $[x]$. Let $x_n$ be the point on $\gamma_n$ such that $[x_n]=[x]$, so $x_n \in C_x$. Now, we have
    \begin{equation*}
        \begin{split}
            \distance(C_x,C_y) &\leq d_X(x_n,y') \\
            &\leq d_X(x_n,y_n)+d_X(y_n,y') \\ 
            &= d_X(p,y_n)-d_X(p,x_n)+d_X(y_n,y') \\
            &=d_X(p,y_n)-d_X(p,x)+d_X(y_n,y').
        \end{split}
    \end{equation*}
    By letting $n \to \infty$, we get
    $$ \distance(C_x,C_y) \leq d_X(p,y')-d_X(p,x)=d_X(p,y)-d_X(p,x).$$
\end{proof}

\begin{lemma}\label{lem:disconnect_graph}
    Let $G$ be a finite graph, and $A$ be finite set containing $k \geq \betti(G)+1$ points, where each point of $A$ is in the interior of an edge of $G$. Then $G\backslash A$ is disconnected.
\end{lemma}
\begin{proof}
    Let $E$ be a neighborhood of $A$ consisting of $k$ disjoint open intervals, where each interval contains exactly one point of $A$. By excision, $H_n(G,G\backslash A) \simeq H_n(E,E\backslash A)$. Here, we are using real coefficients for homology. By long exact sequence of $(E,E\backslash A)$, we get the following exact $$0 \to H_1(G,G\backslash A) \to H_0(E\backslash A) \to H_0(E)$$
    Note that $H_0(E\backslash A) \to H_0(E)$ is surjective, so it has $k$-dimensional kernel. This implies that $H_1(G,G\backslash A)$ is $k$-dimensional. Long exact sequence of $(G,G\backslash A)$ gives us the following exact sequence:
    $$ 0 \to H_1(G) \to H_1(G,G\backslash A) \to \bar{H}_0(G\backslash A) \to 0$$
    so we have
    $$\dim(\bar{H}_0(G\backslash A))=k-b_1(G) > 0.$$
    This shows that $G\backslash A$ is not path connected. Since $G\backslash A$ is locally path connected, this implies that it is not connected.
\end{proof}

\begin{lemma}\label{lem:disconnect_core}
    Let $X$ be a compact geodesic space, $p,x,x' \in X$. Assume there is a core $G$ of $X$ which is a finite graph, containing $p,x,x'$. Let $f: X \to \R$, $x \mapsto d_X(p,x)$. Let $s=f(y)<f(x)=t$, and $t_0 \in (s,t)$ be such that $f$ does not achieve value $t_0$ at any vertex of $G$. Then there exist $k \leq \betti(X)+1$ points $a_1,\dots,a_k$ in $f^{-1}(t_0) \cap G$, such that any path between $x$ and $x'$ in $X$ contains at least one $a_i$ for some $i$.
\end{lemma}
\begin{proof}
    Note that by Lemma \ref{lem:distance_to_point}, $f^{-1}(t_0) \cap G$ is finite. Any path in $G$ between $x$ and $x'$ intersects $f^{-1}(t_0) \cap G$. Let $A$ be a subset of $f^{-1}(t_0) \cap G$ such that any path between $x$ and $x'$ in $G$ intersects $A$, and $A$ has the minimal possible number of elements among the ones satisfying the same condition. Let $A=\{a_1,\dots,a_k\}$. Note that each $a_i$ lives in the interior of an edge of $G$. Let $E$ be the closed neighborhood of $A$ consisting of $k$ disjoint closed intervals, where each interval contains exactly one point of $A$. We can shrink $E$ further if necessary to exclude a given point outside of $A$.

    By minimality, for each $i=1,\dots,k$, there exists a path $\gamma_i$ from $x$ to $x'$ which does not intersect $A-\{a_i\}$. Note that $\gamma_i$ contains $a_i$, as otherwise there would be a path from $x$ to $x'$ not intersecting $A$.
    
    Let us show that $G\backslash A$ has exactly two connected components. Since $G\backslash A$ is locally path connected, it is equivalent to having exactly two path components. It has at least two path components, corresponding to $x$ and $x'$. Let $x''$ be a point in $G\backslash A$. Without loss of generality, we can assume that the closed neighborhood $E$ of $A$ does not contain $x,x',x''$, and $a_1$ is the closest point to $x''$ in $A$. Let $x_1,y_1$ be the first and last points $\gamma_1$ intersects in the closed interval $I$ containing $a_1$ in $E$. Note that $x_1,y_1$ are endpoints of $I$. Let $\alpha$ be the geodesic in $G$ from $x''$ to $a_1$, and $z_1$ be the first point it intersects in $I$. Then $z_1$ is an endpoint, and either $z_1=x_1$ or $z_1=y_1$. Assume $z_1=x_1$. Then we get a path in $G\backslash A$ between $x$ and $x''$, by concatenating the part of $\gamma_1$ from $x$ to $x_1$ and part of $\alpha$ from $z_1$ to $x''$. Hence $x''$ is in the same component of $G\backslash A$ with $x$. This shows that $G\backslash A$ has exactly two components, corresponding to $x$ and $x'$.

    If we let $A'=\{a_2,\dots,a_k \}$, then $G\backslash A'$ is path connected, since $x,x',a_1$ fall into the same path component through $\gamma_1$, and any other point has a path to either $x$ or $x'$ in the smaller set $G\backslash A$. By Lemma \ref{lem:disconnect_graph}, $k-1\leq \betti(G)=\betti(X)$, so $k \leq \betti(X)+1$.

    Now let $\gamma$ be a path between $x$ and $x'$ in $X$. By \cite[Problem~6.3.11]{engelking1989general}, there is a simple path $\gamma'$ from $x$ to $x'$ contained in $\gamma$. By Proposition~\ref{prop:totally_geodesic}, $\gamma'$ is contained in $G$, so it intersects $A$. Hence $\gamma$ intersects $A$. This completes the proof.
\end{proof}

\begin{lemma}\label{lem:closed_cover}
    Let $A$ be a connected metric space, and $A_1,\dots,A_k$ be a cover of $A$ consisting of closed sets. Then $$\diam(A) \leq \sum_{i=1}^k \diam(A_i).$$
\end{lemma}
\begin{proof}
    Without loss of generality, we can assume that $A_i$ is non-empty for all $i$. Let us prove the claim by induction on $k$. If $k=1$ it is trivial.

    Let $k>1$. $A_1$ intersects $A_i$ for some $i>1$, since otherwise $A_1$, $\cup_{i>1}^k A_i$ would be two non-empty disjoints open sets covering $A$. Assume $A_1$ intersects $A_2$. By inductive assumption,
    $$\diam(A) \leq \diam(A_1 \cup A_2) + \sum_{i=3}^k \diam(A_i).$$
    It remains to show that $\diam(A_1 \cup A_2) \leq \diam(A_1)+\diam(A_2).$ Let $a \in Y_1 \cap Y_2$. Let $x,x' \in A_1 \cup A_2$. If $x,x'$ are both contained in $A_i$ for some $i=1,2$, then $d_A(x,x') 
    \leq \diam(A_i)$ for some $i=1,2$. Assume $x \in A_1$, $x' \in A_2$. Then 
    $$d_A(x,x') \leq d_A(x,a)+d_A(a,x') \leq \diam(A_1)+\diam(A_2).$$
    Therefore, $\diam(A_1 \cup A_2) \leq \diam(A_1)+\diam(A_2)$.
\end{proof}

\begin{lemma}\label{lem:component_diameter}
    Let $X$ be a compact geodesic space with $\betti(X)<\infty$, $p \in X$, and $f:X \to \R$, $x \mapsto d_X(p,x)$. Let $t>0$ and $A$ be a connected component of $f^{-1}(t)$. Then for all but finitely many $s \in [0,t]$, there exists $k \leq \betti(X)+1$ distinct connected components $X_1,\dots,X_k$ of $f^{-1}(s)$ such that
    $$\diam(A) \leq \sum_{i=1}^k \diam(X_i) + 2k(t-s).$$
\end{lemma}
\begin{proof}
    Let $x \in A$. By Theorem \ref{thm:structure}, $\Rf$ has a core $G$ which is a finite metric graph with $\betti(G) \leq \betti(X)$. By Lemma \ref{lem:graph_core_extension}, we can assume that $G$ contains $[p],[x]$. Let $s \in [0,t]$ be a value      that is not attained by $\barf$ on the vertices of $G$.
    
    By Lemma \ref{lem:disconnect_core}, there exists distinct $[x_1],\dots,[x_k]$ in $\barf^{-1}(s)$ such that $k \leq \betti(X)+1$, and every path from $[x]$ to $[p]$ in $\Rf$ contains $[x_i]$ for some $i$. Let $X_i:=\pi_f^{-1}([x_i])$. Note that $X_i$'s are distinct components of $f^{-1}(s)$.

    Let $a \in A$, and $\gamma$ be a geodesic from $a$ to $p$. Note that $\pi_f(a)=[x]$. There exists $i$ such that $\pi_f(\gamma)$ contains $[x_i]$. Hence, there is a point $x_i'$ in $X_i$ such that $d_X(a,x_i')=d_X(a,p)-d_X(x_i',p)=t-s$.

    For $i=1,\dots,k$, define
    $$A_i=\{a \in A: \distance(a,X_i) \leq t-s \} .$$
    $A_i$ is a closed subset for all $i$, and by the discussion above $A_1,\dots,A_n$ cover $A$. By Lemma \ref{lem:closed_cover}, we have
    $$\diam(A) \leq \sum_{i=1}^k \diam(A_i). $$
    Hence, it is enough to show that $\diam(A_i) \leq \diam(X_i)+2(t-s).$ Given $a,a' \in A_i$, there exists $x_i',x_i'' \in X_i$ such that $d_X(a,x_i')\leq t-s$ and $d_X(a',x_i'') \leq t-s$. This implies that
    $$d_X(a,a') \leq d_X(a,x_i')+d_X(x_i',x_i'')+d_X(x_i'',a') \leq \diam(X_i)+2(t-s),$$
    and completes the proof.
\end{proof}

Before proving the next lemma, we will introduce the notion of \emph{distance between two subsets} of a metric space.
For any two non-empty subsets $A,B$ of a metric space $X$, we define the distance between $A$ and $B$, $\distance(A, B)$ as $\inf_{a \in A, b \in B} d_X(a,b)$.

\begin{corollary}\label{cor:component_length}
    Let $X$ be a compact geodesic space with $\betti(X)<\infty$, $p \in X$, and $f:X \to \R$, $x \mapsto d_X(p,x)$. Let $t>0$ and $A$ be a connected component of $f^{-1}(t)$. Then, for all but finitely many $s \in [0,t]$,
    $$\diam(A) \leq \mathcal{H}^1(f^{-1}(s))+(2\betti(X)+2)(t-s). $$
    Here $\mathcal{H}^1$ denotes the one-dimensional Hausdorff measure.
\end{corollary}
\begin{proof}
    Let $X_1,\dots,X_k$ be as in Lemma \ref{lem:component_diameter}. By \cite[Lemma~2.6.1]{bbi01}, $\diam(X_i) \leq \mathcal{H}^1(X_i)$. Hence $\sum_i \diam(X_i) \leq \mathcal{H}^1(f^{-1}(s))$, which implies the desired inequality.
\end{proof}

\begin{proof}[Proof of Proposition \ref{prop:reeb-diameter-distortion}]
    Without loss of generality, we can assume that $\betti(X)<\infty$. Let $x,x' \in X$. By Theorem \ref{thm:structure} and Lemma \ref{lem:graph_core_extension}, there is a finite metric graph core $G$ of $\Rf$ containing $[p],[x],[x']$ such that $\betti(G) \leq \betti(X)$. Note that $d_p([p],x)=d_X(p,x)$.
    
    Let $\gamma$ be a geodesic in between $[x],[x']$ in $\Rf$. By Proposition~\ref{prop:totally_geodesic}, $\gamma \in G$. Let $[x_1],\dots,[x_n]$ be points in $G$ described in Lemma \ref{lem:path_partition} with respect to $[p]$. Let $x_0=x,x_{n+1}=x'$. Let $C_i:=\pi_f^{-1}(\pi_f(x_i))$. By Lemma \ref{lem:component-distance}, we get
    $$d_p([x],[x'])=\sum_{i=0}^n \distance(C_i,C_{i+1}) .$$
    Note that all geodesics from $x$ to $p$ (resp. $x'$ to $p$) intersect $C_1$ (resp. $C_n$), as their images are geodesics in $\Rf$. Note that the part of the geodesic from $x$ to its intersection with $C_1$ realizes $\distance(C_0,C_1)$. Same holds for $C_n,C_{n+1}$. Let $a_i,b_i$ be points in $C_i,C_{i+1}$ respectively such that $\distance(C_i,C_{i+1})=d_X(a_i,b_i)$ for $i=0,\dots,n$. We can let $a_0=x$ and $b_n=x'$.
    We have
    \begin{equation*}
        \begin{split}
            d_X(x,x') &\leq \sum_{i=0}^{n-1} d_X(a_i,b_i)+d_X(b_i,a_{i+1})+d_X(a_n,b_n) \\
                     &\leq nD+\sum_{i=0}^n\distance(C_i,C_{i+1}) \leq d_p([x],[x'])+(2\betti(X)+1)D.
        \end{split}
    \end{equation*}
    Therefore $|d_X(x,x')-d_p([x],[x'])|=d_X(x,x')-d_p([x],[x']) \leq (2\betti(X)+1)D$. This completes the proof.
\end{proof}

\begin{proof}[Proof of Theorem \ref{thm:area_bound}]
    Let $t>0$ and $A$ be a connected component of $f^{-1}(t)$. By Proposition \ref{prop:reeb-diameter-distortion}, it is enough to show that $\diam(A) \leq 4\sqrt{(\betti(X)+1)\mathcal{H}^2(X)}$.

    Let $\delta=\min(t,\diam(A)/(4\betti(X)+4)).$ By Corollary \ref{cor:component_length}, for all but finitely many $s$ in $[t-\delta,t]$,  we have
    $$\diam(A) \leq \mathcal{H}^1(f^{-1}(s))+(2\betti(X)+2)\frac{\diam(A)}{4\betti(X)+4},$$
    implying
    $$\mathcal{H}^1(f^{-1}(s)) \geq \diam(A)/2. $$
    By Eilenberg's inequality \cite[Theorem~13.3.1]{burago2013geometric}, we then obtain
    \begin{equation*}
        \begin{split}
            \mathcal{H}^2(X) \geq \frac{\pi}{4}\int^*_{[t-\delta,t]} \mathcal{H}^1(f^{-1}(s))ds \geq \frac{\diam(A)^2}{16\betti(X)+16}
        \end{split}
    \end{equation*}
where $\int^*_{[t-\delta,t]}$ denotes the upper Lebesgue integral (see~\cite{kennedy1931upper}).
    Hence,
    $$\diam(A) \leq 4\sqrt{(\betti(X)+1)\mathcal{H}^2(X)}. $$
\end{proof}

\appendix
\section{Appendix}

\subsection{Preliminaries on metric geometry}\label{app:metric_geometry}

In this section, we recall some basic definitions and results from metric geometry. For more details, we refer the reader to \cite{bbi01}.

\subsubsection{Length spaces}
Let $(X, d_X)$ be a metric space. We use the word \emph{path} to denote continuous maps from a (closed) interval $I\subset \R$ to $X$.
We use $\mathcal{A}$ to denote the set of all continuous paths in $X$.
We now introduce the notion of length structure on $X$.

\begin{definition}[Length structure and length metric, {\cite[Section~2.1]{bbi01}}]\label{def:length_structure}
    A length structure on a Haudorff space \( X \) consists  a length map \( L: \mathcal{A} \rightarrow \mathbb{R}_{+} \cup \{\infty\} \).
    The length map \( L \) is required to satisfy the following conditions:
    \begin{enumerate}
        \item \textbf{Additivity of Path Length}: For any path \( \gamma:[a, b] \rightarrow X \) and any \( c \in[a, b] \),
              \[ L\left(\gamma_{\left.\right|_{[a, b]}}\right)=L\left(\gamma_{\left.\right|_{[a, c]}}\right)+L\left(\gamma_{\left.\left.\right|_{[c, b]}\right]}\right) \]
        \item \textbf{Continuity of Length}: The length of a piece of a path continuously depends on the piece. That is
        $L(\gamma|_{[a,\cdot]}):[a,b] \to [0,\infty]$ is continuous if $L_f(\gamma)<\infty$,
        \item \textbf{Invariance under Reparameterizations}: The length is invariant under reparameterizations,
              \[ L(\gamma \circ \varphi)=L(\gamma) \]
              for any homeomorphism \( \varphi : [a, b] \to [a, b]\).
        \item \textbf{Compatiability with Topology}: $L_f$ is compatible with the topology of $X$ in the sense that for a neighborhood $U$ of a point $x \in X$, the length of paths conneccing $x$ with points of the complement of $U$ is seperated from $0$:
        \[
            \inf \left\{L(\gamma): \gamma(a)=x, \gamma(b) \in X \backslash U\right\}>0
        \]
    \end{enumerate}
The length structure induces a length metric \( d_L \) on \( X \) defined by
    \begin{equation*}
        d_L(x, y)=\inf \{L(\gamma) ; \gamma:[a, b] \rightarrow X, \gamma \in \mathcal{A}, \gamma(a)=x, \gamma(b)=y\}\
    \end{equation*}
\end{definition}

\begin{definition}[Complete length structure, {\cite[Definition~2.1.10]{bbi01}}]\label{def:complete_length}
A length structure $L$ on a Hausdorff space $X$ is called complete if for every two points $x,x'\in X$ there exists a path $\gamma$ connecting $x$ and $x'$ such that $L(\gamma)=d_L(x,x')$.
\end{definition}

\begin{proposition}[{\cite[Exercise~2.1.2, 2.1.5]{bbi01}}]\label{prop:length_metric}
    Let $X$ be a Hausdorff space and $L$ be a length structure on $X$. Then the length metric $d_L$ is a metric on $X$ that is finer than the original topology of $X$.
\end{proposition}

For any metric space $(X,d_X)$, there is a natural length structure on $X$ induced by the metric $d_X$.

\begin{definition}[Induced length structure, {\cite[Definition~2.3.1]{bbi01}}]\label{def:induced_length}
Let $(X, d_X)$ be a metric space and $\gamma$ be a path in $X$. Then one defines the length of $\gamma$ as
\begin{equation*}
    L_{d_X}(\gamma)=\sup \sum_{i=1}^{n} d_X\left(\gamma\left(t_{i-1}\right), \gamma\left(t_{i}\right)\right)
\end{equation*}
where the supremum is taken over all partitions
$a=t_0<t_1<\cdots<t_n=b$ of $[a,b]$.  It is easy to see that $L_{d_X}$
is a length structure on $X$.
A curve is said to be \emph{rectifiable}.
\end{definition}

\begin{definition}[Geodesic]\label{def:geodesic}
Let $(X,d_X)$ be a metric space. A path $\gamma: [a, b]\subset \R$ in $X$ is called a geodesic
if its length realizes the distance between $a$ and $b$, i.e. if $L_{d_X}(\gamma)=d_X(\gamma(a),\gamma(b))$.
\end{definition}

Note that, starting from a metric space $(X,d_X)$, we can define a
length metric $d_{L_{d_X}}$ on $X$ which in turn give rises to a new
metric space $(X,d_{L_{d_X}})$.

The terminology ``length space" is defined when the metric $d_{L_{d_X}}$ coincides with the original metric $d_X$.

\begin{definition}[Length space, {\cite[Proposition~2.4.1]{bbi01}}]\label{def:length_space}
    A metric space $(X,d_X)$ is called a \emph{length space} if the length metric $d_{L_{d_X}}$ coincides with the original metric $d_X$.
    In this case, $d_X$ is called an intrinsic metric on $X$.
    If $L_{d_X}$ is complete, then $(X,d_X)$ is called a \emph{complete length space}.
\end{definition}

\begin{definition}[Geodesic space]\label{def:geodesic_space}
    A metric space $(X,d_X)$ is called a \emph{geodesic space} if for any $x,x'\in X$, there exists a geodesic $\gamma$ connecting $x$ and $x'$.
\end{definition}

In particular, a geodesic space is a length space.
It turns out whether a metric space is a geodesic space can be
characterized by the following \emph{midpoint} property.

\begin{theorem}[{\cite[Theorem~2.4.16]{bbi01}}]\label{thm:midpoint}
  Let $(X, d)$ be a complete metric space. Then $(X, d)$ is a length space if and only if for any $x, x' \in X$, there exists $x'' \in X$ such that $d(x, x'')=d(x'', x')= \frac{1}{2} d(x, x')$.
\end{theorem}

\subsubsection{Quotient Metric Spaces}\label{app:quotient}

\textbf{Notation:} Given a subspace $A$ of a metric space $(X,d_X)$, and $x \in X$, $D_X(x,A)$ denotes the distance of $x$ to $A$. $(X/A,d_A)$ denotes quotient metric space (see \cite[Definition~3.1.12]{bbi01}) of $X$ under the equivalence relation $x \sim x'$ if $x=x'$ or $x,x' \in A$. We denote the $1$-Lipschitz metric quotient map by $\pi_A:X \to X/A$. If $X$ is geodesic, then so is $X/A$ (see the argument in \cite[p.~62,63]{bbi01}).

 The following result gives an explicit expression of the quotient metric.

\begin{lemma}\label{lem:quotient_metric}
    Let $(X,d_X)$ be a compact metric space and $A$ be a closed subspace of $X$. Then, for all $x,x' \in X$, we have
    $$d_A(\pi_A(x),\pi_A(x'))=\min(d_X(x,x'),D_X(x,A)+D_X(x',A)). $$
    Furthermore, $X/A$ is homeomorphic to the topological quotient of $X$ by $A$.    
\end{lemma}
\begin{proof}
    Let $\sim$ denote the equivalence relation on $X$ given by $x \sim x'$ if $x=x'$ or $x,x' \in A$. Let us denote $\pi_A(x)$ by $[x]$. Let $x,x' \in X$. Let $p_1=x$, $q_1$ (resp. $p_2$) be the closest point to $x$ (resp. $x'$) in $A$ whose existence is guaranteed by the compactness of $A$, and $q_2=x'$. Note that $q_1 \sim p_2$, therefore $$d_A([x],[x']) \leq d_X(p_1,q_1)+d_X(p_2,q_2)=D_X(x,A)+D_X(x',A). $$ We already know that $d_A([x],[x']) \leq d_X(x,x')$. Hence
    $$d_A([x],[x']) \leq \min(d_X(x,x'),D_X(x,A)+D_X(x',A)). $$
    Let us show the reverse inequality. Assume $d_A([x],[x'])=0$, and $x \neq x'$. Let $p_1,\dots,p_n,q_1,\dots,q_n$ be points in $X$ so that $p_1=x$, $q_n=x'$, and $q_i \sim p_{i+1}$ for $i=1,\dots,n-1$. Let $p_{n+1}=x'$. If non of $q_i$ is contained in $A$, then $q_i=p_{i+1}$ for $i=1,\dots,n$, so we have
    $$\sum_{i=1}^n d_X(p_i,q_i) = \sum_{i=1}^n d_X(p_i,p_{i+1}) \geq d_X(x,x').$$
    Now, assume some $q_i$'s are in $A$. Let $j$ be the first index, and $k$ be the last index where $q_i$ is in $A$. Note that $p_{i+1},p_{j+1} \in A$. If $i<j$ or $i>k$, then $q_i=p_{i+1}$. Note that $p_{k+1} \in A$. We have
    \begin{equation*}
        \begin{split}
            \sum_{i=1}^n d_X(p_i,q_i) &= \sum_{i=1}^{j-1} d_X(p_i,p_{i+1}) + \sum_{i=j}^{k} d_X(p_i,q_i)+ \sum_{i=k+1}^{n}d_X(p_i,p_{i+1}) \\
            & \geq d_X(x,p_j)+d_X(p_j,q_j)+d_X(p_{k+1},x') \\
            & \geq d_X(x,q_j) + d_X(p_{k+1},x') \geq D_X(x,A)+D_X(x',A).
        \end{split}
    \end{equation*}
    Taking infimum over all $p_1,q_1,\dots,p_n,q_n$ as above, we get
    $$d_A([x],[x']) \geq \min(d_X(x,x'),D_X(x,A)+D_X(x',A)). $$
    Hence,
    $$d_A([x],[x']) = \min(d_X(x,x'),D_X(x,A)+D_X(x',A)). $$
    If $d_A([x],[x'])=0$ if and only if either $d_X(x,x')=0$ or $D_X(x,A)+D_X(x',A)=0$, or equivalently $x=x'$ or $x,x' \in A$. Hence $d_A([x],[x'])=0$ if and only if $x \sim x'$. Therefore, the map from the topological quotient of $X$ by $A$ to $(X/A,d_A)$ sending the equivalence class of $x$ to $[x]$ is a continuous bijection, which implies that it is a homeomorphism since its domain is compact and range is Hausdorff.
\end{proof}

\begin{corollary}\label{cor:metric_quotient_distortion}
    Let $X$ be a compact metric space and $A$ be closed subspace of $X$. Then $\dgh(X,X/A) \leq \diam(A)/2$.
\end{corollary}
\begin{proof}
    By Lemma \ref{lem:quotient_metric}, we have 
    \begin{equation*}
        \begin{split}
            |d_X(x,x')-d_A(\pi_A(x),\pi_A(x')| &= d_X(x,x')-d_A(\pi_A(x),\pi_A(x')) \\
                                             &\leq D_X(x,A)+\diam(A)+D_X(x',A)-D_X(x,A)-D_X(x',A)\\
                                             &=\diam(A).
        \end{split}
    \end{equation*}
    Hence $\dgh(X,X/A) \leq \dis(\pi_A)/2 \leq \diam(A)/2$.
\end{proof}

The following result shows how the Gromov-Hausdorff distance between quotient metric spaces are controlled.

\begin{lemma}\label{lem:quotient_gh}
    Let $(X,d_X)$ and $(Y,d_Y)$ be compact metric spaces, and $A$ and $B$ are compact subspaces of $X$ and $Y$ respectively. Let $f: X \to \R$, $x \mapsto D_X(x,A)$, and $g: Y \to \R$, $y \mapsto D_Y(y,B)$. Then,
    \begin{equation*}
        \begin{split}
            \dgh(X/A, Y/B) &\leq 2 \dgh((X,f),(Y,g)), \\
            \dgh(A,B) &\leq 5 \dgh((X,f),(Y,g)).
        \end{split}
    \end{equation*}
\end{lemma}

\begin{proof}
    Let $r>\dgh((X,f),(Y,g))$. Let $R$ be a correspondence between $X$ and $Y$ such that distortion $\dis_{f,g}(R) <2r$. Let $[x]$ (resp. $[y]$) denote the image of $x \in X$ (resp. $y \in Y$) under the metric quotient map $X \to X/A$ (resp. $Y \to Y/A)$. Let $\bar{R}$ denote the correspondence between $X/A,Y/A$ given by
    $$\bar{R}:=\{([x],[y]): (x,y) \in R \} .$$
    Given $(x,y), (x',y')$ in $R$, by Lemma \ref{lem:quotient_metric} we have
    \begin{equation*}
        \begin{split}
            |d_A([x],[x'])-d_B([y],[y'])| &= |\min(d_X(x,x'), f(x)+f(x'))-\min(d_Y(y,y'), g(y)+g(y'))| \\
                                          &\leq \max(|d_X(x,x')-d_Y(y,y')|,|f(x)-g(y)|+|f(x')-g(y')| ) < 4r.
        \end{split}
    \end{equation*}
    Hence $\dis(\bar{R}) < 4r$, which implies that $\dgh(X/A,Y/B) < 2r$.

    Let $S$ be the relation between $A$ and $B$ given by
    $$S:=\{(a,b) \in A \times B: \exists (x,y) \in R \text{ such that } d_X(x,a), \, d_Y(y,b) \leq 2r \}. $$
    Let us show that $S$ is a correspondence between $A$ and $B$. Given $a \in A$, let $y \in Y$ such that $(a,y) \in R$. Since $f(a)=0$, we have $g(y)=g(y)-f(a) \leq 2r$, hence there exists $b \in B$ such that $d_Y(y,b) \leq 2r$. This implies that $(a,b) \in S$. Similarly, for each $b \in B$, there exists $a \in A$ such that $(a,b) \in S$. Hence, $S$ is a correspondence between $A$ and $B$. Let $(a,b), (a',b') \in S$. There exists $(x,y), (x',y') \in R$ such that $d_X(x,a),d_X(x',a'), d_Y(y,b), d_Y(y',b') \leq 2r$. We have
    \begin{equation*}
        \begin{split}
            |d_X(a,a')-d_Y(b,b')| &\leq |d_X(a,a')-d_X(x,x')|+|d_X(x,x')-d_Y(y,y')|+|d_Y(y,y')-d_Y(b,b')| \\
                                  &< d_X(a,x)+d_X(a',x')+2r+d_Y(y,b)+d_Y(y',b') <10r.
        \end{split}
    \end{equation*}
    Hence $\dis(S) < 10r$, which implies that $\dgh(A,B) < 5r$. This completes the proof since $r>\dgh((X,f),(Y,g))$ was arbitrary.
\end{proof}

\begin{corollary}\label{cor:quotient_hausdorff}
    Let $X$ be a compact metric space and $A$ and $B$ be closed subsets of $X$. Then,
    $$\dgh(X/A,X/B) \leq d_\mathrm{H}(A,B). $$
\end{corollary}
\begin{proof}
    Let $f: X \to \R$, $x \mapsto D_X(x,A)$, and $g: X \to \R$, $x \mapsto D_X(x,B)$. If $R$ is the identity correspondence between $X$ and $Y$, then $\dis_{f,g}(R)=\sup |f-g| \leq d_\mathrm{H}(A,B)$, implying $\dgh((X,f),(X,g)) \leq d_\mathrm{H}(A,B)/2$. The result follows from Lemma \ref{lem:quotient_gh}.
\end{proof}

\subsection{Finite Graphs}\label{app:graph}

A finite graph is a topological space which has a $1$-dimensional CW-complex structure. A finite metric graph is a geodesic space which is homeomorphic to a finite graph.

\begin{lemma}\label{lem:increasing}
    Let $f:[0,1] \to \R$ be a continuous function such that $\max f = f(1)$ and $f$ does not have any local minimum inside $(0,1)$. Then $f$ is strictly increasing.
\end{lemma}
\begin{proof}
    Since $f$ has no local minimum in $(0,1)$, it is enough to show that $f$ is increasing. Let $s \leq t$. The minimum of $f|_{[s,1]}$ is a achieved at either $s$ or $1$. Since $\max(f)=f(1)$, then $\min f|_{[s,1]}=f(s)$, implying $f(s) \leq f(t)$.
\end{proof}

\begin{lemma}\label{lem:edge_injective}
    Let $G$ be a finite graph and $f: G \to \R$ be continuous function with finitely many local minima. Then, there is a $1$-dimensional CW-complex structure on $G$ such that $f$ is injective over each edge.
\end{lemma}

\begin{proof}
    Consider a $1$-dimensional CW-complex structure on $G$ whose vertex set contains all the local minima. Extend the vertex set by adding one point from each edge $e$ where $f|_e$ achieves its maximum. With this new CW complex structure, $f$ does not have any local minima in the interior of any edge, and over each edge it achieves its maximum at one of the endpoints. Now the result follows from Lemma \ref{lem:increasing}.
\end{proof}

\begin{lemma}\label{lem:distance_to_point}
    Let $(G,d_G)$ be a finite metric graph. Let $V$ be the vertex set of a $1$-dimensional CW-complex structure on $G$ containing $p$. Let $f: G \to \R$, $x \mapsto d_G(p,x)$. By adding at most one point from each edge to the vertex set if necessary, we can guarantee that $f$ maps each edge isometrically into $\R$. Furthermore, edges in this new CW-complex structure are geodesics.
\end{lemma}
\begin{proof}
     From each edge, chose a point where $f$ takes its maximum, and add it to the set of vertices. With this $1$-dimensional CW complex structure, on each edge $f$ achieves its maximum at one of the endpoints. Let us show that with this structure, $f$ maps each edge isometrically into $\R$. It is enough to show that this happens at the interior of each edge. Let $e$ be an edge with endpoints $x,y$, where the maximum of $f$ over $e$ is obtained at $y$. Let $x_1,x_2$ be two distinct points in the interior of $e$. Without loss of generality, we can assume that $x_1$ lies in between $x$ and $x_2$. The geodesic $\gamma$ from $p$ to $x_2$ contains either $x$ or $y$, but it cannot contain $y$ since otherwise $f(x_2)=d_G(p,x_2)$ would be greater than $f(y)=d_G(p,y)$. Therefore, $\gamma$ passes through $x$ and contains the part of the edge from $x$ to $x_2$, which includes $x_1$. Therefore $d_G(x_1,x_2)=d_G(p,x_2)-d_G(p,x_1)=f(x_2)-f(x_1)$. This means $f$ isometrically embeds the interior of $e$, hence the whole $e$, into $\R$. Note that the length of $e$ is $|f(x)-f(y)|=|d_G(p,x)-d_G(p,y)| \leq d_G(x,y)$, hence $e$ is a geodesic.
\end{proof}

\begin{lemma}\label{lem:edge_length}
    Let $G$ be a finite metric graph with a 1-dimensional CW-complex structure. Then, length of each edge is less than or equal to $2 \diam(G)$.
\end{lemma}
\begin{proof}
    By Lemma \ref{lem:distance_to_point}, each edge can be divided into at most two geodesics, hence its length is bounded by $2\diam(G)$.
\end{proof}

\begin{lemma}\label{lem:edge_difference}
    Let $G$ be a finite topological graph, with $N$ edges. If $(G,d_1)$ and $(G,d_2)$ are metric graphs such that edge lengths with respect to $d_1,d_2$ differs at most by $\epsilon$, then $\dgh((G,d_1),(G,d_2)) \leq N\epsilon/2$.
\end{lemma}
\begin{proof}
    Let $L_1,L_2$ denote the length structures of $(G,d_1),(G,d_2)$ respectively. For each edge $e$, let $\alpha_e, \beta_e: [0,1] \to G$ be constant speed parametrization of $e$ with respect to $d_1,d_2$ respectively, with same endpoints. Let $R$ be the correspondence between $(G,d_1),(G,d_2)$ given by
    $$R:=\{(\alpha_e(t),\beta_e(t): e \text{ is an edge in G, } t \in [0,1] \} .$$ Note that $R$ is a homeomorphism.
    Let $(x,y), (x',y') \in R$. Let $\gamma$ be a geodesic between $x$ to $x'$ in $(G,d_1)$. If $\gamma$ is contained in a single edge $e$, then $\gamma=\alpha|_{[s,t]}$ for some $s,t \in [0,1]$. Then $\beta|_{[s,t]}$ is a path between $y$ and $y'$, and we have
    $$d_2(y,y')-d_1(x,x') \leq L_2(\beta|_{[s,t]})-L_2(\alpha|_{[s,t]})=(s-t)(L_2(e)-L_1(e)) \leq \epsilon.$$
    
    Now, assume $\gamma$ is not contained in a single edge.  Then $\gamma$ can be decomposed as $\gamma=\gamma_1 \cdot \gamma_n$, where $\gamma_i$ is contained in edge $e_i$, $(e_i)'s$ are distinct, and for $i=2,\dots,n-1$, $\gamma_i$ is exactly the edge $e_i$. Changing the orientation of $\alpha_{e_1}$ and $\beta_{e_1}$ if necessary, we can assume that $\gamma_1$ is equal to $\alpha_{e_1}|_{[s,1]}$ for some $s \in [0,1]$. Then $x=\alpha_{e_1}(s)$, and $y=\beta_{e_1}(s)$. Similarly, $\gamma_n$ is equal to $\alpha_{e_n}|_{[0,t]}$ for some $t \in [0,1]$, and $x'=\alpha_{e_n}(t)$, $y'=\beta_{e_n}(t)$.  We have,
    \begin{equation*}
        \begin{split}
            L_2(\beta_{e_1}|_{[s,1]})-L_1(\alpha_{e_1}|_{[s,1]})| & = (1-s)(L_2(e_1)-L_1(e_1)) \leq \epsilon \\
            L_2(\beta_{e_n}|_{[0,t]})-L_1(\alpha_{e_n}|_{[0,t]})| & = t(L_2(e_n)-L_1(e_n)) \leq \epsilon.
        \end{split}
    \end{equation*}
    Note that $\gamma'=\beta_{e_1}|_{[s,1]} \cdot e_1 \cdot \dots e_{n-1} \cdot \beta_{e_n}|_{[0,t]} $ is a path from $y$ to $y'$. Hence, we have
    \begin{equation*}
        \begin{split}
            d_2(y,y')-d_1(x,x') &\leq L_2(\gamma')-L_1(\gamma) \\
                                &\leq L_2(\beta_{e_1}|_{[s,1]})-L_1(\alpha_{e_1}|_{[s,1]})| + \sum_{i=2}^{n-1} L_2(e_i)-L_1(e_i) + L_2(\beta_{e_n}|_{[0,t]})-L_1(\alpha_{e_n}|_{[0,t]})| \\
                                &\leq N \epsilon.
        \end{split}
    \end{equation*}

    Similarly, starting from a geodesic between $y,y'$ in $(G,d_2)$, we can show that
    $$d_1(x,x')-d_2(y,y') \leq  N \epsilon. $$
    Hence $\dis(R) \leq N \epsilon$, which completes the proof.
\end{proof}

\begin{lemma}\label{lem:graph_convergence}
    Let $(G_n)$ be a sequence of finite metric graphs  with bounded diameter and number of edges. Then $(G_n)$ has a convergent subsequence Gromov-Hausdorff converging to $G$, where $G$ is a finite metric graph and $\betti(G) \leq \limsup \betti(G_n)$.
\end{lemma}
\begin{proof}
    Number of vertices of $(G_n)$ is upper bounded, hence by passing to a subsequence we can assume that it is constant and equal to $V$. Similarly, we can assume that number of edges is constant, and equal to $E$. Let us prove by induction on $E$. If $E=0$, and everything is a point and the result is trivial. Let $\epsilon_n$ denote the length of the shortest edge $e_n$ in $G_n$. By passing to a subsequence, we can assume that $\epsilon_n$ is convergent. Assume $\epsilon_n \to 0$. Then, by inductive hypothesis, $G_n/e_n$ has a subsequence Gromov-Hausdorff converging to a finite graph $G$, where $\betti(G) \leq \limsup \betti(G_n/e_n) \leq \limsup \betti(G_n)$. By Corollary \ref{cor:metric_quotient_distortion}, $\dgh(G_n,G_n/e_n) \leq \epsilon_n/2$. This implies that the corresponding subsequence of $G_n$ also converges to $G$. 
    
    Now, we can assume that $\epsilon_n$ converges to $2\epsilon>0$. By passing to a subsequence, we can assume that the shortest edge in $G_n$ has length at least $\epsilon$ for all $n$. Let $V_n$ denote the vertex set of $G_n$, and order the elements of $V_n$ so that $V_n=\{v_n^1,\dots,v_n^V\}$. Let $E_n^{i,j}$ denote the number of edges between $v_n^i,v_n^j$ in $G_n$. By passing to a subsequence, we can assume that $E_n^{i,j}=E^{i,j}$ is a constant with respect to $n$. For each $n$, order the edges between $v_n^i$ and $v_n^j$, and let $l_n^{i,j} \in \R^{E^{i,j}} $ denote the vector whose entries are given by the lengths of edges between $v_n^i,v_n^j$, with the assumed order. Note that each entry of $l_n^{i,j}$ is at least $\epsilon$ and at most $2D$ by Lemma \ref{lem:edge_length}, where $D$ is the diameter limit. Hence, by passing to a subsequence, we can assume that $l_n^{i,j}$ converges to $l_{i,j}$, where entries of $l^{i,j}$ is at least $\epsilon$ and at most $2D$. Now, let $G$ be the metric graph whose vertes set is $\{v^1,\dots,v^V \}$, and between $v^i,v^j$, there are $E^{i,j}$ edges whose lengths are given by $l^{i,j}$. By Lemma \ref{lem:edge_difference}, $G_n$ converges to $G$. As $G$ has exactly the same graph structure with the constructed subsequence $G_n$, $\beta_1(G)=\beta_1(G_n)$ for all $n$. This completes the proof.
\end{proof}

\subsection{Hyperbolicity and Metric Trees}\label{app:hyperbolicity}

\begin{definition}[Gromov product]
    Let $(X,d_X)$ be a metric space and $p \in X$. The \emph{Gromov product} $g_p: X \times X \to \R$ is the function defined by
    $$g_p(x,x'):=(d_X(p,x)+d_X(p,x')-d_X(x,x'))/2. $$
\end{definition}

\begin{definition}[Hyperbolicity]
    Let $(X,d_X)$ be a metric space and $p \in X$. The \emph{p-hyperbolicity} $\hyp_p(X)$ of $X$ is the infimum of all $\delta \geq 0$ such that
    $g_p(x,z'') \geq \min(g_p(x,x'),g_p(x',x''))-\delta$
    for all $x,x',x''$ in $X$. The hyperbolicity $\hyp(X)$ of $X$ is defined as
    $$\hyp(X):= \sup_{p \in X} \hyp_p(X).$$ By \cite[Corollary~1.1.B]{g87}, $\hyp(X) \leq 2\hyp_p(X)$ for all $p \in X$.
\end{definition}

\begin{lemma}\label{lem:hyperbolicity}
    Let $(X,d_X),(Y,d_Y)$ be metric spaces. Then $|\hyp(X)-\hyp(Y)| \leq 6 \dgh(X,Y)$.
\end{lemma}
\begin{proof}
    Let $r>\dgh(X,Y)$, and $R$ be a correspondence between $X$ and $Y$ such that $\dis(R) < 2r$. Let $p,x_1,x_2,x_3 \in X$. Let $q,y_1,y_2,y_3$ in $Y$ such that $(p,q),(x_i,i) \in R$ for all $i$. We have $$|g_p(x_i,x_j)-g_q(y_i,y_j)| \leq 3r$$ for all $i,j$. Assume $\hyp(Y) < \delta$. Then, we have
    \begin{equation*}
        \begin{split}
            g_p(x_1,x_3) &\geq g_q(y_1,y_3)-3r
                         \geq \min(g_q(y_1,y_2),g_q(y_2,y_3))-\delta - 3r \\
                         &\geq \min(g_p(x_1,x_2),g_p(x_2,x_3))-\delta - 6r.
        \end{split} 
    \end{equation*}
    This shows that $\hyp_p(X) \leq \hyp(Y)+6r$. Taking supremum over $p$, we get $\hyp(X) \leq \hyp(Y) + 6r$. Similarly, $\hyp(Y) \leq \hyp(X)+6r$, so $|\hyp(X)-\hyp(Y)| \leq 6r$. This completes the proof, as $r>\dgh(X,Y)$ was arbitrary.
\end{proof}

    Recall that a \emph{simple path} is a path that is injective.

\begin{definition}[Metric Tree]
    A metric tree is a geodesic space such that up to reparametrization, there is a unique simple path between each pair of points. A geodesic space $X$ is a metric tree if and only if $\hyp(X)=0$ \cite[Example~2.11, Lemma~2.13]{bestvina2002r}. Note that $\hyp(X)=0$ if and only if $\hyp_p(X)=0$ for some $p \in X$.
\end{definition}

\begin{lemma}\label{lem:tree_disconnected}
    Let $T$ be a metric tree and $p,x,y \in T$ contained in the interior of a simple path between $x,y$. Then $T-\{p\}$ is disjoint union of two open sets $U$ and $V$ such that $x \in U$, $y \in V$.
\end{lemma}
\begin{proof}
    Define a relation on $T-\{p\}$ by $w~z$ if the path between $w$ and $z$ does not contain $p$. Note that this is an equivalence relation, as the simple path between $w,z$ is contained in the union of paths between $w,w'$ and $w',z$. If $d_T(w,p)>\epsilon$, then the open ball $B_\epsilon(w)$ is contained in the equivalence class of $w$, hence equivalence classes are open. Since $x,y$ are in distinct equivalence classes, if we let $U$ be the equivalence class of $x$ and $V$ be the union of all other equivalence classes, then $U,V$ satisfies the necessary conditions.
\end{proof}

\begin{lemma}\label{lem:tree_connected}
    Any connected metric space $X$ with $0$-hyperbolicity is a metric tree.
\end{lemma}
\begin{proof}
    \cite[Lemma~2.13]{bestvina2002r}, $X$ isometrically embeds into a metric tree $T$. It is enough to show that for distinct $x,x' \in X$, then $X$ contains the path between $x,x'$ in $T$. Let $p \in T$ be a contained in the interior of a simple path between $x,x'$ in $T$. Let $U,V$ be as in Lemma \ref{lem:tree_disconnected}. We have $X-\{p\}$ disconnected, as it is the union of non-empty open sets in $X$ given by $U \cap X$, $V \cap X$. Since $X$ is connected, this implies that $p \in X$. This completes the proof.
\end{proof}

\begin{lemma}\label{lem:limit_tree}
    Gromov-Hausdorff limit of metric trees is a metric tree.
\end{lemma}
\begin{proof}
    Follows from Lemma \ref{lem:hyperbolicity}.
\end{proof}

\begin{lemma}\label{lem:tree_core}
    Let $(T,d_T)$ be a compact metric tree and $A$ be a closed connected subspace of $T$. Then $T/A$ is a metric tree.
\end{lemma}
\begin{proof}
    Let $p$ be the point in $T/A$ corresponding to $A$. Let us show that $\hyp_p(T/A)=0$.

    For $x \in T$, we denote the image of $x$ in $T/A$ under the quotient map by $[x]$. By Lemma \ref{lem:quotient_metric}, $d_A(p,[x])=D(x,A)$.
    Given $x,y \in T$,
    if the geodesic between $x$ and $y$ intersects $A$, then $d_A([x],[y])=D(x,A)+D(y,A)$, hence $g_p([x],[y])=0$. Given $x,y,z$ in $T$, if the geodesic between $x$ and $y$ or $z$ and $y$ intersects $A$, then we have $g_p([x],[z]) \geq 0 = \min(g_p([x],[y]),g_p([y],[z])$.

    Now, we can assume that the geodesic between $x,y$ and $y,z$ does not intersect $A$. Note that there is at most unique point $a$ in $A$ so that the geodesic from $x$ to $a$ does not intersect $A$ except at the endpoints, as if $a'$ was another such point, then the geodesic from $x$ to $a$ concatenated with the simple path in $A$ between $a,a'$ would be a distinct simple path from the geodesic from $x$ to $a'$. If $a$ is the closest point to $x$ in $A$, then the geodesic from $x$ to $a$ does not intersect $A$ except at $a$. Let $a_x,a_y,a_z$ denote a closest point of $A$ to $x,y,z$ respectively. Let $w$ be the closest point on the geodesic between $x$ and $y$ to $A$. If we concatenate the geodesic from $x$ to $w$ with the geodesic from $w$ to the closest point $a_w$ in $A$, we get the geodesic $x$ to $a_w$, which intersect $A$ only at $a_w$. Hence $a_x=a_w$. Similarly $a_y=a_w$, so $a_x=a_y$. Furthermore, $D(x,A)+D(y,A)=d_T(x,a_x)+d_T(x,a_y) \geq d_T(x,w)+d_T(y,w)=d_T(x,y)$, hence by Lemma \ref{lem:quotient_metric}, $d_A([x],[y])=d_T(x,y)$. The same argument shows that $a_y=a_z$, and $d_A([y],[z])=d_T(y,z)$, $d_A([x],[z])=d_T(x,z)$. Hence, $a_x=a_y=a_z=a$. Now, we have
    $$g_p([x],[z])=(D(x,A)+D(z,A)-d_A([x],[z]))/2=(d_T(x,a)+d_T(z,a)-d_T(x,z))/2=g_a(x,z).$$
    Similarly, $g_p([x],[y])=g_a(x,y)$, $g_p([y],[z])=g_a(y,z)$. Since $T$ is a metric tree, we have
    $$g_p([x],[z])=g_a(x,z) \geq \min(g_a(x,y),g_a(y,z))=\min(g_p([x],[y]),g_p([y],[z])).$$
    This shows that $\hyp_p(T/A)=0$ and completes the proof.
\end{proof}

\subsection{Cores of Geodesic Spaces}\label{app:core}

Let us review the notion of \emph{core} which we introduced in Section \ref{sec:structure}
\begin{definition}[Core of a metric space]
    Let $(X,d_X)$ be a compact metric space. We call a closed subspace $A$ of $X$ a \emph{core} of $X$ if $X/A$ is a metric tree. 
\end{definition}

\begin{lemma}\label{lem:core_extension}
    If $A$ is a core of compact metric space $X$, and $B$ is a closed connected space containing $A$, then $B$ is a core of $X$.
\end{lemma}
\begin{proof}
    By Lemma \ref{lem:tree_core}, it is enough to show that $X/B=(X/A)/(B/A)$, as $B/A$ is a connected subspace of the metric tree $X/A$. By $1$-Lipschitzness of $\pi_A: X \to X/A$, we have $D(x,B) \geq D(\pi_A(x),B/A)$ for all $x$ in $X$. Let $b_0$ be the closest point in $B$ to $x$. For any $b \in B$, by Lemma \ref{lem:quotient_metric}, we have
    \begin{equation*}
        \begin{split}
            d_A(\pi_A(x),\pi_A(b))&=\min(d_X(x,b),D(x,A)+D(b,A)) \\ 
                                  &\geq \min(d_X(x,b_0),d_X(x,b_0)+D(b,A))=d_X(x,b_0)=D(x,B).
        \end{split}     
    \end{equation*}
    This shows that $D(\pi_A(x),B/A)=D(x,B).$ Therefore, for any $x, x'$ in $X$, we have
    \begin{equation*}
        \begin{split}
            d_{B/A}(\pi_{B/A}(\pi_A(x)), \pi_{B/A}(\pi_A(x'))) &= \min(d_A(\pi_A(x),\pi_A(x')), D(x,B)+D(x',B)) \\
                                        &= \min(d_X(x,x'),D(x,A)+D(x',A),D(x,B)+D(x',B)) \\
                                        &= \min(d_X(x,x'),D(x,B)+D(x',B))=d_B(\pi_B(x),\pi_B(x')).
        \end{split}
    \end{equation*}
\end{proof}

\begin{lemma}\label{lem:unique_path}
    Let $(X,d_X)$ be a geodesic space and $A$ be a core of $X$. Up to reparametrization, there is a unique path $\gamma:[0,1] \to X$ such that $\gamma(0)=x$, $\gamma(1) \in A$, $\gamma([0,1)) \cap A = \emptyset.$ The endpoint $\gamma(1)$ is the unique closest point to $x$ in $A$, and $\gamma$ is a geodesic. Furthermore, for each $a$ in $A$, $r^{-1}(a)-\{ a\}$ is open in $X$.
\end{lemma}
\begin{proof}
    Let $\gamma_1,\gamma_2$ be simple paths from $x$ to points in $A$ as in the statement. By Lemma \ref{lem:quotient_metric}, $\pi_A \circ \gamma_1,\pi_A \circ \gamma_2$ are simple paths in $X/A$ with same endpoints. Since $X/A$ is a metric tree, there is a reparametrization $\tau: [0,1] \to [0,1]$ such that $\pi_A \circ \gamma_1 = \pi_A \circ \gamma_2 \circ \tau$. This implies that $\gamma_1=\gamma_2 \circ \tau$ over $[0,1)$. By continuity, $\gamma_1=\gamma_2 \circ \tau$ over their whole domain $[0,1]$.

    Let $a_x$ be the closest point to $x$ in $A$, and let $\gamma$ be the geodesic from $a_x$ to $a$. Note that $\gamma$ satisfies the conditions of the statement. Hence, up to reparametrization, $\gamma$ is the unique path satisfying these conditions. This also shows that there is a unique closest point to $x$ in $A$.

    Let $a$ be any point in $A$ and $x$ be a point in $r^{-1}(a)-\{a\}$. For any $x' \notin r^{-1}(a)$, the geodesic $\gamma$ from $x$ to $x'$ intersects $A$, as otherwise $r(x')=r(x)=a$. Furthermore, $a=r(x)$ is contained in $\gamma$, and it is the first point $\gamma$ intersects $A$. Therefore, $d_X(x,x')\geq d_X(x,a)$. If we let $\epsilon=d_X(x,a)$, then we have the open $\epsilon$ ball $B^\epsilon(x)$ of $x$ contained in $r^{-1}(a)-\{a\}$. This completes the proof.

\end{proof}

\begin{lemma}\label{lem:retraction}
    Let $X$ be a geodesic space and $A$ be a core of $X$. Let $r: X \to A$ be the map sending $x$ to the closest point in $A$. Then, $r$ is a $1$-Lipschitz retraction, and fibers of $r$ are metric trees, which are isometrically mapped into $X/A$ under $\pi_A: X \to X/A$. 
\end{lemma}
\begin{proof}
    Let $x,y \in X$, and $\gamma$ be a geodesic from $x$ to $y$.

    Assume $\gamma$ does not intersect $A$. Let $z$ be a closest point $A$ on $\gamma$. Let $\alpha$ be the geodesic from $z$ to the its closest point in $A$. If we let $\gamma_1$ (resp. $\gamma_2$) be the part of $\gamma$ from $x$ (resp. $y$) to $z$. Then $\gamma_1 \circ \alpha$ and $\gamma_2 \circ \alpha$ satisfy the conditions of Lemma \ref{lem:unique_path}. Hence $a:=r(x)=r(y)=r(z)$.

    Assume $\gamma$ intersect $A$. Let $a$ be the first and $b$ be the last point it intersects $A$. As the part of $\gamma$ from $x$ to $a$ and $y$ to $b$ satisfies the conditions of Lemma \ref{lem:unique_path}, $r(x)=a$, $r(y)=b$, so we have
    $$d_X(x,y)=d_X(x,r(x))+d(x,r(y))+d_X(r(x),r(y)). $$
    This completes the proof of $1$-Lipschitzness of $r$.

    Let us show that $\pi_A: X \to X/A$ maps the fibers of $r$ isometrically into $X/A$. Assume $r(x)=r(y)=a$. Note that $D(x,A)=d_X(x,a)$ and $D(y,A)=d_X(y,a)$. Hence, by Lemma \ref{lem:quotient_metric}, we have
    $$d_A(\pi_A(x),\pi_A(y))=\min(d_X(x,y),d_X(x,a)+d_X(y,a))=d_X(x,y). $$
    Since any connected subspace of a metric tree is a metric tree, it remains to show that fibers of $r$ are connected. This is true since for any $x$, the geodesic from $x$ to $r(x)$ is contained in the fiber $r^{-1}(r(x))$.
\end{proof}

{\color{orange}

}

\begin{proposition}\label{prop:deformation_retraction}
    Let $(X,d_X)$ be a geodesic space and $A$ be a core of $X$. Then there is a deformation retraction $H:X \times [0,1] \to A$ such that $H(\cdot,t)$ is $1$-Lipschitz for all $t$ and $H(\cdot,1)=r$, where $r:X \to A$ is the retraction described in Lemma \ref{lem:retraction}, mapping each point to the unique closest point in $A$.
\end{proposition}

{
\begin{proof}
        Let \( r: X \to A \) be the retraction map given in Lemma \ref{lem:retraction}. Let \(\gamma_{x}\) denote the unique constant speed geodesic between \(x\) and \(r(x)\). The uniqueness is given by Lemma \ref{lem:unique_path}. Let \(H: X \times [0,1] \to X\), \((x,t)\) map to \(\gamma_{x}(t)\). It is evident that \(H(\cdot,0) = Id_X\) and \(H(\cdot,1) = r\). What's left to prove is the continuity of \(H\).

        Consider \(\gamma\), a geodesic between \(x\) and \(x'\) in \(X\). If \(\gamma\) intersects \(A\), then by Lemma \ref{lem:unique_path}, the segment of \(\gamma\) starting from \(x\) up to its first intersection point with \(A\) agrees with \(\gamma_{x}\) up to a reparametrization. Analogously, the segment from \(x'\) up to its last intersection point with \(A\) is \(\gamma_{x'}\). Thus,
        \(d_X(H(x,s),H(x',t)) \leq d_X(x,x')\) for all \(s,t \in I\) and the continuity of \(H\) follows.
    
        Now, if \(\gamma\) doesn't intersect \(A\), let \(x''\) be the closest point to \(A\) on \(\gamma\). Let's denote the segment of \(\gamma\) from \(x\) to \(x''\) as \(\gamma_1\) and from \(x'\) to \(x''\) as \(\gamma_2\). Notably, both \(\gamma_1 \cdot \gamma_{x''}\) and \(\gamma_2 \cdot \gamma_{x''}\) meet the requirements of Lemma \ref{lem:unique_path}. Thus, after suitable reparametrization, they coincide with \(\gamma_{x}\) and \(\gamma_{x'}\) respectively.
        Let $s,t \in I$.
        If both \(\gamma_{x}(s)\) and \(\gamma_{x'}(t)\) lie within \(\gamma\), then \(d_X(H(x,s),H(x',t)) \leq d_X(x,x')\) and we're done. If one of them isn't part of \(\gamma\), then it must be part of \(\gamma_{x''}\). Without loss of generality, let's assume \(\gamma_{x}(s)\) is part of \(\gamma_{x''}\). Therefore, both \(\gamma_{x}(s)\) and \(\gamma_{x'}(t)\) are part of \(\gamma_{x'}\). Let \(a = r(x'')\). We have that
        \begin{equation*}
            \begin{split}
                d_X(H(x,s),H(x',t)) &\leq |d_X(\gamma_{x}(s),a)-d_X(\gamma_{x'}(t),a)| \\
                                   &=|(1-t)\,d_X(x,r(x''))-(1-s)\,d_X(x',r(x''))| \\
                                   &\leq |s-t|\,d_X(x,r(x''))+|1-t|\,d_X(x,x') \leq |s-t|\,\diam(X)+d_X(x,x').
            \end{split}
        \end{equation*}
    
        In conclusion, for any \(x,x'\) in \(X\), it holds:
        $$d_X(H(x,s),H(x',t)) \leq |s-t|\,\diam(X)+d_X(x,x'). $$
        This confirms that \(H:X \times I \to X\) is continuous and \(H(\cdot,t)\) is $1$-Lipschitz.
    \end{proof}

}

\begin{lemma}\label{lem:local_minima}
    Let $(X,d_X)$ be a geodesic space and $A$ be a core of $X$. Let $f: X \to \R$ be a continuous function with finitely many local minima. Then number of local minima of $f|_A: A \to \R$ is less than or equal to the number of local minima of $f$.
\end{lemma}
\begin{proof}
    Let $r:X \to A$ be the retraction given Lemma \ref{lem:retraction}. Given $a$ in $A$, let $T_a:=r^{-1}(a)$, and $x_a$ be a point in in $T_a$ where $f|_{T_a}$ achieves its global minimum. Note that if $x_a \neq a$, then by Lemma \ref{lem:retraction}, $x_a$ is a local minimum of $f$, so $x_a=a$ for all but finitely many $a$. Now assume $a$ is local minimum of $f|_A$. Let us show that $x_a$ is a local minimum of $a$. If $x_a \neq a$, we already know that it is a local minimum. Assume $x_a=a$. Let $V$ be a neighborhood of $a$ in $A$ such that for all $a' \in A$, $f(a')\geq f(a)$ and $x_{a'}=a'$. Let $U = r^{-1}(V)$. Note that $U$ is a neighborhood of $a$ in $X$. Let $x \in U$. Let us show that $f(x) \geq f(a)$. Let $a'=r(x)$. Since $x \in T_{a'}$, we have $f(x) \geq f(x_{a'})=f(a') \geq f(a)$.

    Let $C$ be the set of local minima of $f|_A$ and $D$ be the set of local minima of $f$. The argument above shows that we have a function $\iota: C \to D$, $a \mapsto x_a$. If $a \neq a'$, then $x_a \neq x_{a'}$ since $x_a,x_{a'}$ are contained in distinct fibers of $r^{-1}(a), r^{-1}(a')$ of $r$ respectively. Hence, $\iota$ is injective.
\end{proof}

The following proposition will be used in Section~\ref{sec:reeb_distortion}.
\begin{proposition}\label{prop:totally_geodesic}
    Let $X$ be a compact geodesic space and $A$ be a closed core of $X$. Let $x,y \in A$ and $\gamma$ be a simple path between $x,y$ in $X$. Then $\gamma$ is contained in $A$. In particular, $A$ is a totally geodesic subspace of $X$.
\end{proposition}
\begin{proof}
    Let us prove by contradiction. Assume there exists $t \in (0,1)$ such that $\gamma(t) \notin A$. Let $t_0:=\sup\{s < t: \gamma(s) \in A \}$, and $t_1:=\inf\{s > t_0: \gamma(s) \in A \}$. Note that $0 \leq t_0 < t < t_1 \leq 1$. The parts of $\gamma$ between $[t_0,t]$ and $[t,t_1]$ provides distinct simple paths $\alpha,\beta$ from $\gamma(t)$ to $A$ such that $\alpha((0,1))\cap \beta((0,1))=\emptyset$. This is a contradiction by Lemma \ref{lem:unique_path}, since $\alpha((0,1)),\beta((0,1) \subseteq X-A$.
\end{proof}

Now, let us analyze cores of finite metric graphs.

\begin{lemma}\label{lem:graph_core}
    Let $G$ be a finite metric graph and $H$ be a subgraph of $G$ with $\betti(H)=\betti(G)$. Then $H$ is a core of $G$.
\end{lemma}
\begin{proof}
    Note that $H_1(H) \to H_1(G)$ is a bijection. By the long exact homology sequence of $(G,H)$ we get $H_1(G/H)=H_1(G,H)=0$. Since $G/H$ is also a finite metric graph, $H_1(G/H)=0$ implies that it is a metric tree.
\end{proof}

\noindent\textbf{Notation:} Give a finite graph $G$, let $V(G), E(G)$ denote the minimum number of vertices and edges a $1$-dimensional CW complex structure on $G$ can have respectively.
The following result shows that finite metric graphs has \emph{simple} cores.

\begin{lemma}\label{lem:simple_core}
    Let $G$ be a finite metric graph. Then $G$ has a core $G_0$, which is itself a metric graph and $E(G_0) \leq 3\betti(G)$.
\end{lemma}
\begin{proof}
    Let $G_0$ be a core of $G$ which is a finite metric graph, has minimum $V(G_0)$ among such cores. Let us show that $E(G_0) \leq C$. Let $V$ be a vertex set of $G_0$ achieving $|V|=V(G_0)$. Assume $E(G_0) \geq 1$, as otherwise the result is already true. Let $p$ be the vertex in $V$ with minimal degree. Note that $\deg(p)>1$, as otherwise we could remove the vertex $p$ and the interior of the the unique edge coming to it from $G_0$, and the resulting subgraph $G_1$ would still be a core with less number of vertices by Lemma \ref{lem:graph_core}. Assume $\deg(p)=2$. If two distinct edges were coming into $p$, then we could remove $p$ from the vertex of $G_0$, which would contradict the minimality. Hence $p$ has a self loop, and no other edges coming to it. Since $G_0$ is connected, this implies that $G_0$ is the cycle, and in this case $E(G_0)=1<3=3\betti(G)$. Now, we can assume that $\deg(p) \geq 3$. Let $E$ denote the number of edges of $G_0$ with vertex set $V$. As each edge contributes by one to the degree of both of its endpoints (two if it is a loop), we have
    $$2E=\sum_{v \in V} \deg(v) \geq 3|V|.$$
    By Euler's formula, we get
    $$E=\betti(G_0)-1+|V| \leq \betti(G)-1+2E/3, $$
    implying $E \leq 3\betti(G)$.
\end{proof}

The following lemma shows we can extend a finite metric graph core.
\begin{lemma}\label{lem:graph_core_extension}
    Let $(X,d_X)$ be a compact geodesic space, $x \in X$, and $G$ be a core of $X$ which is a finite metric graph. Then there is a core $G'$ of $X$ which is a finite metric graph and contains $G$ and $x$.
\end{lemma}
\begin{proof}
    Assume $x \notin G$. Let $v$ be the point in $G$ to  $x$. Let $e$ denote the geodesic from $x$ to $v$. Let $V$ be the vertex set of a $1$-dimensional $CW$-complex structure on $G$ containing $v$. Let $G'=G \cup e$. By Lemma \ref{lem:core_extension}, $G'$ is a core. Let us show that $G'$ is a finite metric graph. $G'$ is a geodesic space, since it is retract of the geodesic space $X$ by Lemma \ref{lem:retraction}. Note that $G'-{V \cup x}$ is the disjoint union of $G-V$ and of $e-\{v, x \}$. By \cite[Exercise~3.2.17]{bbi01}, $G'$ is a metric graph.
\end{proof}

The following result shows how paths in $X$ are related to paths in its core.
\begin{lemma}\label{lem:core_path}
    Let $X$ be a compact geodesic space and $A$ be a core of $X$. Let $r: X \to A$ be the retraction given in Lemma \ref{lem:retraction}. Let $\beta:[0,1] \to X$ be a path whose endpoints are in $A$. For each $0 = t_0 \leq t_1 \dots \leq t_n = 1$, there exists $0 = s_0 \leq \dots \leq s_n = 1$ such that $\beta(s_i)=r(\beta(t_i))$ for all $i$.
\end{lemma}
\begin{proof}
    Let $s_i$ be the first time $\gamma$ enters into the core after $t_i$. Note that $s_0=0$, $s_n=1$, and $s_i \leq s_{i+1}$ for all $i=1,\dots,n-1$. Let us show that $\beta(s_i)=r(\beta(t_i))$. Assume $s_i<t_i$. Note that $\beta|_{(t_i,s_i)}$ is in $X-A$, hence it is contained in the disjoint union of open sets $\coprod_{a \in A} r^{-1}(a)-\{a\}$. This implies that there exists $a$ such that $\beta|_{(t_i,s_i)} \subseteq r^{-1}(a)$, hence $\beta|_{[t_i,s_i]} \subseteq r^{-1}(a)$, therefore $a=\beta(s_i)=r(\beta(s_i))=r(\beta(t_i))$.
\end{proof}

\bibliography{graph.bib}
\bibliographystyle{plain}

\end{document}